%% file: Symmetries_Outer_espace_Wn.tex
\documentclass[a4paper,oneside,11pt]{article}
 
\usepackage[a4paper]{geometry}
\usepackage{aeguill}                 
\usepackage{graphicx}
\usepackage{caption}
\usepackage{amsmath,amsfonts,amssymb,amsthm}
\usepackage{mathabx}
\usepackage{pstricks,comment} 
\usepackage{pst-plot}
\usepackage{pstricks-add}
\usepackage{multido}
\usepackage{url}
\usepackage{epstopdf,enumerate}
\usepackage{srcltx}
\usepackage[utf8]{inputenc} 
\usepackage[T1]{fontenc} 
\usepackage[english]{babel} 
\usepackage{hyperref} 
\usepackage[all]{xy} 
\usepackage{titlesec}
\usepackage{mathabx}
\usepackage{tikz}
\usepackage{fancyhdr}
\usepackage[us,12hr]{datetime}
\fancypagestyle{plain}{
\fancyhf{}
\rfoot{ {\ddmmyyyydate\today}}
\lfoot{\thepage}
}

    \newtheoremstyle{TheoremNum}
        {\topsep}{\topsep}              
        {\itshape}                      
        {}                              
        {\bfseries}                     
        {.}                             
        { }                             
        {\thmname{#1}\thmnote{ \bfseries #3}}

{\theoremstyle {definition} \newtheorem {defi} {Definition}[section]}
{\theoremstyle {plain}  \newtheorem {theo} [defi] {Theorem}}
{\theoremstyle {plain}  \newtheorem {coro} [defi] {Corollary}}
{\theoremstyle {plain} \newtheorem {prop} [defi] {Proposition}}
{\theoremstyle {plain} \newtheorem {lem}[defi] {Lemma}}
{\theoremstyle {plain} \newtheorem {rmq}[defi] {Remark}}
{\theoremstyle {plain} }
{\theoremstyle{TheoremNum} }
{\theoremstyle{TheoremNum} }
{\theoremstyle{TheoremNum} }

\newcommand{\Aut}{\mathrm{Aut}}
\newcommand{\Out}{\mathrm{Out}}
\newcommand{\Inn}{\mathrm{Inn}}

\newcommand{\PO}{\mathbb{P}\mathcal{O}}
\newcommand{\Stab}{\mathrm{Stab}}
\newcommand{\lk}{\mathrm{lk}}
\newcommand{\ZZ}{\mathbb{Z}}
\newcommand{\NN}{\mathbb{N}}

\newcommand{\dem}{\noindent{\bf Proof. }}

\title{The symmetries of the Outer space of a universal Coxeter group}
\author{Yassine Guerch}
\date{\today}
\begin{document}
\maketitle
\renewcommand*\labelenumi{(\theenumi)}

\begin{abstract}
This paper studies the geometric rigidity of the universal Coxeter
group of rank $n$, which is the free product $W_n$ of $n$ copies of
$\ZZ/2\ZZ$. We prove that for $n\geq 4$ the group of symmetries of the
spine of the Guirardel-Levitt outer space of $W_n$ is reduced to the
outer automorphism group $\Out(W_n)$.
\end{abstract}

\section{Introduction}
Let $n$ be an integer greater than $1$. Let $F=\ZZ/2\ZZ$ be a
cyclic group of order $2$ and $W_n=\bigast_n F$ be a universal Coxeter
group, which is a free product of $n$ copies of $F$. The geometric
study of automorphisms groups of free products is currently in strong
expansion, see for instance \cite{mcculloughmiller96,Guirardel,Piggott2012,HandelMosher14,gupta18,
Guirardelhorbez19,DahmaniLi19}.  This paper proves a major geometric rigidity result of the outer automorphism group $\Out(W_n)$ of $W_n$.

The study of $\Out(W_n)$ benefits from analogies with algebraic groups, with $\Out(F_N)$, the outer automorphism group of a free group of rank $N$, and with the mapping class group of a connected compact surface. As usual in geometric group theory, the understanding of the group $\Out(W_n)$ is related to the construction of  geometric spaces on which it acts nicely (properly or cocompactly). Such constructions appear in the study of $\Out(F_N)$, which involves the study of its action on the spine of the Outer space introduced by Culler and Vogtmann in \cite{Vogtmann1986}. Similarly, the study of the mapping class group of a connected compact surface involves the study of its action on the Teichmüller space and on the curve graph of the surface, while the study of algebraic groups implies the study of their actions on buildings. 

The spaces introduced in these cases are \emph{rigid geometric models} in the following sense: the symmetries of these spaces are induced by elements of the group itself. Indeed, for algebraic groups, Tits showed that, if the rank of a spherical building associated with a simple connected algebraic group is at least $2$, then the full group of simplicial automorphisms of the building is isomorphic to the algebraic group itself (\cite{Tits74}). In the context of a connected orientable compact surface of genus at least $3$, Royden proved that the group of isometries of the Teichmüller space with respect to the Teichmüller metric coincides with the extended mapping class group of the surface (\cite{royden71}). Moreover, Ivanov (\cite[Theorem~1]{Ivanov97}) showed that the group of simplicial automorphisms of the graph of curves is isomorphic to the extended mapping class group. In the context of $\Out(F_N)$, Bridson and Vogtmann proved that, if $N \geq 3$, the group of simplicial automorphisms of the spine of Outer space is isomorphic to $\Out(F_N)$ (\cite{bridson2001symmetries}).

In the case of $\Out(W_n)$, spaces on which $\Out(W_n)$ acts properly or cocompactly include the McCullough-Miller space \cite{mcculloughmiller96} or $\PO(W_n)$, the outer space of $W_n$ introduced by Guirardel and Levitt in \cite{Guirardel}. These two spaces are $\Out(W_n)$-equivariantly homotopy equivalent (see \cite[Theorem~8.5.]{mcculloughmiller96}). Moreover, it was proved by Piggott (\cite[Theorem~1.1]{Piggott2012}) that, for $n \geq 4$, the McCullough-Miller space is a rigid geometric model for $\Out(W_n)$: the group of simplicial automorphisms of the McMullough-Miller space is isomorphic to $\Out(W_n)$.

In this article, we study the action of $\Out(W_n)$ on a simplicial flag complex on which $\PO(W_n)$ retracts $\Out(W_n)$-equivariantly, called \emph{the spine of $\PO(W_n)$} and denoted by $K_n$. Vertices of $K_n$ are homothety classes of marked graphs of groups whose fundamental group is isomorphic to $W_n$. Two homothety classes $\mathcal{X}$ and $\mathcal{Y}$ of marked graphs of groups are adjacent in $K_n$ if they have representatives $X$ and $Y$ such that one can obtain $Y$ from $X$ by collapsing a forest in the underlying graph of $X$, or conversely. The group $\Out(W_n)$ naturally acts on $K_n$ by precomposing the marking. The aim of this article is to prove that $K_n$ is a rigid geometric model for $\Out(W_n)$ in the following sense. Here we denote by $\Aut(K_n)$ the group of simplicial automorphisms of $K_n$.

\begin{theo}\label{Rigidity Kn}
Let $n \geq 4$. The natural homomorphism $$\Out(W_n) \to \Aut(K_n)$$ is an isomorphism.
\end{theo}

This question is first motivated by the aforementioned examples, but also by algebraic results on $\Out(W_n)$. Indeed, for instance in the case of the mapping class group of a connected orientable compact surface of genus at least $3$, the fact that the curve complex is a rigid geometric model for the extended mapping class group is used by Ivanov in order to prove that any automorphism of the extended mapping class group is in fact a conjugation (see \cite[Theorem 2]{Ivanov97}). Similarly, the fact that the spine of Outer space is a rigid geometric model for $\Out(F_N)$ with $N \geq 3$ is related to the fact that any automorphism of $\Out(F_N)$ is a conjugation (\cite{bridson2000automorphisms}). As, for $n \geq 4$, any automorphism of $\Out(W_n)$ is a conjugation (see \cite[Théorème~1.1]{Guerch2020out}) and as the proof relies on the study of the action of $\Out(W_n)$ on $K_n$, it was natural to expect that $K_n$ is a rigid geometric model for $\Out(W_n)$. Even though the McCullough-Miller space and $\PO(W_n)$ are $\Out(W_n)$-equivariantly homotopy equivalent, the author does not know how to deduce the rigidity of $K_n$ out of the rigidity of the McCullough-Muller space. Indeed, there is no canonical graph isomorphism between $K_n$ and the McCullough-Miller space, and corresponding vertices in the McCullough-Miller space and in $K_n$ do not share the same properties of minimality. For instance the negative link of a $\{0\}$-star (see Sections~\ref{def Ln} and \ref{section rigidity K_n} for precise definitions) is nontrivial in $K_n$, whereas it is trivial in the McCullough-Miller space.

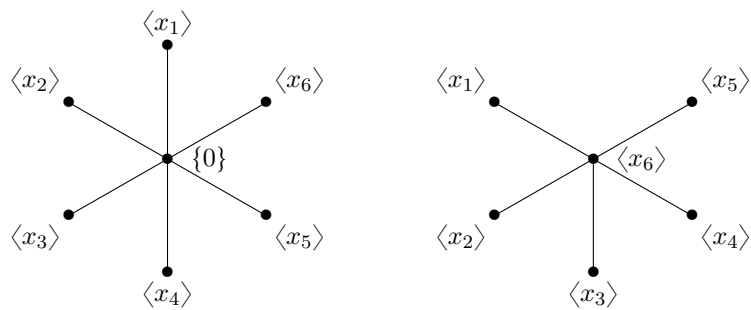
\begin{figure}[ht]
\centering
\captionsetup{justification=centering}
\begin{tikzpicture}[scale=1.5]
\draw (0:0) -- (90:1);
\draw (0:0) -- (150:1);
\draw (0:0) -- (210:1);
\draw (0:0) -- (270:1);
\draw (0:0) -- (330:1);
\draw (0:0) -- (30:1);
\draw (0:0) node {$\bullet$};
\draw (90:1) node {$\bullet$};
\draw (210:1) node {$\bullet$};
\draw (150:1) node {$\bullet$};
\draw (270:1) node {$\bullet$};
\draw (330:1) node {$\bullet$};
\draw (30:1) node {$\bullet$};
\draw (0:0) node[right, scale=0.9] {$\;\;\{0\}$};
\draw (90:1) node[above, scale=0.9] {$\left\langle x_1 \right\rangle$};
\draw (150:1) node[above left, scale=0.9] {$\left\langle x_2 \right\rangle$};
\draw (210:1) node[below left, scale=0.9] {$\left\langle x_3 \right\rangle$};
\draw (270:1) node[below, scale=0.9] {$\left\langle x_4 \right\rangle$};
\draw (330:1) node[below right, scale=0.9] {$\left\langle x_5 \right\rangle$};
\draw (30:1) node[above right, scale=0.9] {$\left\langle x_6 \right\rangle$};
\end{tikzpicture}
\hspace{1cm}
\begin{tikzpicture}[scale=1.5]
\draw (0:0) -- (270:1);
\draw (0:0) -- (150:1);
\draw (0:0) -- (210:1);
\draw (0:0) -- (330:1);
\draw (0:0) -- (30:1);
\draw (0:0) node {$\bullet$};
\draw (270:1) node {$\bullet$};
\draw (210:1) node {$\bullet$};
\draw (150:1) node {$\bullet$};
\draw (330:1) node {$\bullet$};
\draw (30:1) node {$\bullet$};
\draw (0:0) node[right, scale=0.9] {$\;\;\left\langle x_6 \right\rangle$};
\draw (270:1) node[below, scale=0.9] {$\left\langle x_3 \right\rangle$};
\draw (150:1) node[above left, scale=0.9] {$\left\langle x_1 \right\rangle$};
\draw (210:1) node[below left, scale=0.9] {$\left\langle x_2 \right\rangle$};

\draw (330:1) node[below right, scale=0.9] {$\left\langle x_4 \right\rangle$};
\draw (30:1) node[above right, scale=0.9] {$\left\langle x_5 \right\rangle$};
\end{tikzpicture}
\caption{A $\{0\}$-star (left) and an $F$-star (right).}\label{0 star F star}
\end{figure}

The proof of Theorem~\ref{Rigidity Kn} relies on the study of the action of $\Out(W_n)$ on a subgraph of $K_n$ called \emph{the graph of $\{0\}$-stars and $F$-stars} and denoted by $L_n$. Vertices of $L_n$ are $\{0\}$-stars and $F$-stars (see Section~\ref{def Ln} and Figure~\ref{0 star F star}). Two vertices of $L_n$ are adjacent if and only if they are adjacent in $K_n$. We first prove that $L_n$ is a rigid geometric model for $W_n$ (see Theorem~\ref{rigidity Ln}). This relies on studying systoles of $L_n$, that is, embedded cycles of minimal length. For this, we introduce (see Section~\ref{Section rigidity Ln}) a new complexity associated with an edge of $L_n$, and a relative complexity associated with pairs of $\{0\}$-stars. For $n=3$, the same study is not possible as the $\{0\}$-stars are no longer the vertices with minimal degree in $L_n$. We do not know whether Theorem~\ref{rigidity Ln} holds for $n=3$.

The rest of the proof consists in showing that there exists a homomorphism from $\Aut(K_n)$ to $\Aut(L_n)$ defined by restriction which turns out to be injective. We note that the characterization of the vertices of $L_n$ in $K_n$ is only based on the study of the possible decompositions of the link of the vertices of $K_n$. This differs from the proof of the similar result by Bridson and Vogtmann in the case of $\Out(F_n)$ since they used homological arguments in order to characterize some vertices of the spine of Outer space. Another major difference is that the strictly local rigidity properties of $L_n$ are much weaker than the ones of the spine of Outer space, and we need to explore the combinatorial balls of radius $4$ in $L_n$ in order to acquire a sufficient rigidity. Note that in the case of algebraic groups, Tits only needed to explore the combinatorial balls of radius $2$.

In Section~\ref{section FS}, we study the simplicial completion of $K_n$, denoted by $\overline{K}_n$. The simplicial complex $\overline{K}_n$ is also known as \emph{the free splitting complex of $W_n$} (see \cite{AramayonaSouto2011,HandelMosher13} and Section~\ref{section FS}). This complex has an analogue in the case of a free group of rank $N$, called \emph{the free splitting complex of $F_N$}. It was proved by Aramayona and Souto that the free splitting complex of $F_N$ is also a rigid geometric model for $\Out(F_N)$ when $N \geq 3$ (see \cite[Theorem~1]{AramayonaSouto2011}). In Section~\ref{section FS}, we prove the following theorem:

\begin{theo}\label{Rigidity FS}
Let $n \geq 4$. The natural homomorphism $$\Out(W_n) \to \Aut(\overline{K}_n)$$ is an isomorphism.
\end{theo}

Theorem~\ref{Rigidity FS} can be deduced from Theorem~\ref{Rigidity Kn} as follows. The spine $K_n$ has a natural embedding into $\overline{K}_n$. We first show that any automorphism of $\overline{K}_n$ preserves the image of $K_n$. This gives a homomorphism $\Aut(\overline{K}_n) \to \Aut(K_n)$ and the main point, using techniques of Scott-Swarup and Horbez-Wade, is to prove its injectivity. We then conclude using Theorem~\ref{Rigidity Kn}.

\bigskip

{\small{\bf Acknowledgments. } I warmly thank my advisors, Camille Horbez and Frédéric Paulin, for their precious advices and for carefully reading the different versions of this article.}

\section{Preliminaries}
\subsection{Background on the outer space of $W_n$}

Let $n$ be an integer greater than $1$. Let $F=\ZZ/2\ZZ$ be a cyclic group of order $2$ and $W_n= \bigast_n F$ be the universal Coxeter group of order $n$. We recall the definition of the outer space $\PO(W_n)$ introduced by Guirardel and Levitt in \cite{Guirardel}. A point in $\PO(W_n)$ is a homothety class of metric graph of groups $X$ whose fundamental group is $W_n$, equipped with a group isomorphism $\rho \colon W_n \to \pi_1(X)$ called a \emph{marking}, which satisfies~:
\begin{enumerate}
\item the underlying graph of $X$ is a finite tree~;
\item every edge group is trivial~;
\item there are exactly $n$ vertices whose associated group is isomorphic to $F$~;
\item all the other vertices have trivial associated group~;
\item if $v$ is a vertex whose associated group is trivial, then $\operatorname{deg}(v) \geq 3$.
\end{enumerate}

Two metric graphs of groups $(X,\rho)$ and $(X',\rho')$ are in the same homothety class if there exists a homothety $f \colon X \to X'$ (meaning an application multiplying all edge lengths by the same scalar) and such that $f_{\ast} \circ \rho=\rho'$. We denote by $[X,\rho]$ the homothety class of such a metric graph of groups $(X,\rho)$. If the marking is implicit, we denote by $\mathcal{X}$ the homothety class. The group $\Aut(W_n)$ acts by precomposing the marking. As, for any $\alpha \in \Inn(W_n)$, and for any $\mathcal{X} \in \PO(W_n)$, we have $\alpha(\mathcal{X})=\mathcal{X}$, the action of $\Aut(W_n)$ induces an action of $\Out(W_n)$.

\bigskip

The set $\PO(W_n)$ is equipped with a topology which we recall now. Suppose that \mbox{$[X,\rho] \in \PO(W_n)$} and let $(X,\rho)$ be the representative of $[X,\rho]$ such that the sum of the edge lengths is equal to $1$. To $(X,\rho)$ we associate a simplex by varying the lengths of the edges, so that the sum of the edge lengths is still equal to $1$. A homothety class $[X',\rho'] \in \PO(W_n)$ defines a codimension $1$ face of the simplex associated with $(X,\rho)$ if we can obtain $(X',\rho')$ from $(X,\rho)$ by contracting an edge of the underlying graph of $X$. The \emph{weak topology} is then defined in the following way: a set is open if and only if its intersection with every open simplex is open. 

\bigskip

We now recall the definition of a deformation retract of $\PO(W_n)$ known as \emph{the spine of $\PO(W_n)$} and denoted by $K_n$. It is a flag complex whose vertices are the open simplices associated with each homothety class $[X,\rho] \in \PO(W_n)$. Two vertices corresponding to two homothety classes $[X,\rho]$ and $[X',\rho']$ are adjacent if $[X,\rho]$ defines a face of the simplex associated with $[X',\rho']$ and conversely. There is an embedding $F \colon K_n \hookrightarrow \PO(W_n)$ whose image is the barycentric spine of $\PO(W_n)$. We will from then on identify $K_n$ with $F(K_n)$.

We now give a description of the stabilizer of a point in $\overline{K}_n$ due to Levitt. If $\mathcal{X} \in V\overline{K}_n$, we denote by $\Stab(\mathcal{X})$ the stabilizer of $\mathcal{X}$ under the action of $\Out(W_n)$. Let $X$ be a representative of $\mathcal{X}$. We denote by $\Stab^0(\mathcal{X})$ the subgroup of $\Stab(\mathcal{X})$ made of all elements $F \in \Out(W_n)$ such that the automorphism induced by $F$ on $X$ is the identity. We write the next proposition in a more general context where the nontrivial vertex groups are not necessarly isomorphic to $F$ (see Section~\ref{section FS}).

\begin{prop}\cite[Proposition 4.2]{levitt2005}\label{Levitt stab}
Let $n \geq 4$ and $\mathcal{X} \in V\overline{K}_n$. Let $X$ be a representative of $\mathcal{X}$ and let $v_1,\ldots,v_k$ be the vertices of $X$ with nontrivial associated groups. For $i \in \{1,\ldots,k\}$, let $G_i$ be the group associated with $v_i$. Then $\Stab^0(\mathcal{X})$ is isomorphic to 
$$ \prod\limits_{i=1}^k G_i^{\operatorname{deg}(v_i)-1} \rtimes \Aut(G_i)~,$$
where $\Aut(G_i)$ acts on $G_i^{\operatorname{deg}(v_i)-1}$ diagonally.
\end{prop}

\subsection{The graph of $\{0\}$-stars and $F$-stars.}\label{def Ln}

In order to prove Theorem~\ref{Rigidity Kn}, we introduce a graph included in the spine $K_n$ called \emph{the graph of $\{0\}$-stars and $F$-stars}.

\begin{defi}
\noindent{$(1)$ } A \emph{$\{0\}$-star} is the equivalence class in $K_n$ of a metric graph of groups whose underlying graph has $n+1$ vertices and $n$ leaves. 

\medskip

\noindent{$(2)$ } A \emph{$F$-star} is the equivalence class in $K_n$ of a metric graph of groups whose underlying graph has $n$ vertices and $n-1$ leaves. 

\medskip

\noindent{$(3)$ } The \emph{graph of $\{0\}$-stars and $F$-stars}, denoted by $L_n$, is the full subgraph of $K_n$ whose vertices are exactly the $\{0\}$-stars and the $F$-stars. There is an edge between two vertices of $L_n$ if and only if there is an edge between the corresponding vertices in $K_n$.
\end{defi}

As $\Aut(W_n)$ acts on $K_n$ by precomposition of the action, the graph $L_n$ is invariant by $\Out(W_n)$.

Since any two $\{0\}$-stars are at distance at least $2$ in $K_n$, the neighbors of a $\{0\}$-star in $L_n$ are $F$-stars. Conversely, since any two $F$-stars are at distance at least $2$ in $K_n$, the neighbors of an $F$-star in $L_n$ are $\{0\}$-stars. The number of neighbors in $L_n$ of a $\{0\}$-star is equal to $n$. They correspond to collapsing exactly one edge of the underlying graph. The number of neighbors in $L_n$ of an $F$-star is equal to $2^{n-2}$. They correspond to blowing-up the central vertex of the underlying graph while applying a partial conjugation by the generator of the preimage by the marking of the group associated with the center. As $\Aut(W_n)$ acts transitively on the set of free bases of $W_n$, we see that $\Aut(W_n)$ acts transitively on the set of $\{0\}$-stars. Thus, as partial conjugations and permutations generate $\Aut(W_n)$ by \cite[Theorem B]{muhlherr1997}, it follows that the graph $L_n$ is connected.

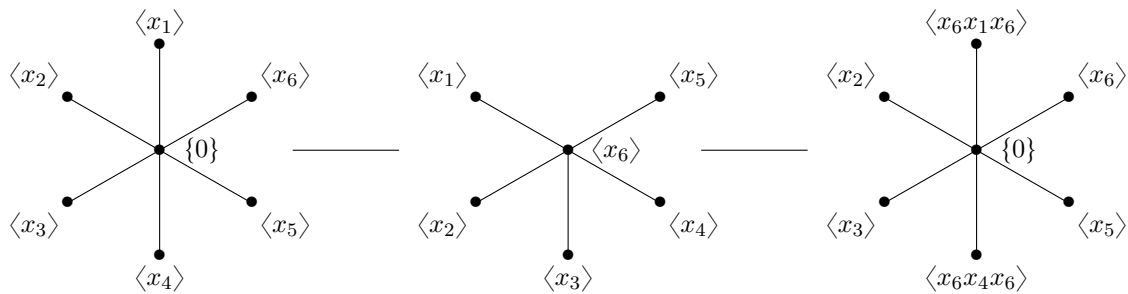
\begin{figure}[ht]
\centering
\captionsetup{justification=centering}
\hspace{-1cm}
\begin{tikzpicture}[scale=1.4]
\draw (0:0) -- (90:1);
\draw (0:0) -- (150:1);
\draw (0:0) -- (210:1);
\draw (0:0) -- (270:1);
\draw (0:0) -- (330:1);
\draw (0:0) -- (30:1);
\draw (0:0) node {$\bullet$};
\draw (90:1) node {$\bullet$};
\draw (210:1) node {$\bullet$};
\draw (150:1) node {$\bullet$};
\draw (270:1) node {$\bullet$};
\draw (330:1) node {$\bullet$};
\draw (30:1) node {$\bullet$};
\draw (0:0) node[right, scale=0.9] {$\;\;\{0\}$};
\draw (90:1) node[above, scale=0.9] {$\left\langle x_1 \right\rangle$};
\draw (150:1) node[above left, scale=0.9] {$\left\langle x_2 \right\rangle$};
\draw (210:1) node[below left, scale=0.9] {$\left\langle x_3 \right\rangle$};
\draw (270:1) node[below, scale=0.9] {$\left\langle x_4 \right\rangle$};
\draw (330:1) node[below right, scale=0.9] {$\left\langle x_5 \right\rangle$};
\draw (30:1) node[above right, scale=0.9] {$\left\langle x_6 \right\rangle$};

\draw (0:1.25) -- (0:2.25);
\end{tikzpicture}
\begin{tikzpicture}[scale=1.4]

\draw (0:0) -- (270:1);
\draw (0:0) -- (150:1);
\draw (0:0) -- (210:1);
\draw (0:0) -- (330:1);
\draw (0:0) -- (30:1);
\draw (0:0) node {$\bullet$};
\draw (270:1) node {$\bullet$};
\draw (210:1) node {$\bullet$};
\draw (150:1) node {$\bullet$};
\draw (330:1) node {$\bullet$};
\draw (30:1) node {$\bullet$};
\draw (0:0) node[right, scale=0.9] {$\;\;\left\langle x_6 \right\rangle$};
\draw (270:1) node[below, scale=0.9] {$\left\langle x_3 \right\rangle$};
\draw (150:1) node[above left, scale=0.9] {$\left\langle x_1 \right\rangle$};
\draw (210:1) node[below left, scale=0.9] {$\left\langle x_2 \right\rangle$};

\draw (330:1) node[below right, scale=0.9] {$\left\langle x_4 \right\rangle$};
\draw (30:1) node[above right, scale=0.9] {$\left\langle x_5 \right\rangle$};
\draw (0:1.25) -- (0:2.25);
\end{tikzpicture}
\begin{tikzpicture}[scale=1.4]

\draw (0:0) -- (90:1);
\draw (0:0) -- (150:1);
\draw (0:0) -- (210:1);
\draw (0:0) -- (270:1);
\draw (0:0) -- (330:1);
\draw (0:0) -- (30:1);
\draw (0:0) node {$\bullet$};
\draw (90:1) node {$\bullet$};
\draw (210:1) node {$\bullet$};
\draw (150:1) node {$\bullet$};
\draw (270:1) node {$\bullet$};
\draw (330:1) node {$\bullet$};
\draw (30:1) node {$\bullet$};
\draw (0:0) node[right, scale=0.9] {$\;\;\{0\}$};
\draw (90:1) node[above, scale=0.9] {$\left\langle x_6x_1x_6 \right\rangle$};
\draw (150:1) node[above left, scale=0.9] {$\left\langle x_2 \right\rangle$};
\draw (210:1) node[below left, scale=0.9] {$\left\langle x_3 \right\rangle$};
\draw (270:1) node[below, scale=0.9] {$\left\langle x_6x_4x_6 \right\rangle$};
\draw (330:1) node[below right, scale=0.9] {$\left\langle x_5 \right\rangle$};
\draw (30:1) node[above right, scale=0.9] {$\left\langle x_6 \right\rangle$};
\end{tikzpicture}
\caption{Examples of two neighbours of an $F$-star in $L_6$.}\label{F stars in L6}
\end{figure}

\section{Rigidity of the graph of $\{0\}$-stars and $F$-stars}\label{Section rigidity Ln}

In this section, we prove the following theorem.

\begin{theo}\label{rigidity Ln}
Let $n \geq 4$. Let $f$ be an automorphism of $L_n$ preserving $O_n$ and $F_n$. Then $f$ is induced by the action of a unique element $\gamma$ of $\Out(W_n)$.
\end{theo}

For $n \geq 5$, any $F$-star has $2^{n-2}$ neighbours in $L_n$ and any $\{0\}$-star has $n$ neighbours in $L_n$. As $2^{n-2}>n$ precisely when $n \geq 5$, we see that every automorphism of $L_n$ preserves the set of $\{0\}$-stars and the set of $F$-stars. We thus have the following corollary.

\begin{coro}\label{corollary rigidity Ln}
Let $n \geq 5$. The natural homomorphism $$\Out(W_n) \to \Aut(L_n)$$ is an isomorphism.
\hfill\qedsymbol
\end{coro}

Before proving Theorem~\ref{rigidity Ln}, we first prove a lemma which characterises the number of paths in a ball of radius $4$ centered at a $\{0\}$-star. 

\bigskip

Let $\mathcal{X}$ be a $\{0\}$-star, and $(X,\rho)$ a representative of $\mathcal{X}$. Let $v_1,\ldots,v_n$ be the $n$ leaves of the underlying graph of $X$. For $i \in \{1,\ldots,n\}$, let $x_i$ be the preimage by $\rho$ of the generator of the group associated with $v_i$, and let $\mathcal{Y}_i$ be the $F$-star adjacent to $\mathcal{X}$ such that a representative of $\mathcal{Y}_i$ is obtained from $X$ by contracting the edge adjacent to $v_i$. For distinct $i,j \in \{1,\ldots,n\}$, let $\sigma_{j,i} \colon W_n \to W_n$ be the automorphism sending $x_j$ to $x_ix_jx_i$ and, for $k \neq j$, fixing $x_k$. In this context we will call $x_i$ the \emph{twistor of $\sigma_{j,i}$}. For distinct $i,j \in \{1,\ldots,n\}$, let $(i\;j)$ be the automorphism of $W_n$ switching $x_i$ and $x_j$ and, for $k \neq i,j$, fixing $x_k$. A theorem of Mühlherr (c.f. \cite[Theorem B]{muhlherr1997}) implies that $\{\sigma_{i,j} \, |i \neq j \} \cup \{(i\;j) \,| i \neq j\}$ is a generating set of $\Aut(W_n)$. Note that, for every integers $i,j,k,\ell$, there exist $p,q$ such that $(i\;j)\sigma_{k,\ell}(i\;j)=\sigma_{p,q}$.

We now fix $i \in \{1,\ldots,n\}$. Let $\mathcal{X}'$ be a $\{0\}$-star adjacent to $\mathcal{Y}_i$ and distinct from $\mathcal{X}$. Let $(X',\rho')$ be a representative of $\mathcal{X}'$. Let $w_1,\ldots,w_n$ be the leaves of the underlying graph of $X'$, and, for $j \in \{1,\ldots,n\}$, let $y_j$ be the preimage by $\rho'$ of the generator of the group associated with $w_j$. Up to composition by an inner automorphism and reordering, either $y_j=x_j$ or $y_j=x_iy_jx_i$ (see Figure~\ref{F stars in L6} with $i=n=6$). Thus, there exist $k \in \{1,\ldots,n-1\}$ and $i_1,\ldots,i_k \in \{1,\ldots, \widehat{i},\ldots,n\}$ pairwise distinct such that, for all $j \in \{1,\ldots,n\}$, 

$$\Big(\prod\limits_{l=1}^k \sigma_{i_l,i} \Big)(x_j)=y_j.$$

\begin{figure}
\centering
\captionsetup{justification=centering}
\begin{tikzpicture}[scale=1.5]
\draw (0:0) -- (270:1);
\draw (0:0) -- (150:1);
\draw (0:0) -- (210:1);
\draw (0:0) -- (330:1);
\draw (0:0) -- (30:1);
\draw (0:0) node {$\bullet$};
\draw (270:1) node {$\bullet$};
\draw (210:1) node {$\bullet$};
\draw (150:1) node {$\bullet$};
\draw (330:1) node {$\bullet$};
\draw (30:1) node {$\bullet$};
\draw (0:0) node[right, scale=0.9] {$\;\;\{0\}$};
\draw (270:1) node[below, scale=0.9] {$\left\langle x_3 \right\rangle$};
\draw (150:1) node[above left, scale=0.9] {$\left\langle x_3x_1x_3 \right\rangle$};
\draw (210:1) node[below left, scale=0.9] {$\left\langle x_3x_2x_3 \right\rangle$};

\draw (330:1) node[below right, scale=0.9] {$\left\langle x_4 \right\rangle$};
\draw (30:1) node[above right, scale=0.9] {$\left\langle x_5 \right\rangle$};
\end{tikzpicture}
\hspace{1cm}
\begin{tikzpicture}[scale=1.5]
\draw (0:0) -- (270:1);
\draw (0:0) -- (150:1);
\draw (0:0) -- (210:1);
\draw (0:0) -- (330:1);
\draw (0:0) -- (30:1);
\draw (0:0) node {$\bullet$};
\draw (270:1) node {$\bullet$};
\draw (210:1) node {$\bullet$};
\draw (150:1) node {$\bullet$};
\draw (330:1) node {$\bullet$};
\draw (30:1) node {$\bullet$};
\draw (0:0) node[right, scale=0.9] {$\;\;\{0\}$};
\draw (270:1) node[below, scale=0.9] {$\left\langle x_3 \right\rangle$};
\draw (150:1) node[above left, scale=0.9] {$\left\langle x_1 \right\rangle$};
\draw (210:1) node[below left, scale=0.9] {$\left\langle x_2 \right\rangle$};

\draw (330:1) node[below right, scale=0.9] {$\left\langle x_3x_4x_3 \right\rangle$};
\draw (30:1) node[above right, scale=0.9] {$\left\langle x_3x_5x_3 \right\rangle$};
\end{tikzpicture}
\caption{Two representatives of the same homothety class $\mathcal{X}'$ realizing $k_{\mathcal{X},i}(\mathcal{X}')$.}\label{example two different k(X')}
\end{figure}
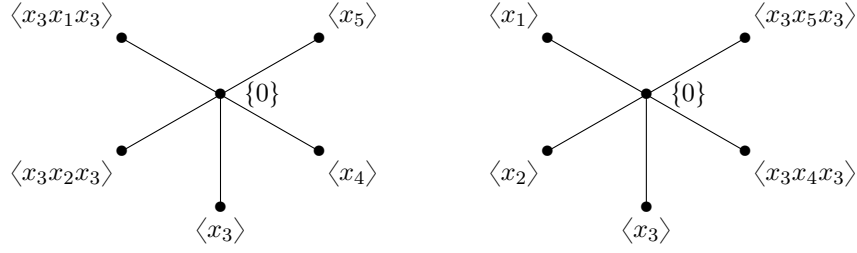

\noindent Let $\Inn^{\#}(W_n)=\left\langle\Inn(W_n),\{\sigma_{i,j}\;|\; i \neq j\}\right\rangle$. We define the \emph{first term complexity of $\mathcal{X}'$} by 
$$
k_{\mathcal{X},i}(\mathcal{X}')=\min \left\{
k \;\Bigg|\; 
\begin{array}{c} \exists i_1,\ldots,i_k \in \{1,\ldots, \hat{i},\ldots,n\}, I \in \Inn^{\#}(W_n) \text{ such that } \\
\forall j \in \{1,\ldots,n\}, \; I \circ \Big(\prod\limits_{l \in \{1,\ldots,k\}} \sigma_{i_l,i}\Big)(x_j)=y_j
\end{array}
\right\}.
$$
This definition does not depend on the choice of a representative of $\mathcal{X}'$. Note that the sequence $i_1,\ldots,i_k$ realizing the minimum is not necessarily unique (see Figure~\ref{example two different k(X')} with $n=5$ and $i=3$). However, if $k_{\mathcal{X},i}(\mathcal{X}') \neq n-k_{\mathcal{X},i}(\mathcal{X}')-1$, such a sequence is unique.

\bigskip

We now define a notion of relative complexity in $O_n$, the set of $\{0\}$-stars in $L_n$. Let $\mathcal{Z}$ be a $\{0\}$-star in $L_n$ distinct from $\mathcal{X}$ and let $(Z,\psi)$ be a representative of $\mathcal{Z}$. Let $w_1,\ldots,w_n$ be the leaves of the underlying graph of $Z$, and, for $j \in \{1,\ldots,n\}$, let $z_j$ be the preimage by $\psi$ of the generator of the group associated with $w_j$. As $\{\sigma_{i,j} \, |i \neq j \} \cup \{(i\;j) \,| i \neq j\}$ is a generating set of $\Aut(W_n)$ (c.f. \cite[Theorem B]{muhlherr1997}), we see that, up to composition by an inner automorphism and reordering, 
$$ \begin{array}{c}
\exists k \in \NN, \; \exists(i_1,j_1),\ldots,(i_k,j_k) \in \{1,\ldots,n\}^2-\{(x,x)\;|\; x \in \{1,\ldots,n\}\}, \\
\forall j \in \{1,\ldots,n\}, \; \Big(\prod\limits_{m=1}^k\sigma_{i_m,j_m}\Big)(x_j)=z_j.
\end{array}
$$

We now define the \emph{second term complexity of $\mathcal{Z}$} by 

$$ \ell_{\mathcal{X}}(\mathcal{Z})=\min
\left\{
\ell \;\Bigg|\;
\begin{array}{c}
\exists a_1,\ldots,a_{\ell} \in \{1,\ldots,n\} \text{ pairwise distinct, } \\
\exists I \in \Inn^{\#}(W_n) \text{ such that, } \; \exists k \in \NN \\
\exists (i_1,j_1),\ldots,(i_k,j_k) \in \{1,\ldots,n\} \times \{a_1,\ldots,a_{\ell}\}-\{(x,x)\;|\; x \in \{a_1,\ldots,a_{\ell}\}\}, \\
\forall j \in \{1,\ldots,n\}, \; I \circ \Big(\prod\limits_{m=1}^k\sigma_{i_m,j_m}\Big)(x_j)=z_j
\end{array}
\right\}.
$$

The intuition behind the second term complexity is the following. We want to count the minimal number $N$ of elements of $\{x_1,\ldots,x_n\}$ such that, for all $j \in \{1,\ldots,n\}$, the generator $z_j$ can be obtained from $x_j$ using partial conjugations by at most $N$ twistors. Note that, in the definition of $\ell_{\mathcal{X}}(\mathcal{Z})$, if $\Big(\prod\limits_{m=1}^k\sigma_{i_m,j_m}\Big)(x_j)=z_j$,
and if $m_1,m_2 \in \{1,\ldots,k\}$ are distinct, we do not require that $j_{m_1} \neq j_{m_2}$, so that the same twistor can appear in distinct partial conjugations. Note also that $\ell_{\mathcal{X}}(\mathcal{Z})$ does not depend on the choice of a representative of $\mathcal{Z}$. 

\begin{lem}\label{properties of second term complexity}
Let $\mathcal{X}$, $(X,\rho)$, $v_1,\ldots,v_n$ and $(\mathcal{Y}_i)_{i=1,\ldots,n}$ be as above. 

\noindent{$(1)$ } Fix $i \in \{1,\ldots,n\}$ and let $\mathcal{X}'$ be a $\{0\}$-star adjacent to $\mathcal{Y}_i$ and distinct from $\mathcal{X}$. Then $\ell_{\mathcal{X}}(\mathcal{X}')=1$ and a set $\{a_1,\ldots,a_{\ell}\}$ realizing the minimum defining $\ell_{\mathcal{X}}(\mathcal{X}')$ is $\{i\}$.

\medskip

\noindent{$(2)$ } Let $B(\mathcal{X},r)$ be the closed ball in $L_n$ of radius $r$ centered at $\mathcal{X}$. Let $\mathcal{Z} \in  B(\mathcal{X},4) \cap O_n$. Then $\ell_{\mathcal{X}}(\mathcal{Z}) \leq 2$. Moreover, the set realizing the minimum defining $\ell_{\mathcal{X}}(\mathcal{Z})$ is unique.
\end{lem}

\dem
Let $(X',\rho')$ be a representative of $\mathcal{X}'$. Let $y_1,\ldots,y_n$ be the preimage by $\rho'$ of the generators of the nontrivial vertex groups of $X'$. Then, up to composing by an inner automorphism and reordering, for all $j \in \{1,\ldots,n\}$, either $y_j=x_j$ or $y_j=x_ix_jx_i$. Thus, for all $j \in \{1,\ldots,n\}$, the only twistor that we need in order to obtain $y_j$ from $x_j$ using partial conjugations is $x_i$. Since $\mathcal{X}' \neq \mathcal{X}$, it follows that $\ell_{\mathcal{X}}(\mathcal{X}')=1$ and that a set realizing the minimum defining it is $\{i\}$. 

\bigskip

For the second assertion, let $Z$ be a representative of $\mathcal{Z}$, and let $z_1,\ldots,z_n$ be the preimage by the marking of the generators of the vertex groups. Then, there exist $j,k \in\{1,\ldots,n\}$ such that, for all $m \in \{1,\ldots,n\}$, one of the following holds:

\noindent{$(1)$ } $z_m=x_m$,

\smallskip

\noindent{$(2)$ } $z_m=x_jx_mx_j$,

\smallskip

\noindent{$(3)$ } $z_m=x_kx_mx_k$,

\smallskip

\noindent{$(4)$ } $z_m=x_kx_jx_kx_mx_kx_jx_k$,

\smallskip

\noindent{$(5)$ } $z_m=x_kx_jx_mx_jx_k$.

\smallskip

Thus, for all $m \in \{1,\ldots,n\}$, as we only need $x_j$ and $x_k$ as twistors to obtain $z_m$ from $x_m$, we see that $\ell_{\mathcal{X}}(\mathcal{Z}) \leq 2$.

Moreover, the twistors $x_j$ and $x_k$ are the unique elements of $\{x_1,\ldots,x_n\}$ such that, for all $i \in \{1,\ldots,n\}$, the generator $z_i$ is obtained from $x_i$ by partial conjugations using $x_j$ and $x_k$ as twistors. Thus, for all $\mathcal{Z} \in B(\mathcal{X},4) \cap O_n$, the set $\{a_1,\ldots,a_l\}$ realizing the minimum defining $\ell_{\mathcal{X}}(\mathcal{Z})$ is unique.
\hfill\qedsymbol

\bigskip

We isolate here a technical argument that will appear frequently in the proof of Lemma~\ref{number of path in Ln}.

\begin{lem}\label{second term complexity and 3 letters}
Let $\mathcal{X}$, $(X,\rho)$, $v_1,\ldots,v_n$ and $(\mathcal{Y}_i)_{i=1,\ldots,n}$ be as above.

Fix $i \in \{1,\ldots,n\}$ and let $\mathcal{X}'$ be a $\{0\}$-star adjacent to $\mathcal{Y}_i$ and distinct from $\mathcal{X}$. 
Let $k,\ell \in \{1,\ldots,n\}-\{i\}$ be distinct. Let $\mathcal{X}_k^{(2)}$ be a $\{0\}$-star such that:

\noindent{$\bullet$} $d(\mathcal{X}',\mathcal{X}_k^{(2)})=2$,

\noindent{$\bullet$} $\ell_{\mathcal{X}}(\mathcal{X}_k^{(2)})=2$ and a set realizing the minimum defining it is $\{i,k\}$. 

Let $\mathcal{X}_k^{(3)}$ be a $\{0\}$-star at distance $2$ of $\mathcal{X}_k^{(2)}$ and such that any set realizing $\ell_{\mathcal{X}}(\mathcal{X}_k^{(3)})$ contains $\ell$. Then $\ell_{\mathcal{X}}(\mathcal{X}_k^{(3)}) \geq 3$.
\end{lem}

\dem Let $(X',\rho')$ be a representative of $\mathcal{X}'$. Let $w_1,\ldots,w_n$ be the leaves of the underlying graph of $X'$, and, for $m \in \{1,\ldots,n\}$, let $y_m$ be the preimage by $\rho'$ of the generator of the group associated with $w_m$. For $j \in \{2,3\}$, let $(X_k^{(j)},\psi^{(j)})$ be a representative of $\mathcal{X}_k^{(j)}$, let $w_1^{(j)},\ldots,w_n^{(j)}$ be the $n$ leaves of the underlying graph of $X_k^{(j)}$ and, for $m \in \{1,\ldots,n\}$, let $y_m^{(j)}$ be the preimage by $\psi^{(j)}$ of the generator of the group associated with $w_m^{(j)}$. Note that, up to composition by an inner automorphism and reordering, for all $m \in \{1,\ldots,n\}$, $$y_m^{(2)}=x_i^{\gamma_m}x_k^{\beta_m}x_i^{\alpha_m}x_mx_i^{\alpha_m}x_k^{\beta_m}x_i^{\gamma_m}, \quad \alpha_m,\beta_m,\gamma_m \in \{0,1\}.$$

Note also that $\gamma_m=1$ precisely when $y_k=x_ix_kx_i$ and $\beta_m=1$. Thus, for all $m \in \{1,\ldots,n\}$, the element $y_m^{(2)}$ is obtained from $x_m$ using partial conjugations with twistors $x_i$ and $x_k$. Moreover, as $k \neq i$, and as $\mathcal{X}' \neq \mathcal{X}$, there exists $n_1 \in \{1,\ldots,n\}$ such that $\alpha_{n_1} \neq 0$ or $\gamma_{n_1} \neq 0$. Since $\mathcal{X}_k^{(2)} \neq \mathcal{X}'$, there exists $n_2$ such that $\beta_{n_2} \neq 0$.

As $\ell$ is contained in any set realizing the minimum defining $\ell_{\mathcal{X}}(\mathcal{X}_k^{(3)})$, there exists $p \in \NN$ and $m_1 \in \{1,\ldots,n\}$ such that 

$$\Big(\prod\limits_{m=1}^p\sigma_{i_m,j_m}\Big)(x_{m_{1}})=y_{m_{1}}^{(3)},$$

\noindent and there exists $m$ such that $j_m=\ell$.

\medskip

\noindent{\bf Claim. } The elements $x_k$ and $x_i$ are twistors of any set realizing the minimum defining $\ell_{\mathcal{X}}(\mathcal{X}_k^{(3)})$. 

\medskip

\dem As $\ell$ is contained in any set realizing $\ell_{\mathcal{X}}(\mathcal{X}_k^{(3)})$, a representative of $\mathcal{X}_k^{(3)}$ is obtained from $X_k^{(2)}$ by contracting the edge adjacent to $w_{\ell}^{(2)}$ and then blowing-up an edge at the central vertex. We then distinguish different cases according to the value of $y_{\ell}^{(2)}$. 

If $y_{\ell}^{(2)}=x_{\ell}$, then for all $m \in \{1,\ldots,n\}$, $$y_m^{(3)}=x_{\ell}^{\delta_m}x_i^{\gamma_m}x_k^{\beta_m}x_i^{\alpha_m}x_mx_i^{\alpha_m}x_k^{\beta_m}x_i^{\gamma_m}x_{\ell}^{\delta_m}, \quad \alpha_m,\beta_m,\gamma_m,\delta_m \in \{0,1\}.$$ 
Since $n_1$ and $n_2$ are such that $\alpha_{n_1} \neq 0$ or $\gamma_{n_1} \neq 0$, and $\beta_{n_2} \neq 0$, the claim follows. 

If $y_{\ell}^{(2)}=x_ix_{\ell}x_i$, then for all $m \in \{1,\ldots,n\}$, we have $$y_m^{(3)}=x_i^{\delta_m}x_{\ell}^{\delta_m}x_i^{\gamma_m+\delta_m}x_k^{\beta_m}x_i^{\alpha_m}x_mx_i^{\alpha_m}x_k^{\beta_m}x_i^{\gamma_m+\delta_m}x_{\ell}^{\delta_m}x_i^{\delta_m}, \quad \alpha_m,\beta_m,\gamma_m,\delta_m \in \{0,1\}.$$
Since $n_2$ is such that $\beta_{n_2} \neq 0$ and $n_1$ is such that $\gamma_{n_1}+\delta_{n_1} \neq 0$, or $\alpha_{n_1} \neq 0$ or $\delta_{n_1} \neq 0$, the claim follows. 

If $y_{\ell}^{(2)}=x_kx_{\ell}x_k$, then for all $m \in \{1,\ldots,n\}$, we have $$y_m^{(3)}=x_k^{\delta_m}x_{\ell}^{\delta_m}x_k^{\delta_m}x_i^{\gamma_m}x_k^{\beta_m}x_i^{\alpha_m}x_mx_i^{\alpha_m}x_k^{\beta_m}x_i^{\gamma_m}x_k^{\delta_m}x_{\ell}^{\delta_m}x_k^{\delta_m}, \quad \alpha_m,\beta_m,\gamma_m,\delta_m \in \{0,1\}.$$
Since $n_1$ and $n_2$ are such that $\alpha_{n_1} \neq 0$ or $\gamma_{n_1} \neq 0$, and $\beta_{n_2} \neq 0$, the result follows. 

Finally, if $y_{\ell}^{(2)}=x_kx_ix_{\ell}x_ix_k$, or if $y_{\ell}^{(2)}=x_ix_kx_{\ell}x_kx_i$, or if  $y_{\ell}^{(2)}=x_ix_kx_ix_{\ell}x_ix_kx_i$, then $y_{m_1}$ is obtained from $x_{m_1}$ using $x_i$, $x_k$ and $x_{\ell}$ as twistors. This concludes the proof of the claim.
\hfill\qedsymbol

\medskip

Thus, $i$, $k$ and $\ell$ are contained in any set realizing the minimum defining $\ell_{\mathcal{X}}(\mathcal{X}_k^{(3)})$, and this implies that $\ell_{\mathcal{X}}(\mathcal{X}_k^{(3)}) \geq 3$.
\hfill\qedsymbol

\bigskip

We are now ready to prove a lemma concerning the number of embedded paths in $B(\mathcal{X},4)$.

\begin{lem}\label{number of path in Ln}
Let $\mathcal{X}$, $(X,\rho)$, $v_1,\ldots,v_n$ and $(\mathcal{Y}_i)_{i=1,\ldots,n}$ be as above.

Fix $i \in \{1,\ldots,n\}$ and let $\mathcal{X}'$ be a $\{0\}$-star adjacent to $\mathcal{Y}_i$ and distinct from $\mathcal{X}$. Let $(X',\rho')$ be a representative of $\mathcal{X'}$ and let $\overline{X}'$ be the underlying graph of $X'$. Let $$\{i_1,\ldots,i_{k_{\mathcal{X},i}(\mathcal{X}')}\} \subseteq \{1,\ldots,\hat{i},\ldots,n\}$$ be a set realizing the minimum defining $k_{\mathcal{X},i}(\mathcal{X}')$, and $j \in \{1,\ldots,\hat{i},\ldots,n\}$. Let $x_1',\ldots,x_n'$ be the preimages by $\rho'$ of the generators of the nontrivial vertex groups. Up to reordering, suppose that, for all $k \in \{1,\ldots,n\}$, $x_k'$ is obtained from $x_k$ by a conjugation.
\begin{enumerate}
\item If $x_j'=x_j$, the number of distinct injective edge paths in $B(\mathcal{X},4)-\{\mathcal{X}\}$ of length at most $5$ between $\mathcal{X}'$ and $\mathcal{Y}_j$ is equal to $2^{n-k_{\mathcal{X},i}(\mathcal{X}')-2}-1$.

\item If $x_j'=x_ix_jx_i$, the number of distinct injective edge paths in $B(\mathcal{X},4)-\{\mathcal{X}\}$ of length at most $5$ between $\mathcal{X}'$ and $\mathcal{Y}_j$ is equal to $2^{k_{\mathcal{X},i}(\mathcal{X}')-1}-1$.

\item Let $\mathcal{Z}$ be a $\{0\}$-star distinct from $\mathcal{X}$ and adjacent to $\mathcal{Y}_j$ and such that $k_{\mathcal{X},j}(\mathcal{Z})=1$. Let $\{t\}$ be a set realizing the minimum defining $k_{\mathcal{X},j}(\mathcal{Z})$. Suppose that $x_j'=x_j$. 

If $t \in \{i_1,\ldots,i_{k_{\mathcal{X},i}(\mathcal{X}')}\}$, then there is no path between $\mathcal{X}'$ and $\mathcal{Z}$ of length at most $4$ in $B(\mathcal{X},4)-\{\mathcal{X}\}$. If $t \notin \{i_1,\ldots,i_{k_{\mathcal{X},i}(\mathcal{X}')}\}$, then there is at least one path between $\mathcal{X}'$ and $\mathcal{Z}$ of length at most $4$ in $B(\mathcal{X},4)-\{\mathcal{X}\}$.
\end{enumerate}
\end{lem}

\begin{figure}[ht]
\centering
\captionsetup{justification=centering}
\input{figurecheminLn}
\caption{Example of a path in Lemma~\ref{number of path in Ln} between $\mathcal{X}'$ (adjacent to $\mathcal{Y}_i$ with $i=2$) and $\mathcal{Y}_j$ with $j=6$.}\label{example path lemma number of path}
\end{figure}
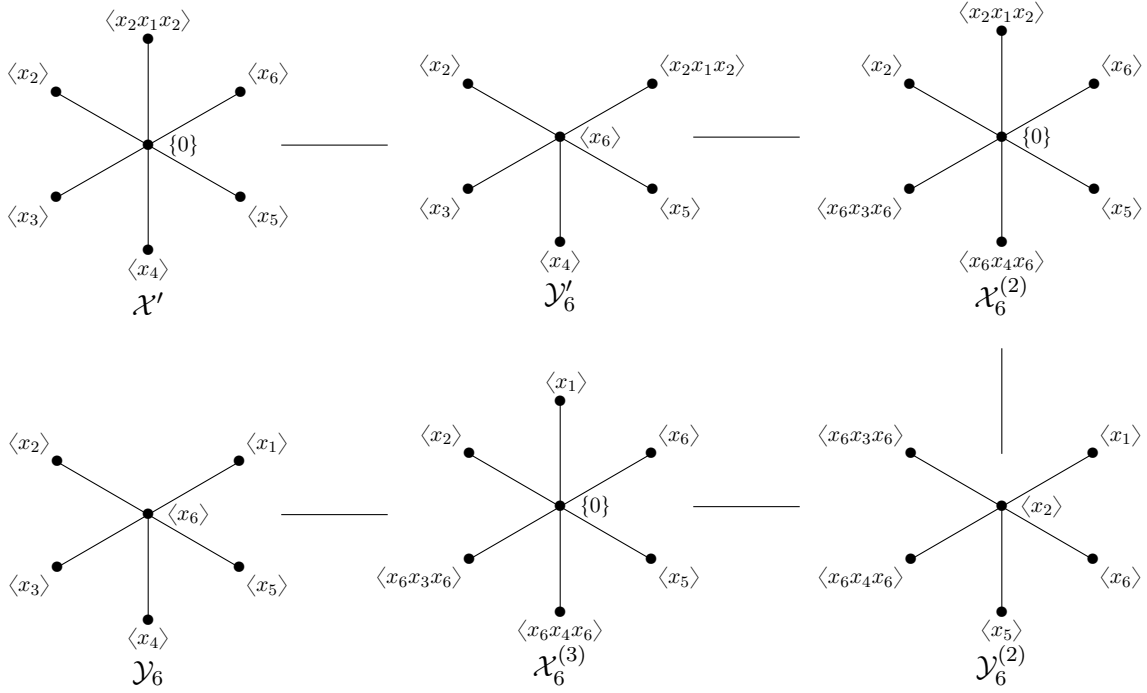

\dem We prove the case $x_j'=x_j$. The proof of the case $x_j'=x_ix_jx_i$ is similar. The proof consists in showing that the possible arcs $P$ are as represented in Figure~\ref{example path lemma number of path}.

Let $P$ be an arc (that is, an injective edge path), in $B(\mathcal{X},4)-\{\mathcal{X}\}$ between $\mathcal{X}'$ and $\mathcal{Y}_j$ of length at most $5$. Let $w_1,\ldots,w_n$ be the leaves of $\overline{X}'$, and, for $k \in \{1,\ldots,n\}$, let $y_k$ be the preimage by $\rho'$ of the generator of the group associated with $w_k$. Note that, up to reordering and composing by an inner automorphism, if $k \notin \{i_1,\ldots,i_{k_{\mathcal{X},i}(\mathcal{X}')}\}$, then $y_k=x_k$, and, if $k \in \{i_1,\ldots,i_{k_{\mathcal{X},i}(\mathcal{X}')}\}$, then $y_k=x_ix_kx_i$. For $k \neq i$, let $\mathcal{Y}_k'$ be the $F$-star such that a representative of $\mathcal{Y}_k'$ is obtained from $X'$ by contracting the edge adjacent to $w_k$.

\bigskip

\noindent{\bf Claim. } If $k \notin \{i,j\}$, the path $P$ cannot contain $\mathcal{Y}_k'$.

\bigskip

\dem Suppose towards a contradiction that $\mathcal{Y}_k' \in P$, with $k \notin \{i,j\}$. Since $P$ is an arc, there exists in $P$ a $\{0\}$-star $\mathcal{X}_k^{(2)}$ adjacent to $\mathcal{Y}_k'$ and distinct from $\mathcal{X}'$. By Lemma~\ref{properties of second term complexity}~$(2)$, we see that $l_{\mathcal{X}}(\mathcal{X}_k^{(2)}) \leq 2$. We claim that $l_{\mathcal{X}}(\mathcal{X}_k^{(2)})=2$.
Indeed, let $(X_k^{(2)},\psi)$ be a representative of $\mathcal{X}_k^{(2)}$, let $w_1^{(2)},\ldots,w_n^{(2)}$ be the $n$ leaves of the underlying graph of $X_k^{(2)}$ and, for $m \in \{1,\ldots,n\}$, let $y_m^{(2)}$ be the preimage by $\psi$ of the generator of the group associated with $w_m^{(2)}$. Note that, up to composition by an inner automorphism and reordering, for all $m \in \{1,\ldots,n\}$, $$y_m^{(2)}=x_i^{\gamma_m}x_k^{\beta_m}x_i^{\alpha_m}x_mx_i^{\alpha_m}x_k^{\beta_m}x_i^{\gamma_m}, \quad \alpha_m,\beta_m,\gamma_m \in \{0,1\}.$$

Note also that $\gamma_m=1$ precisely when $y_k=x_ix_kx_i$ and $\beta_m=1$. Thus, for all $m \in \{1,\ldots,n\}$, the element $y_m^{(2)}$ is obtained from $x_m$ using partial conjugations with twistors $x_i$ and $x_k$. Moreover, as $k \neq i$, and as $\mathcal{X}' \neq \mathcal{X}$, there exists $n_1$ such that either $\alpha_{n_1} \neq 0$ or $\gamma_{n_1} \neq 0$. As $\mathcal{X}_k^{(2)} \neq \mathcal{X}'$, there exists $n_2$ such that $\beta_{n_2} \neq 0$. It implies that $l_{\mathcal{X}}(\mathcal{X}_k^{(2)})=2$.

Therefore, by Lemma~\ref{properties of second term complexity}~$(1)$, the $\{0\}$-star $\mathcal{X}_k^{(2)}$ is not adjacent to $\mathcal{Y}_j$ since any $\{0\}$-star $\mathcal{Z}$ adjacent to $\mathcal{Y}_j$ is such that $l_{\mathcal{X}}(\mathcal{Z})=1$. 

So $P$ contains an $F$-star $\mathcal{Y}_k^{(2)}$ adjacent to $\mathcal{X}_k^{(2)}$ and distinct from $\mathcal{Y}_k'$ (see Figure~\ref{example path lemma number of path} with $k=6$). We claim that a representative of $\mathcal{Y}_k^{(2)}$ is obtained from $X_k^{(2)}$ by contracting the edge adjacent to $w_i^{(2)}$. Indeed, if it is not the case, one of the following two possibilities holds.

\medskip

\noindent{$(i)$ } A representative of $\mathcal{Y}_k^{(2)}$ is obtained from $X_k^{(2)}$ by contracting the edge adjacent to $w_k^{(2)}$. But then we go back to $\mathcal{Y}_k'$, which contradicts the fact that $P$ is an arc. 

\medskip

\noindent{$(ii)$ } A representative of $\mathcal{Y}_k^{(2)}$ is obtained from $X_k^{(2)}$ by contracting the edge adjacent to $w_{\ell}^{(2)}$, with $\ell \neq i,k$. Let $\mathcal{X}_k^{(3)}$ be the $\{0\}$-star in $P$ adjacent to $\mathcal{Y}_k^{(2)}$ and distinct from $\mathcal{X}_k^{(2)}$, and let $X_k^{(3)}$ be a representative of $\mathcal{X}_k^{(3)}$. Then, there exist $p \in \NN$, $m_1 \in \{1,\ldots,n\}$ and $i_1,\ldots,i_p,j_1,\ldots,j_p \in \{1,\ldots,n\}$ such that, if $y_{m_{1}}$ is the preimage by the marking of a nontrivial vertex group of $X_k^{(3)}$, we have 

$$\Big(\prod\limits_{m=1}^p\sigma_{i_m,j_m}\Big)(x_{m_{1}})=y_{m_{1}},$$

\noindent and there exists $m$ such that $j_m=\ell$. Therefore, by Lemma~\ref{second term complexity and 3 letters}, we see that $\ell_{\mathcal{X}}(\mathcal{X}_k^{(3)}) \geq 3$. But, by Lemma~\ref{properties of second term complexity}~$(2)$, we have $\mathcal{X}_k^{(3)} \notin B(\mathcal{X},4)$ and this contradicts the fact that $P \subseteq B(\mathcal{X},4)$. 

\medskip

Therefore a representative of $\mathcal{Y}_k^{(2)}$ is obtained from $X_k^{(2)}$ by contracting the edge adjacent to $w_i^{(2)}$ (see Figure~\ref{example path lemma number of path}).

But then, for every $\{0\}$-star $\mathcal{Z}$ adjacent to $\mathcal{Y}_k^{(2)}$, the set realizing $\ell_{\mathcal{X}}(\mathcal{Z})$ must contain $k$. Indeed, let $Z$ be a representative of $\mathcal{Z}$ and let $z_1,\ldots,z_n$ be the preimages by the marking of the generators of the nontrivial vertex groups. Up to composition by an inner automorphism and reordering, $z_i=y_i^{(2)}$ and there exists $m \in \{1,\ldots,n\}$ such that $z_m \neq y_m^{(2)}$. Thus, if $y_i^{(2)}$ is obtained from $x_i$ using partial conjugations such that one of the twistors is $x_k$, then a set realizing the minimum defining $\ell_{\mathcal{X}}(\mathcal{Z})$ contains $k$. Moreover, if $y_i^{(2)}=x_i$, then, since $n_2$ is such that $\beta_{n_2} \neq 0$, we see that $z_{n_2}$ is obtained from $x_{n_2}$ using partial conjugations such that one of the twistors is $x_k$. In any case, the set realizing $\ell_{\mathcal{X}}(\mathcal{Z})$ must contain $k$. On the other hand, by Lemma~\ref{properties of second term complexity}~$(1)$, if $\mathcal{Z}'$ is a $\{0\}$-star adjacent to $\mathcal{Y}_j$, then the set realizing $\ell_{\mathcal{X}}(\mathcal{Z}')$ only contains $j$. Since a set realizing the second term complexity is unique by Lemma~\ref{properties of second term complexity}~$(2)$, we see that $\mathcal{Z}$ and $\mathcal{Y}_j$ are not adjacent in $L_n$. This leads to a contradiction since we suppose that the length of $P$ is at most $5$. 
\hfill\qedsymbol

\bigskip

So the above claim implies that the path $P$ either contains $\mathcal{Y}_i$ or $\mathcal{Y}_j'$ (note that $\mathcal{Y}_i=\mathcal{Y}_i'$). The case $\mathcal{Y}_i \in P$ cannot occur by the following claim.

\medskip

\noindent{\bf Claim. } There does not exist a $\{0\}$-star adjacent to $\mathcal{Y}_i$ and distinct from $\mathcal{X}$ at distance $3$ from $\mathcal{Y}_j$ in $B(\mathcal{X},4)-\{\mathcal{X}\}$.

\medskip

\dem Let $\mathcal{X}_i^{(2)}$ be a $\{0\}$-star adjacent to $\mathcal{Y}_i$ and distinct from $\mathcal{X}$ and $\mathcal{X}'$. Let $X_i^{(2)}$ be a representative of $\mathcal{X}_i^{(2)}$. Let $v_1^{(2)},\ldots,v_n^{(2)}$ be the leaves of the underlying graph of $X_i^{(2)}$ and, for $m \in \{1,\ldots,n\}$, let $a_m^{(2)}$ be the preimage by the marking of the generator of the groups associated with $v_m^{(2)}$. Up to reordering and composing by an inner automorphism, we can suppose that, for all $m \in \{1,\ldots,n\}$, either $a_m^{(2)}=x_m$ or $a_m^{(2)}=x_ix_mx_i$ (this is possible by Lemma~\ref{properties of second term complexity}~$(1)$). Let $\mathcal{X}_i^{(3)}$ be a $\{0\}$-star distinct from $\mathcal{X}$ and at distance $2$ of $\mathcal{X}_i^{(2)}$, let $X_i^{(3)}$ be a representative of $\mathcal{X}_i^{(3)}$ and let $a_1^{(3)},\ldots,a_n^{(3)}$ be the preimages by the marking of the generators of the nontrivial vertex groups. Then one of the following holds:

\smallskip

\noindent{$(a)$ } The $\{0\}$-star $\mathcal{X}_i^{(3)}$ is adjacent to $\mathcal{Y}_i$. By Lemma~\ref{properties of second term complexity}~$(1)$, a set realizing $\ell_{\mathcal{X}}(\mathcal{X}_i^{(3)})$ is equal to $\{i\}$. On the other hand, the set realizing the minimum defining the second term complexity of every $\{0\}$-star adjacent to $\mathcal{Y}_j$ contains $j$. As $i \neq j$, we see that $\mathcal{X}^{(3)}$ cannot be adjacent to $\mathcal{Y}_j$.

\smallskip

\noindent{$(b)$ } There exist $p \in \NN$, $k \in \{1,\ldots,n\}-\{i\}$, \mbox{$\ell,i_1,\ldots,i_p,j_1,\ldots,j_p \in \{1,\ldots,n\}$} and \mbox{$s \in \{1,\ldots,p\}$} such that $$\Big(\prod\limits_{m=1}^p\sigma_{i_m,j_m}\Big)(x_{\ell})=a_{\ell}^{(3)}$$ and $j_s=k$. Thus $k$ is contained in any set realizing the minimum defining $\ell_{\mathcal{X}}(\mathcal{X}_i^{(3)})$. Moreover, we claim that $i$ is contained in any set realizing the minimum defining $\ell_{\mathcal{X}}(\mathcal{X}_i^{(3)})$. Indeed, as $k$ is contained in a set realizing the minimum defining $\ell_{\mathcal{X}}(\mathcal{X}_i^{(3)})$, a representative of $\mathcal{X}_i^{(3)}$ is obtained from $X_i^{(2)}$ as follows. We first contract the edge adjacent to the vertex $v_k^{(2)}$. This gives an $F$-star denoted by $Y_i^{(3)}$. Then, a representative of $\mathcal{X}_i^{(3)}$ is obtained from $Y_i^{(3)}$ by blowing-up an edge. If $a_k^{(2)}=x_k$, then, as $\mathcal{X}_i^{(2)}$ is adjacent to $\mathcal{Y}_i$, a set realizing the minimum defining $\ell_{\mathcal{X}}(\mathcal{X}_i^{(2)})$ is equal to $\{i\}$ by Lemma~\ref{properties of second term complexity}~$(1)$. As $a_k^{(2)}=x_k$, we see that either $a_m^{(3)}=a_m^{(2)}$ or $a_m^{(3)}=x_ka_m^{(2)}x_k$. Thus, as $i \neq k$, we see that $i$ is contained in a set realizing the minimum defining $\ell_{\mathcal{X}}(\mathcal{X}_i^{(3)})$. If $a_k^{(2)}=x_ix_kx_i$, then $a_k^{(3)}=a_k^{(2)}=x_ix_kx_i$ and any set realizing the minimum defining $\ell_{\mathcal{X}}(\mathcal{X}_i^{(3)})$ must contain $i$. Therefore, in any case, we have that $\{i,k\}$ is contained in any set realizing the minimum defining $\ell_{\mathcal{X}}(\mathcal{X}_i^{(3)})$. This shows that $\ell_{\mathcal{X}}(\mathcal{X}_i^{(3)}) \geq 2$. However, since the $\{0\}$-stars adjacent to $\mathcal{Y}_j$ have second term complexity equal to $1$ by Lemma~\ref{properties of second term complexity}~$(1)$, we see that $\mathcal{X}_i^{(2)}$ cannot be such that $d_{B(\mathcal{X},4)-\{\mathcal{X}\}}(\mathcal{X}_i^{(2)},\mathcal{Y}_j)=3$.
\hfill\qedsymbol

\bigskip

Thus, $P$ contains $\mathcal{Y}_j'$. As any two distinct $F$-stars are at distance at least $2$ in $L_n$, the path $P$ contains a $\{0\}$-star $\mathcal{X}_j^{(2)}$ adjacent to $\mathcal{Y}_j'$ and distinct from $\mathcal{X}'$ (see Figure~\ref{example path lemma number of path}). Let $(X_j^{(2)},\psi)$ be a representative of $\mathcal{X}_j^{(2)}$, let $w_1^{(2)},\ldots,w_n^{(2)}$ be the $n$ leaves of the underlying graph of $X_j^{(2)}$ and, for $m \in \{1,\ldots,n\}$, let $y_m^{(2)}$ be the preimage by $\psi$ of the generator of the group associated with $w_m^{(2)}$. Note that, up to composition by an inner automorphism and reordering, for all $m \in \{1,\ldots,n\}$,
$$y_m^{(2)}=x_j^{\beta_m}x_i^{\alpha_m}x_mx_i^{\alpha_m}x_j^{\beta_m}, \quad \alpha_m,\beta_m \in \{0,1\}.$$

As $\mathcal{X}_j^{(2)} \neq \mathcal{X}'$, there exist $k,l \in \{1,\ldots,n\}$ such that $\alpha_k \neq 0$ and $\beta_l \neq 0$. Thus, $\ell_{\mathcal{X}}(\mathcal{X}_j^{(2)})=2$ and a set realizing the minimum defining $\ell_{\mathcal{X}}(\mathcal{X}_j^{(2)})$ is $\{i,j\}$. This also implies that the $\{0\}$-star $\mathcal{X}_j^{(2)}$ is not adjacent to $\mathcal{Y}_j$ by Lemma~\ref{properties of second term complexity}~$(1)$. So $P$ contains an $F$-star $\mathcal{Y}_j^{(2)}$ adjacent to $\mathcal{X}_j^{(2)}$ and distinct from $\mathcal{Y}_j$ and $\mathcal{Y}_j'$ (see Figure~\ref{example path lemma number of path}).
We claim that a representative of $\mathcal{Y}_j^{(2)}$ is obtained from $X_j^{(2)}$ by contracting the edge that contains $w_i^{(2)}$. Indeed, if it is not the case, then one of the following holds.

\medskip

\noindent{$(i)$ } A representative of $\mathcal{Y}_j^{(2)}$ is obtained from $X_j^{(2)}$ by contracting the edge that contains $w_j^{(2)}$. Then $\mathcal{Y}_j^{(2)}=\mathcal{Y}_j'$, and this contradicts the fact that $P$ is an arc.

\medskip

\noindent{$(ii)$ } A representative of $\mathcal{Y}_j^{(2)}$ is obtained from $X_j^{(2)}$ by contracting the edge that contains $w_{\ell}^{(2)}$, with $\ell \neq i,j$. Let $\mathcal{X}_{\ell}^{(3)}$ be a $\{0\}$-star adjacent to $\mathcal{Y}_j^{(2)}$ and distinct from $\mathcal{X}_j^{(2)}$, and let $X_{\ell}^{(3)}$ be a representative of $\mathcal{X}_{\ell}^{(3)}$. As $\mathcal{X}_{\ell}^{(3)} \neq \mathcal{X}_j^{(2)}$ any set realizing the minimum defining $\ell_{\mathcal{X}}(\mathcal{X}_{\ell}^{(3)})$ must contain $\ell$. Accordingly, since $\ell \neq i,j$ we see by Lemma~\ref{second term complexity and 3 letters} that $\ell_{\mathcal{X}}(\mathcal{X}_{\ell}^{(3)}) \geq 3$. This contradicts the fact that $P \subseteq B(\mathcal{X},4)-\{\mathcal{X}\}$ by Lemma~\ref{properties of second term complexity}~$(2)$.

\medskip

Therefore, a representative of $\mathcal{Y}_j^{(2)}$ is obtained from $X_j^{(2)}$ by contracting the edge adjacent to $w_i^{(2)}$. We now distinguish two cases, according to the value of $\beta_i$.
\bigskip

\noindent{\bf Claim. } \noindent{$(1)$ } If $\beta_i=0$, then, for all $m \in \{1,\ldots,\hat{i},\ldots,n\}$, we have $(\alpha_m,\beta_m) \neq (1,1)$.

\medskip

\noindent{$(2)$ } If $\beta_i=1$, then, for all $m \in \{1,\ldots,\hat{i},\ldots,n\}$ such that $\alpha_m=1$, the pair $(\alpha_m,\beta_m)$ equals $(1,1)$.

\bigskip

\dem Let $\mathcal{Z}$ be a $\{0\}$-star adjacent to $\mathcal{Y}_j^{(2)}$, let $Z$ be a representative of $\mathcal{Z}$, and let $z_1,\ldots,z_n$ be the preimages by the marking of the generators of the nontrivial vertex groups of $Z$. 

\medskip

\noindent{$(1)$ } Suppose that $\beta_i=0$ and that there exists $m \in \{1,\ldots,\hat{i},\ldots,n\}$ such that $(\alpha_m,\beta_m)=(1,1)$. Then any set realizing the minimum defining $\ell_{\mathcal{X}}(\mathcal{Z})$ must contain $i$ because, as $y_m^{(2)}=x_jx_ix_mx_ix_j$, and as $y_i^{(2)}=x_i$, we see that, up to composing by an inner automorphism and reordering, we have that $z_i=x_i$ and either $z_m= x_jx_ix_mx_ix_j$ or $z_m= x_ix_jx_ix_mx_ix_jx_i$.

\medskip

\noindent{$(2)$ } Suppose now that $\beta_i=1$ and that there exists $m \in \{1,\ldots,\hat{i},\ldots,n\}$ such that $\alpha_m=1$ and such that the pair $(\alpha_m,\beta_m)$ equals $(1,0)$. Then any set realizing the minimum defining $\ell_{\mathcal{X}}(\mathcal{Z})$ must contain $i$ because, as $y_m^{(2)}=x_ix_mx_i$, and as $y_i^{(2)}=x_jx_ix_j$, we see that, up to composing by an inner automorphism and reordering, we have that $z_i=x_jx_ix_j$ and either $z_m= x_ix_mx_i$ or $z_m= x_jx_ix_jx_ix_mx_ix_jx_ix_j$. 

\medskip

So, in both cases, for every $\{0\}$-star $\mathcal{Z}$ adjacent to $\mathcal{Y}_j^{(2)}$, the set realizing the minimum defining $\ell_{\mathcal{X}}(\mathcal{Z})$ must contain $i$. 

Let $\mathcal{Z}$ be  the $\{0\}$-star in $P$ adjacent to $\mathcal{Y}^{(2)}$ and distinct from $\mathcal{X}^{(2)}$. Then the set realizing the minimum defining $\ell_{\mathcal{X}}(\mathcal{Z})$ must contain $i$ by the above. Since the length of $P$ is at most $5$, the $\{0\}$-star $\mathcal{Z}$ is adjacent to $\mathcal{Y}_j$.
But then, by Lemma~\ref{properties of second term complexity}~$(1)$, the set realizing the minimum defining $\ell_{\mathcal{X}}(\mathcal{Z})$ is equal to $j$. Since a set realizing the second term complexity is unique by Lemma~\ref{properties of second term complexity}~$(2)$, we see that $\mathcal{Z}$ and $\mathcal{Y}_j$ cannot be adjacent and this leads to a contradiction.
\hfill\qedsymbol

\bigskip

So if $\beta_i=0$, then, for all $m \in \{1,\ldots,\widehat{i},\ldots,n\}$, the pair $(\alpha_m,\beta_m) \neq (1,1)$ and if $\beta_i =1$, for all $m \in \{1,\ldots,\hat{i},\ldots,n\}$ such that $\alpha_m=1$, the pair $(\alpha_m,\beta_m)$ equals $(1,1)$. 

We now claim that there are exactly $2^{n-k_{\mathcal{X},i}(\mathcal{X}')-2}-1$ possible values for the sequence $(\beta_1,\ldots,\widehat{\beta_i},\ldots, \widehat{\beta_j},\ldots,\beta_n)$.

First, if $\beta_i=1$, then by the above claim, for all $m \in \{1,\ldots,\widehat{i},\ldots,n\}$ such that $\alpha_m=1$, we have $\beta_m=1$. Using a global conjugation by $x_j$, it then follows that every marked graph of groups whose associated sequence $(\beta_1,\ldots,\beta_n)$ satisfies the above claim and is such that $\beta_i=1$ is equivalent to a marked graph of groups whose associated sequence $(\beta_1',\ldots,\beta_n')$ satisfies the above claim and is such that $\beta_i'=0$. Thus we can suppose that, for such a sequence $(\beta_1,\ldots,\beta_n)$, we have $\beta_i=0$.

Moreover, by the above claim, all the pairs $(\alpha_m,\beta_m)$ such that $\alpha_m=1$ have the same value for $\beta_m$. Thus the sequence  $(\beta_1,\ldots,\hat{\beta_i},\ldots,\hat{\beta_j},\ldots,\beta_n)$ is determined by the pairs $(\alpha_m,\beta_m)$ such that $\alpha_m=0$ and the choice of $\beta_m$. By hypothesis, there are exactly $k_{\mathcal{X},i}(\mathcal{X}')$ values of $m \in \{1,\ldots,\widehat{i},\ldots,n\}$ such that $\alpha_m=1$ since $\alpha_m=1$ if and only if $x_m'=x_ix_mx_i$. It then suffices to choose whether $\beta_m=0$ or $\beta_m=1$. Furthermore, let $(\beta_1,\ldots,\beta_n)$ and $(\beta_1',\ldots,\beta_n')$ be two distinct sequences satisfying the above claim and such that $\beta_i=\beta_i'=0$. Then there exists $m \in \{1,\ldots,\widehat{i},\ldots,n\}$ such that $\beta_m=1$ and $\beta_m'=0$. Thus, since $\beta_i=\beta_i'=0$, the associated marked graph of groups are not equivalent and the two sequences give rise to two distinct equivalence classes of marked graph of groups. 
Finally, since $\mathcal{X}^{(2)} \neq \mathcal{X}'$, there exists $k \in \{1,\ldots,\widehat{i},\ldots,n\}$ such that $\beta_k=1$. Hence there are $2^{n-k_{\mathcal{X},i}(\mathcal{X}')-2}-1$ possible values for the sequence $(\beta_1,\ldots,\widehat{\beta_i},\ldots,\widehat{\beta_j},\ldots,\beta_n)$.

Let $\mathcal{Z}$ be a $\{0\}$-star adjacent to $\mathcal{Y}_j^{(2)}$ and distinct from $\mathcal{X}_j^{(2)}$ and let $Z$ be a representative of $\mathcal{Z}$. Let $z_1,\ldots,z_n$ be the preimage by the marking of the nontrivial associated groups. Then, for every sequence $(\beta_1,\ldots,\widehat{\beta_i},\ldots,\widehat{\beta_j},\ldots,\beta_n)$ satisfying the above claim, there exists exactly one such $\mathcal{Z}$ such that, up to composing by an inner automorphism and reordering, for all $\ell \in \{1,\ldots,n\}$, we have either $z_{\ell}=x_jx_{\ell}x_j$ or $z_{\ell}=x_{\ell}$. Such a $\{0\}$-star is adjacent to both $\mathcal{Y}_j^{(2)}$ and $\mathcal{Y}_j$. We call this $\{0\}$-star $\mathcal{X}_j^{(3)}$.

Thus, there exists a unique $\{0\}$-star $\mathcal{X}_j^{(3)}$ adjacent to both $\mathcal{Y}_j^{(2)}$ and $\mathcal{Y}_j$ . Since $P$ is an arc of length at most $5$, it must contain $\mathcal{X}_j^{(3)}$. Thus an arc in $B(\mathcal{X},4)-\{\mathcal{X}\}$ with length at most $5$ between $\mathcal{X}'$ and $\mathcal{Y}_j$ is completely determined by a sequence $(\beta_1,\ldots,\hat{\beta_i},\ldots,\hat{\beta_j},\ldots,\beta_n)$ satisfying the above claim. This concludes the proof of Lemma~\ref{number of path in Ln}~$(1)$.

In order to prove the third assertion of the lemma, let $P$ be an arc between $\mathcal{X}'$ and $\mathcal{Z}$ of length at most $4$. Then there exists an arc $P'$ between $\mathcal{X}'$ and $\mathcal{Y}_j$ of length at most $5$ which contains $P$. Thus, $P$ is contained in one of the paths constructed in the proof of the first assertion of the lemma. Therefore, using the notations of the proof of Lemma~\ref{number of path in Ln}~$(1)$, we see that $\mathcal{Z}=\mathcal{X}_j^{(3)}$. Let $X_j^{(3)}$ be a representative of $\mathcal{X}_j^{(3)}$ and let $y_1^{(3)},\ldots,y_n^{(3)}$ be the preimages by the marking of the generators of the nontrivial vertex groups. Then for all $m \in \{1,\ldots,n\}$, if $\alpha_m=0$, then $y_m^{(3)}=x_j^{\beta_m}x_mx_j^{\beta_m}$, and if $\alpha_m=1$, then either $y_m^{(3)}=x_m$ or $y_m^{(3)}=x_jx_mx_j$. Moreover, by construction, we know that there exists $m$ such that $\alpha_m=0$ and $\beta_m \neq 0$. As $k_{\mathcal{X},j}(\mathcal{Z})=1$, and as $\alpha_t=0$ if and only if $t \notin \{i_1,\ldots,i_{k_{\mathcal{X},i}}(\mathcal{X}')\}$, we see that there is an arc between $\mathcal{X}'$ and $\mathcal{Z}$ of length at most $4$ if and only if $t \notin \{i_1,\ldots,i_{k_{\mathcal{X},i}}(\mathcal{X}')\}$. This concludes the proof.
\hfill\qedsymbol

\bigskip

\begin{prop}\label{auto fixing star}
Let $n \geq 4$. Let $\mathcal{X} \in O_n$. Let $f \in \Aut(L_n)$ be such that $f$ restricted to the star of $\mathcal{X}$ is the identity. Then $f=\mathrm{id}_{L_n}$.
\end{prop}

\dem In order to prove Proposition~\ref{auto fixing star}, we prove that $f$ fixes the star of all  $\{0\}$-stars at distance $2$ from $\mathcal{X}$. This concludes by propagation since $L_n$ is connected. 

First, we prove that $f$ fixes $B(\mathcal{X},2) \cap O_n-\{\mathcal{X}\}$. Let $\mathcal{X}_1,\mathcal{X}_2 \in B(\mathcal{X},2) \cap O_n$ be distinct $\{0\}$-stars. If there exist distinct $i,j \in \{1,\ldots,n\}$ such that $\mathcal{X}_1$ is adjacent to $\mathcal{Y}_i$ and $\mathcal{X}_2$ is adjacent to $\mathcal{Y}_j$, then $f(\mathcal{X}_1) \neq \mathcal{X}_2$ because $f(\mathcal{Y}_i)=\mathcal{Y}_i$, $f(\mathcal{Y}_j)=\mathcal{Y}_j$ and there is no $\{0\}$-star adjacent to both $\mathcal{Y}_i$ and $\mathcal{Y}_j$ apart from $\mathcal{X}$. 

Suppose that there exists $i$ such that $\mathcal{Y}_i$ is adjacent to both $\mathcal{X}_1$ and $\mathcal{X}_2$. For $\alpha \in \{1,2\}$, let $X_{\alpha}$ be a representative of $\mathcal{X}_{\alpha}$ and let $y_1^{\alpha},\ldots,y_n^{\alpha}$ be the preimages by the marking of the generators of the nontrivial vertex groups of $X_{\alpha}$. Since $\mathcal{X}_1 \neq \mathcal{X}_2$, we see that, up to reordering and composing by an inner automorphism, there exist $j,k \in \{1,\ldots,n\}$ such that $y_i^1=y_i^2=x_i$, such that $y_j^1=x_j$ and $y_j^2=x_ix_jx_i$ and such that $y_k^1=y_k^2$. By Lemma~\ref{number of path in Ln}~$(1)$, if $k_{\mathcal{X},i}(\mathcal{X}_1) \neq k_{\mathcal{X},i}(\mathcal{X}_2)$, then the number of arcs of length at most $5$ in $B(\mathcal{X},4)-\{\mathcal{X}\}$ between $\mathcal{X}_1$ and $\mathcal{Y}_k$ is distinct from the number of arcs of length at most $5$ in $B(\mathcal{X},4)-\{\mathcal{X}\}$ between $\mathcal{X}_2$ and $\mathcal{Y}_k$. Suppose that $k_{\mathcal{X},i}(\mathcal{X}_1) + k_{\mathcal{X},i}(\mathcal{X}_2) \neq n-1$. In particular, we have that $n-k_{\mathcal{X},i}(\mathcal{X}_1)-2 \neq k_{\mathcal{X},i}(\mathcal{X}_2) -1$. Therefore, by Lemma~\ref{number of path in Ln}~$(1)$ and $(2)$, the number of arcs of length at most $5$ in $B(\mathcal{X},4)-\{\mathcal{X}\}$ between $\mathcal{X}_1$ and $\mathcal{Y}_j$ is distinct from the number of arcs of length at most $5$ in $B(\mathcal{X},4)-\{\mathcal{X}\}$ between $\mathcal{X}_2$ and $\mathcal{Y}_j$. Thus $f(\mathcal{X}_1) \neq \mathcal{X}_2$ since $f$ restricted to the star of $\mathcal{X}$ is the identity. In particular, since $n \geq 4$, if $\mathcal{X}' \in B(\mathcal{X},2)-\{\mathcal{X}\}$ is such that $\mathcal{X}'$ is adjacent to $\mathcal{Y}_i$ and that $k_{\mathcal{X},i}(\mathcal{X}')=1$, then $f(\mathcal{X}')=\mathcal{X'}$. 

It remains the case where $k_{\mathcal{X},i}(\mathcal{X}_1)=k_{\mathcal{X},i}(\mathcal{X}_2)=\frac{n-1}{2}$. Let $\mathcal{X}^k$ be the $\{0\}$-star adjacent to $\mathcal{Y}_k$ such that $k_{\mathcal{X},i}(\mathcal{X}^k)=1$ and such that the set realizing $k_{\mathcal{X},i}(\mathcal{X}^k)$ is $\{j\}$. As $k_{\mathcal{X},i}(\mathcal{X}^k)=1$, we have that $f(\mathcal{X}^k)=\mathcal{X}^k$. Moreover, as $y_j^1=x_j$ and as $y_j^2=x_ix_jx_i$, Lemma~\ref{number of path in Ln}~$(3)$ implies that there is no path of length at most $4$ between $\mathcal{X}_2$ and $\mathcal{X}^k$ in $B(\mathcal{X},4)-\{\mathcal{X}\}$ while there is one such path between $\mathcal{X}_1$ and $\mathcal{X}^k$. Thus $f(\mathcal{X}_1)=\mathcal{X}_1$. Hence $f$ fixes $B(\mathcal{X},2) \cap O_n-\{\mathcal{X}\}$.

Now let $\mathcal{X}' \in B(\mathcal{X},2) \cap O_n-\{\mathcal{X}\}$ and let $\mathcal{Y}$ be the $F$-star adjacent to both $\mathcal{X}$ and $\mathcal{X}'$ (the uniqueness of this $F$-star follows from the uniqueness of the set realizing the minimum defining $\ell_{\mathcal{X}}$, see Lemma~\ref{properties of second term complexity}~$(2)$). Let $X'$ be a representative of $\mathcal{X}'$ and let $v_1',\ldots,v_n'$ be the leaves of the underlying graph of $X'$ and, for $i \in \{1,\ldots,n\}$, let $x_i'$ be the preimage by the marking of $X'$ of the generators of the group associated with $v_i'$. Then, up to reordering, we can suppose that a representative of $\mathcal{Y}$ is obtained from $X'$ by contracting the edge adjacent to $v_n'$. Let $\mathcal{Y}^1$ and $\mathcal{Y}^2$ be two distinct $F$-stars adjacent to $\mathcal{X}'$ and distinct from $\mathcal{Y}$.  We prove that $f(\mathcal{Y}^1) \neq \mathcal{Y}^2$. Up to reordering, we can suppose that, for $\alpha \in \{1,2\}$, a representative of $\mathcal{Y}^{\alpha}$ is obtained from $X'$ by contracting the edge adjacent to $v_{\alpha}'$. Let $\mathcal{Z}$ be a $\{0\}$-star adjacent to $\mathcal{Y}$ such that~:

\smallskip

\noindent{$(1)$ } $k_{\mathcal{X}',n}(\mathcal{Z})=1$~;

\smallskip

\noindent{$(2)$ } a set realizing the minimum defining $k_{\mathcal{X}',n}(\mathcal{Z})$ is $\{1\}$.

\smallskip

Then Lemma~\ref{number of path in Ln}~$(1)$ and $(2)$ tells us that the number of paths of length at most $5$ in $B(\mathcal{X}',4)-\{\mathcal{X}'\}$ between $\mathcal{Z}$ and $\mathcal{Y}^1$ is equal to  $2^{k_{\mathcal{X}',n}(\mathcal{Z})-1}-1$ while the number of paths of length at most $5$ in $B(\mathcal{X}',4)-\{\mathcal{X}'\}$ between $\mathcal{Z}$ and $\mathcal{Y}^2$ is equal to $2^{n-k_{\mathcal{X}',n}(\mathcal{Z})-2}-1$. Since $k_{\mathcal{X}',n}(\mathcal{Z})=1$, since $n \geq 4$ and since $f$ restricted to the star of $\mathcal{Y}$ is the identity, we see that $f(\mathcal{Y}^1) \neq \mathcal{Y}^2$ and the proposition follows.
\hfill\qedsymbol

\bigskip

\noindent{\bf Proof of Theorem~\ref{rigidity Ln}. } The uniqueness of $\gamma$ is immediate since no automorphism of $W_n$ fixes the conjugacy class of each element appearing in every free generating set of $W_n$. It thus suffices to prove that every automorphism preserving $O_n$ and $F_n$ is induced by an element of $\Out(W_n)$. Let $f$ be an automorphism of $L_n$ preserving $O_n$ and $F_n$. Since $\Out(W_n)$ acts transitively on $O_n$, we can suppose, up to composing by an element of $\Out(W_n)$, that $f$ fixes a $\{0\}$-star $\mathcal{X}$. Now $\Stab_{\Out(W_n)}(\mathcal{X})$ is isomorphic to $\mathfrak{S}_n$ and every element of $\Stab_{\Out(W_n)}(\mathcal{X})$ acts on the underlying graph of a representative $X$ of $\mathcal{X}$ by permuting the leaves. As a representative of any $F$-star adjacent to $\mathcal{X}$ is obtained from $X$ by contracting the edge adjacent to a leaf, we see that $\Stab_{\Out(W_n)}(\mathcal{X})$ acts transitively on the link of $\mathcal{X}$. Thus, we can suppose, up to composing by an element of $\Out(W_n)$, that $f$ fixes the star of $\mathcal{X}$. Proposition~\ref{auto fixing star} then implies that $f$ is the identity. This concludes the proof of Theorem~\ref{rigidity Ln}.
\hfill\qedsymbol

\section{Rigidity of the outer space of $W_n$}\label{section rigidity K_n}

The aim of this section is to prove Theorem~\ref{Rigidity Kn}, by constructing an injective homomorphism $\Aut(K_n) \hookrightarrow \Aut(L_n)$. We  first give a characterization of the $\{0\}$-stars and the $F$-stars which is preserved under automorphisms of $K_n$. This characterization relies on a study of the link of the vertices of $K_n$. We begin with some definitions.

\begin{defi}
Let $X$ be a graph. A \emph{join} of $X$ is a decomposition of $X$ into two nontrivial subgraphs $A$ and $B$ such that $VA \cap VB = \varnothing$ and, for all $a \in VA$ and $b\in VB$, the vertices $a$ and $b$ are adjacent in $X$. We then write $X=A \ast B$.
\end{defi}

The fact of being decomposed as a join is preserved under automorphisms of graphs. In the case of a vertex $x \in VK_n$, there is a natural decomposition of the link $\lk(x)$ of $x$ in $K_n$.

\begin{defi}
Let $x=\mathcal{X} \in VK_n$. Let $X$ be a representative of $\mathcal{X}$.

\smallskip

\noindent{(1) } The \emph{positive link} of $x$, denoted by $\lk_+(x)$, is the maximal subgraph of $\lk(x)$ whose set of vertices consists in the homothety classes which have a representative that collapses onto $X$.

\smallskip

\noindent{(2) } The \emph{negative link} of $x$, denoted by $\lk_-(x)$, is the maximal subgraph of $\lk(x)$ whose set of vertices consists of homothety classes which have a representative $Y$ such that $X$ collapses onto $Y$.

\end{defi}

For all vertices $x$ of $K_n$, by definition of the adjacency in $K_n$, we have $$\lk(x)=\lk_+(x) \ast \lk_-(x).$$ It is in fact, as we will prove in Proposition~\ref{decomposition join} below, the only decomposition of $\lk(x)$ as a join. 

\begin{lem}\label{quasi maximal vertices lk_+=3}
Let $n \geq 4$. Let $x=\mathcal{X} \in VK_n$ be such that $\lk_+(x) \neq \varnothing$. Let $X$ be a representative of $\mathcal{X}$ and let $\overline{X}$ be its underlying graph.

\medskip

\noindent{$(1)$ } If $\lk_+(x)$ is nontrivial and has no edge, then $2 \leq |\lk_+(x)| \leq 3$. Moreover, $|\lk_+(x)|=3$ if and only if the underlying graph of any representative of $x$ has $n$ leaves.

\medskip

\noindent{$(2)$ } Let $\lk_+^1(x)$ be the set of vertices of $K_n$ such that any element of $\lk_+^1(x)$ has a representative that can be obtained from $X$ by blowing-up exactly one edge. Then $|\lk_+^1(x)| \geq 2$.
\end{lem}

\dem Suppose that $\lk_+(x)$ is nontrivial and has no edge. Then the graph $\overline{X}$ has at least $n-1$ leaves. Otherwise, one can blow-up two distinct edges at two distinct vertices of $\overline{X}$ with nontrival vertex groups which are not leaves. This gives rise to two vertices in the positive link of $x$ that are linked by an edge. This contradicts the fact that $\lk_+(x)$ has no edge. 

Moreover, if $\overline{X}$ has exactly $n-1$ leaves, then all vertices of $\overline{X}$ with trivial associated groups have valence $3$ since otherwise one can blow-up an edge at a non-leaf vertex of $\overline{X}$ with nontrivial vertex group and another edge at a valence-four vertex of $\overline{X}$ with trivial vertex group. This gives rise to two vertices in the positive link of $x$ that are linked by an edge. Moreover, the only non-leaf vertex with nontrivial associated group has valence equal to $2$ since otherwise one can blow-up two edges at this vertex, giving rise to two vertices in the positive link of $x$ that are linked by an edge.

If $\overline{X}$ has $n$ leaves, then at most one vertex of $\overline{X}$ has degree at least $4$ since otherwise one can blow-up two edges at two distinct vertices of $\overline{X}$. This gives rise to two vertices in the positive link of $x$ that are linked by an edge. Thus $\overline{X}$ has at most one vertex $v$ with degree at least $4$. The degree of $v$ is in fact equal to $4$ since otherwise one can blow up a two-edge graph at $v$, which gives rise to two vertices in the positive link of $x$ that are linked by an edge. 

Thus, there are two possibilities for $\overline{X}$.

\medskip

\noindent{$(i)$ } The graph $\overline{X}$ has $n$ leaves. Moreover, there are exactly one vertex $v$ of valence $4$ and $|V\overline{X}|-(n+1)$ vertices of valence $3$. In this case, the number of possible vertices in $\lk_+(x)$ corresponds to partitioning the set of edges adjacent to $v$ into two subsets of order $2$. This shows that $|\lk_+(x)|=3$.

\medskip

\noindent{$(ii)$ } The graph $\overline{X}$ has $n-1$ leaves. Moreover, there are exactly one vertex $v$ of valence $2$ and $|V\overline{X}|-n$ vertices of valence $3$. In this case, the group associated with $v$ is nontrivial and it is the only vertex of $\overline{X}$ that has nontrivial associated group and is not a leaf. In that case, the number of possible vertices in $\lk_+(x)$ corresponds to blowing-up an edge $e$ at $v$ so that one of the endpoint of $e$ is a leaf. Since $v$ has valence $2$, Proposition~\ref{Levitt stab} implies that $\Stab^0(x)$ is isomorphic to $F$. Thus, there are two possibilities for blowing-up the edge $e$ (either blowing it up while applying the nontrivial element of $\Stab^0(x)$ or blowing it up such that the preimages by the marking of the generators of the nontrivial vertex groups of the new graph of groups are the same as the preimages by the marking of $X$ of the generators of the nontrivial vertex groups). This shows that $|\lk_+(x)|=2$.

\bigskip

We now prove the second part of the lemma. Suppose that $\lk_+(x)$ is nontrivial (it might have edges). Suppose first that $\overline{X}$ has at most $n-2$ leaves. Let $v_1$ and $v_2$ be two vertices of $\overline{X}$ with nontrivial associated groups that are not leaves. Then one can find two elements of $\lk_+^1(x)$ by blowing up an edge at either $v_1$ or $v_2$. Thus, $|\lk_+^1(v)| \geq 2$.

Finally, if $\overline{X}$ has at least $n-1$ leaves, then the constructions of distinct elements of $\lk_+^1(x)$ are similar to the case where $\lk_+(x)$ is nontrivial and has no edge.
\hfill\qedsymbol

\begin{lem}\label{quasi minimal number of vertices}
Let $n \geq 4$. Suppose that $x=\mathcal{X} \in VK_n$ is such that $\lk_-(x)$ is nontrivial and has no edge. Let $X$ be a representative of $\mathcal{X}$ and $\overline{X}$ be its underlying graph.

\medskip

\noindent{$(1)$ } There exists a unique vertex in $\overline{X}$ with trivial associated group.

\medskip

\noindent{$(2)$ } The negative link satisfies $3 \leq |\lk_-(x)| \leq n$. Moreover, $|\lk_-(x)|=n$ if and only if $x$ is a $\{0\}$-star.
\end{lem}

\dem $(1)$ The graph $\overline{X}$ contains at least one vertex with trivial associated group since otherwise there would not exist an element $\mathcal{Y} \in VK_n$ such that a representative of $\mathcal{Y}$ is obtained from $X$ by collapsing a forest. This would contradict the fact that $\lk_-(x)$ is nontrivial. Thus $\overline{X}$ contains at least one vertex with trivial associated group. 

Suppose towards a contradiction that $\overline{X}$ contains two vertices with trivial associated groups. Then, since the degree of any vertex of $\overline{X}$ with trivial associated group is at least $3$, there exists two distinct edges $e_1$ and $e_2$ in $\overline{X}$ that can be simultaneously collapsed to get a new element in $VK_n$. Moreover, if $i \in \{1,2\}$, and if $\mathcal{Y}_i$ is the homothety classes of the marked graph of groups obtained from $X$ by collapsing $e_i$, then $\mathcal{Y}_1,\mathcal{Y}_2 \in \lk_-(x)$ and $\mathcal{Y}_1$ and $\mathcal{Y}_2$ are adjacent in $\lk_-(x)$ and distinct. This contradicts the fact that $\lk_-(x)$ has no edge. Thus, there exists a unique vertex in $\overline{X}$ with trivial associated group. 

\medskip

\noindent{$(2)$ } Let $v$ be the unique vertex in $\overline{X}$ with trivial associated group guaranteed by the first assertion. It follows that $\operatorname{deg}(v) \geq 3$. Thus $|\lk_-(x)| \geq 3$. Since $\overline{X}$ contains exactly $n$ vertices with nontrivial associated group, $\operatorname{deg}(v) \leq n$. Thus $|\lk_-(x)| \leq n$.

Now, if $|\lk_-(x)|=n$, then $\operatorname{deg}(v)=n$. Thus $\overline{X}$ contains exactly $n$ leaves and $n+1$ vertices and $\mathcal{X}$ is a $\{0\}$-star.
Conversely, if $\mathcal{X}$ is a $\{0\}$-star, then there exists exactly one vertex in $\overline{X}$ with trivial associated group. Moreover, its degree is equal to $n$. Thus $|\lk_-(x)|=n$.
\hfill\qedsymbol

\begin{lem}\label{partition link_-(x)}
Let $n \geq 4$. Let $x=\mathcal{X} \in VK_n$ be such that $\lk_-(x)$ is nontrivial. Let $X$ be a representative of $\mathcal{X}$ and let $\overline{X}$ be its underlying graph. Let $v_1,\ldots,v_n$ be the vertices of $\overline{X}$ with nontrivial associated group.  Let $e \in E\overline{X}$ and let $\{v_{i_1},\ldots,v_{i_k}\} \amalg \{v_{j_1},\ldots,v_{j_l}\}$ be the partition of $\{v_1,\ldots,v_n\}$ obtained by considering the vertices contained in each connected component of $\overline{X}-\mathring{e}$.

\begin{enumerate}
\item Let $F_0 \subseteq \overline{X}$ be a forest (that may be empty) such that the homothety class of the marked graph of groups $Y$ obtained from $X$ by collapsing $F_0$ is a vertex of $K_n$. Let $p \colon \overline{X} \to \overline{Y}$ be the canonical projection. 

Then, if $p(e)$ is not a vertex, it is the unique edge $f$ of $\overline{Y}$ such that the partition of $\{p(v_1),\ldots,p(v_n)\}$ induced by $\overline{Y}-\mathring{f}$ is $\{p(v_{i_1}),\ldots,p(v_{i_k})\} \amalg \{p(v_{j_1}),\ldots,p(v_{j_l})\}$.

\item Let $y,z \in \lk_-(x)$ be distinct vertices. Let $Y$ and $Z$ be representatives of $y$ and $z$ respectively, and let $\overline{Y}$ and $\overline{Z}$ be their underlying graphs. Let $p_{y} \colon \overline{X} \to \overline{Y}$ and $p_{z} \colon \overline{X} \to \overline{Z}$ be the natural projections. 

If one can obtain $Z$ from $Y$ by collapsing a forest of $\overline{Y}$, and if $p_{z}(e)$ is not a point, there exists a unique edge $\widetilde{p_{z}(e)} \in E\overline{Y}$ such that the partition of $\{p_y(v_1),\ldots,p_y(v_n)\}$ induced by $\widetilde{p_{z}(e)}$ is $$\{p_y(v_{i_1}),\ldots,p_y(v_{i_k})\} \amalg \{p_y(v_{j_1}),\ldots,p_y(v_{j_l})\}.$$
\end{enumerate}
\end{lem}

\begin{rmq}\label{free factor decomposition and canonical lift}
\rm{ Let $X$, $Y$ and $Z$ be as in the above statement. Let $G$ be the forest of $Y$ such that $Z$ is obtained from $Y$ by collapsing $G$.

\medskip

\noindent{(1) } The statements of the lemmas can be reinterpreted in terms of decompositions in free factors of $W_n$. Indeed, a partition of the vertices with nontrivial associated groups $\{v_1,\ldots,v_n\}=A \amalg B$ induced by an edge of $\overline{X}$ gives rise to a decomposition of $W_n$ as $W_n= W_k^A \ast W_{n-k}^B$ well-defined up to global conjugation. In this case, $W_k^A$ is generated by the groups associated with the vertices in $A$, and $W_{n-k}^B$ is generated by the groups associated with the vertices in $B$. In particular, Lemma~\ref{partition link_-(x)}~$(1)$ can be stated as follows. 

\medskip

{\it If $X$ is a graph of groups whose fundamental group is $W_n$, and if $e$ and $f$ are distinct edges of the underlying graph of $X$, then $e$ and $f$ induce distinct free factor decompositions. 

Moreover, if $Y$ is a graph of groups obtained from $X$ by collapsing a forest, and if $g$ is an edge of the underlying graph of $Y$, then there exists a unique edge $\widetilde{g}$ in the underlying graph of $X$ which induces the same free factor decomposition as $g$.}

\medskip

\noindent{(2) } Let $H$ be a forest in $\overline{Z}$. The second statement of the lemma gives a unique minimal forest $\widetilde{H}$ in $\overline{Y}$ that lifts $H$. Indeed, if $h \in EH$, let $\widetilde{h}$ be the unique edge of $\overline{Y}$ given by Lemma~\ref{partition link_-(x)}~$(2)$. Then $\{\widetilde{h}\}_{h \in H}$ is a lift of $H$. This lift has the property that $\widetilde{H} \cap G$ is contained in the leaves of $G$ and that every vertex of $\widetilde{H}$ is adjacent to an edge in $\widetilde{H}$. We call it \emph{the canonical lift} of $H$.
}
\end{rmq}

\dem  For the first statement, we only need to prove the uniqueness result. Let $f$ be an edge of $\overline{Y}$ distinct from $p(e)$. Let $A_1 \amalg A_2$ be the partition of $\{p(v_1),\ldots,p(v_n)\}$ induced by $p(e)$, and let $B_1 \amalg B_2$ be the partition of $\{p(v_1),\ldots,p(v_n)\}$ induced by $f$. We prove that there exist two vertices $v$ and $w$ of $\overline{Y}$ with nontrivial associated groups such that $v$ and $w$ are in the same connected component of $\overline{Y}-\mathring{f}$ while they are not in the same connected component of $\overline{Y}-\mathring{p(e)}$, or conversely. This will imply that there exists $\alpha \in \{1,2\}$ such that $B_{\alpha} \cap A_1 \neq \varnothing$ and that $B_{\alpha} \cap A_2 \neq \varnothing$, or conversely. This will conlcude the proof. There are two cases to distinguish, according to the endpoints of $p(e)$. 

If both of the endpoints of $p(e)$ have nontrivial associated groups, then, since $\overline{Y}$ is a tree, $p(e)$ is necessarily the unique edge of $\overline{Y}$ such that the endpoints of $p(e)$ are in distinct connected components of $\overline{Y}-\mathring{p(e)}$.

\begin{figure}
\centering
\captionsetup{justification=centering}

\begin{tikzpicture}[scale=0.90]
\draw (-3.5,0) -- (0.5,0);
\draw (-2,0) node[below, scale=1] {$f$};
\draw (-2.5,0) node {$\bullet$};
\draw (-1.5,0) node {$\bullet$};
\draw (-0.5,0) -- (0.5,0);
\draw (0,0) node[below, scale=1] {$p(e)$};
\draw (-0.5,0) node {$\bullet$};
\draw (-0.5,0) node[above, scale=1] {$p(v_{\ell})$};
\draw (0.5,0) node {$\bullet$};
\draw (0.5,0) node[above, scale=1] {$v_0$};
\draw (-3.5,0) node {$\bullet$};
\draw (-3.5,0) node[left, scale=1] {$p(v_j)$};
\draw (2.5,1) node {$\bullet$};
\draw (2.5,1) node[right, scale=1] {$p(v_i)$};
\draw (2.5,-1) node {$\bullet$};
\draw (2.5,-1) node[right, scale=1] {$p(v_k)$};
\draw (0.5,0) -- (2.5,1);
\draw (0.5,0) -- (2.5,-1);
\draw (-1,-1) node {Case $2$};

\draw (0.5,0+4) -- (4.5,0+4);
\draw (2,0+4) node[below, scale=1] {$f$};
\draw (1.5,0+4) node {$\bullet$};
\draw (2.5,0+4) node {$\bullet$};
\draw (3.5,0+4) -- (4.5,0+4);
\draw (4,0+4) node[below, scale=1] {$p(e)$};
\draw (3.5,0+4) node {$\bullet$};
\draw (3.5,0+4) node[above left, scale=1] {$v_0$};
\draw (4.5,0+4) node {$\bullet$};
\draw (3.5,0+4) -- (3.5,1.5+4);
\draw (3.5,1.5+4) node {$\bullet$};
\draw (3.5,1.5+4) node[above, scale=1] {$p(v_{k})$};
\draw (0.5,0+4) node {$\bullet$};
\draw (0.5,0+4) node[left, scale=1] {$p(v_j)$};
\draw (6.5,0+4) node {$\bullet$};
\draw (6.5,0+4) node[right, scale=1] {$p(v_i)$};
\draw (4.5,0+4) -- (6.5,0+4);
\draw (3.5,3) node {Case $1$};

\draw (-3.5+8,0) -- (0.5+8,0);
\draw (-2+8,0) node[below, scale=1] {$f$};
\draw (-2.5+8,0) node {$\bullet$};
\draw (-1.5+8,0) node {$\bullet$};
\draw (-0.5+8,0) -- (0.5+8,0);
\draw (0+8,0) node[below, scale=1] {$p(e)$};
\draw (-0.5+8,0) node {$\bullet$};
\draw (-0.5+8,0) node[above left, scale=1] {$v_1$};
\draw (-0.5+8,0) -- (-0.5+8,1.5);
\draw (-0.5+8,1.5) node {$\bullet$};
\draw (-0.5+8,1.5) node[above, scale=1] {$p(v_{\ell})$};
\draw (0.5+8,0) node {$\bullet$};
\draw (0.5+8,0) node[above, scale=1] {$v_0$};
\draw (-3.5+8,0) node {$\bullet$};
\draw (-3.5+8,0) node[left, scale=1] {$p(v_j)$};
\draw (2.5+8,1) node {$\bullet$};
\draw (2.5+8,1) node[right, scale=1] {$p(v_i)$};
\draw (2.5+8,-1) node {$\bullet$};
\draw (2.5+8,-1) node[right, scale=1] {$p(v_k)$};
\draw (0.5+8,0) -- (2.5+8,1);
\draw (0.5+8,0) -- (2.5+8,-1);
\draw (7.5,-1) node {Case $3$};

\end{tikzpicture}
\caption{The arcs constructed in Lemma~\ref{partition link_-(x)}.}\label{partition lk_-(x) figure}
\end{figure}
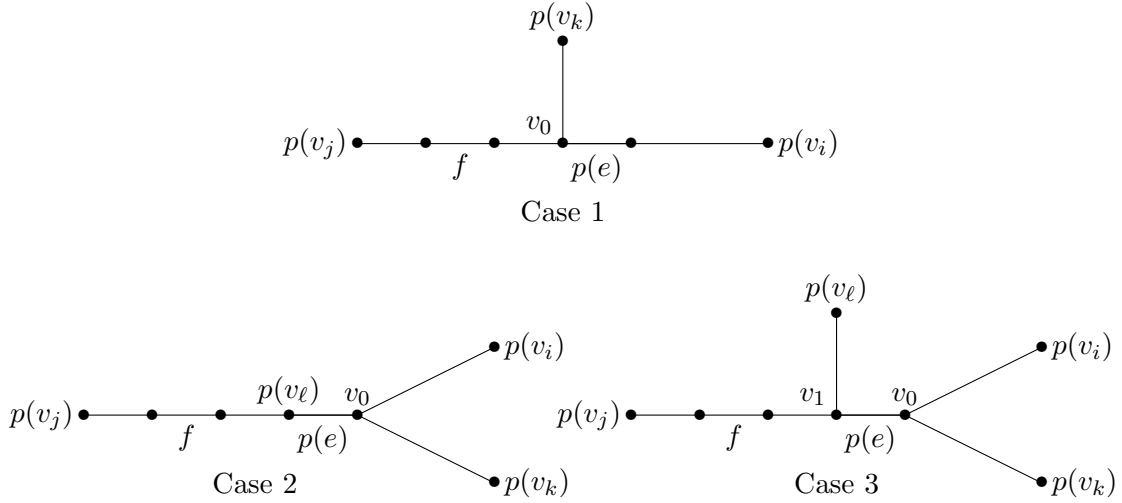

Suppose that one of the endpoints of $p(e)$, denoted by $v_0$, has trivial associated group. Then there exists an arc $P$ between two distinct leaves of $\overline{Y}$, say $p(v_i)$ and $p(v_j)$, such that $p(e)$ and $f$ are (up to replacing them by their opposite edges) contained in this path and in this order. Since $v_0$ has trivial associated group, $\operatorname{deg}(v_0) \geq 3$. Thus, up to exchanging the roles of $p(v_i)$ and $p(v_j)$, there exists a path $P'$ between $p(v_i)$ and a leaf of $\overline{Y}$, say $p(v_k)$, distinct from both $p(v_i)$ and $p(v_j)$, such that $P'$ contains $v_0$ (see Figure~\ref{partition lk_-(x) figure}). 

So if $P'$ contains $p(e)$ (see Figure~\ref{partition lk_-(x) figure}, Case~$1$), then $p(v_i)$ and $p(v_k)$ are not contained in the same connected component of $\overline{X}-\mathring{p(e)}$ while they are contained in the same connected component of  $\overline{X}-\mathring{f}$. If $P'$ does not contain $p(e)$, then there are two cases to distinguish. 

Let $v_1$ be the other endpoint of $p(e)$. If there exists $\ell \in \{1,\ldots,n\}$ such that we have $v_1=p(v_{\ell})$ (see Figure~\ref{partition lk_-(x) figure}, Case~$2$), then $p(v_j)$ and $p(v_{\ell})$ are contained in the same connected component of $\overline{X}-\mathring{p(e)}$ while they are not contained in the same connected component of  $\overline{X}-\mathring{f}$. 

If $v_1$ has trivial associated group (see Figure~\ref{partition lk_-(x) figure}, Case~$3$), then $\operatorname{deg}(v_1) \geq 3$. So there exists $\ell \in \{1,\ldots,n\}$ such that $v_1$ is contained in the arc $P^{(2)}$ between $p(v_j)$ and $p(v_{\ell})$ and such that $p(e)$ is not contained in $P^{(2)}$. Thus $p(v_j)$ and $p(v_{\ell})$ are contained in the same connected component of $\overline{X}-\mathring{p(e)}$ while they are not contained in the same connected component of  $\overline{X}-\mathring{f}$. 
In any case, $p(e)$ and $f$ do not generate the same partition of $\{p(v_1),\ldots,p(v_n)\}$.

\bigskip

Let $Y$ and $Z$ be as in the second statement of the lemma. By the first statement of the lemma, there exists a unique edge $\widetilde{p_{z}(e)} \in E\overline{Y}$ such that the partition of $\{p_y(v_1),\ldots,p_y(v_n)\}$ induced by $\widetilde{p_z(e)}$ is $\{p_y(v_{i_1}),\ldots,p_y(v_{i_k})\} \amalg \{p_y(v_{j_1}),\ldots,p_y(v_{j_l})\}$ (namely it is $p_y(e)$), and we take this edge to be our lift.
\hfill\qedsymbol

\medskip

\begin{prop}\label{decomposition join}
Let $n \geq 4$, and $x=\mathcal{X} \in VK_n$. Suppose that both $\lk_-(x)$ and $\lk_+(x)$ are nontrivial. The only nontrivial decomposition of $\lk(x)$ as a join is $\lk(x)=\lk_+(x) \ast \lk_-(x)$.
\end{prop}

\dem Let $X$ be a representative of $\mathcal{X}$ and let $\overline{X}$ be its underlying graph. Let $\lk(x)=A \ast B$ be a nontrivial decomposition as a join of $\lk(x)$ such that $A \neq \lk_+(x),\lk_-(x)$. Then there exist $x_1,x_2 \in \lk_+(x)$ or $x_1,x_2 \in \lk_-(x)$ such that $x_1 \in A$ and $x_2 \in B$. For $i \in \{1,2\}$, let $\mathcal{X}_i$ be the homothety class corresponding to $x_i$ and let $X_i$ be a representative. Let $\overline{X}_i$ be the underlying graph of $X_i$. Since $x_1$ and $x_2$ are joined by an edge, up to renumbering and changing the representatives, there exists a forest $F_0$ in $\overline{X}_1$ such that $\overline{X}_2$ is obtained from $\overline{X}_1$ by collapsing $F_0$.  We now investigate both cases. 

\medskip

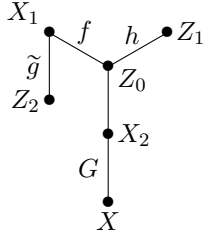
\begin{figure}
\centering
\captionsetup{justification=centering}

\begin{tikzpicture}[scale=0.90]

\draw (0,0) node {$\bullet$};
\draw (0,0) node[below, scale=0.9] {$X$};
\draw (0,0) -- (0,1);
\draw (0,0.5) node[left, scale= 0.9] {$G$};
\draw (0,1) node {$\bullet$};
\draw (0,1) node[right, scale=0.9] {$X_2$};
\draw (0,1) -- (0,2);
\draw (0,2) node {$\bullet$};
\draw (0,2.1) node[below right, scale=0.9] {$Z_0$};
\draw (0,2) -- (0.86,2.5);
\draw (0.86,2.5) node {$\bullet$};
\draw (0.86,2.5) node[right, scale=0.9] {$Z_1$};
\draw (0,2) -- (-0.86,2.5);
\draw (-0.35,2.2) node[above, scale=0.9] {$f$};
\draw (0.35,2.2) node[above, scale=0.9] {$h$};
\draw (-0.86,2.5) node {$\bullet$};
\draw (-0.86,2.5) node[above left,scale=0.9] {$X_1$};
\draw (-0.86,2.5) -- (-0.86,1.5);
\draw (-0.86,2) node[left, scale=0.9] {$\widetilde{g}$};
\draw (-0.86,1.5) node {$\bullet$};
\draw (-0.86,1.5) node[left,scale=0.9] {$Z_2$};
\end{tikzpicture}
\caption{The adjacency of the homothety classes constructed in the first case of Lemma~\ref{decomposition join}.}\label{decomposition join figure case 1}
\end{figure}

Suppose first that $x_1,x_2 \in \lk_+(x)$. We are going to construct two other vertices $z_1$ and $z_2$ such that $z_1 \in A$, $z_2 \in B$ and $z_1$ and $z_2$ are not linked by an edge, which will lead to a contradiction (see Figure~\ref{decomposition join figure case 1}).

Since $x_2 \in \lk_+(x)$, up to changing the representative $X$ of $\mathcal{X}$, there exists a forest $G$ in $\overline{X}_2$ such that $\overline{X}$ is obtained from $\overline{X}_2$ by collapsing $G$. Let $\widetilde{G}$ be the canonical lift of $G$ in $\overline{X}_1$.

Let $f \in EF_0$. Let $\mathcal{Z}_0$ be the homothety class of the marked graph of groups $Z_0$ obtained from $X_1$ by collapsing $f$ and let $z_0$ be the corresponding vertex in $K_n$. Since a representative of $\mathcal{Z}_0$ is obtained from $X_1$ by collapsing an edge, we see that $x_1 \in \lk_+(z_0)$. Moreover, since $f \in EF_0$, we see that $z_0 \in \lk_+(x_2)$ and $z_0 \in \lk_+(x)$. Lemma~\ref{quasi maximal vertices lk_+=3}~$(2)$ applied to $z_0$ then implies that there exists $z_1=\mathcal{Z}_1 \in \lk_+(z_0)$ distinct from $x_1$ such that the underlying graph of any representative of $z_1$ has the same number of edges as $\overline{X}_1$. Since $z_1 \in \lk_+(z_0)$ and $z_0 \in \lk_+(x)$, we have $z_1 \in \lk_+(x)$. As $z_1$ has a representative that has the same number of edges as $X_1$, and as $x_1 \neq z_1$, we see that $z_1 \notin \lk(x_1)$. Therefore we have $z_1 \in A$. 

In order to construct $z_2$, let $\widetilde{g}$ be an edge in $\widetilde{G}$. Let $\mathcal{Z}_2$ be the homothety class of the marked graph of group $Z_2$ obtained from $X_1$ by collapsing $\widetilde{g}$. Let $z_2$ be the corresponding vertex in $K_n$. Then, since $\widetilde{g} \in \widetilde{G}$, we see that $z_2 \in \lk_+(x)$. As $\widetilde{g} \in E\widetilde{G}$, and as two distinct edges induce distinct free factor decompositions by Remark~\ref{free factor decomposition and canonical lift}~$(1)$, there exists an edge $g \in E\overline{X}_2$ (namely the edge whose lift in $\overline{X}_1$ is $\widetilde{g})$ such that the free factor decomposition induced by $g$ is distinct from the free factor decomposition induced by any edge of the underlying graph of $Z_2$. Thus we see that $x_2$ and $z_2$ cannot be adjacent. Indeed, if it was the case, then as $Z_2$ is obtained from $X_1$ by collapsing exactly one edge, either $|EX_2|=|EZ_2|$ or there would exist a representative $Z_2'$ of $\mathcal{Z}_2$ such that $X_2$ is obtained from $Z_2'$ by collapsing a forest. Both cases would contradict Remark~\ref{free factor decomposition and canonical lift}~$(1)$ because the edge $g$ of $\overline{X}_2$ induces a free factor decomposition that is not induced by any edge of the underlying graph of $Z_2'$. This implies that $z_2 \notin \lk(x_2)$ and that $z_2 \in B$. 

\medskip

\noindent{\bf Claim. } The vertices $z_1$ and $z_2$ are not adjacent in $\lk(v)$. 

\medskip

\dem Suppose towards a contradiction that $z_1$ and $z_2$ are adjacent. Let $Z_1$ be a representative of the homothety class corresponding to $z_1$, and, for $i \in \{1,2\}$, let $\overline{Z}_i$ be the underlying graph of $Z_i$. As $|E\overline{Z}_1|=|E\overline{Z}_2|+1$, up to changing the representatives $Z_1$ and $Z_2$, we can suppose that $Z_2$ is obtained from $Z_1$ by collapsing an edge $e \in E\overline{Z}_1$. Let $h$ be the edge in $\overline{Z}_1$ such that the marked graph of groups obtained from $Z_1$ by collapsing $h$ is in $z_0$. As $\mathcal{Z}_1$ is distinct from $\mathcal{X}_1$, the edge $h$ is such that the free factor decomposition of $W_n$ induced by $h$ is distinct from the one induced by any edge of $\overline{X}_1$. Thus, by Lemma~\ref{partition link_-(x)}~$(2)$, the free factor decomposition of $W_n$ induced by $h$ is distinct from the one induced by any edge of $\overline{Z}_2$. Therefore, by Remark~\ref{free factor decomposition and canonical lift}~$(1)$, the marked graph of groups $Z_2$ is obtained from $Z_1$ by collapsing $h$. This implies that $z_2=z_0$ by the choice of $h$. However, $\widetilde{G} \cap F_0$ does not contain any edge by the properties of the canonical lift of $G$ (see Remark~\ref{free factor decomposition and canonical lift}~$(2)$). Thus there exists an edge in $\overline{Z}_0$ which induces the same free factor decomposition of $W_n$ as $\widetilde{g}$. As $Z_2$ is obtained from $X_1$ by collapsing $\widetilde{g}$, Lemma~\ref{partition link_-(x)}~$(1)$ implies that there is no edge in $\overline{Z}_2$ that induces the same free factor decomposition as $\widetilde{g}$. Thus $z_0 \neq z_2$, and this leads to a contradiction.
\hfill\qedsymbol

\medskip

Therefore $z_1$ and $z_2$ are not adjacent in $\lk(x)$. However, $z_1 \in A$ and $z_2 \in B$. This contradicts the fact that $\lk(x)=A \ast B$ is a join decomposition.

\bigskip

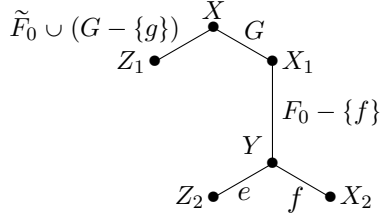
\begin{figure}
\centering
\captionsetup{justification=centering}

\begin{tikzpicture}[scale=0.90]

\draw (0,0) node {$\bullet$};
\draw (0,0) node[above, scale=0.9] {$X$};
\draw (0,0) -- (0.86,-0.5);
\draw (0.6,-0.3) node[above, scale= 0.9] {$G$};
\draw (0.86,-0.5) node {$\bullet$};
\draw (0.86,-0.5) node[right, scale=0.9] {$X_1$};
\draw (0,0) -- (-0.86,-0.5);
\draw (-0.35,-0.35) node[above left, scale= 0.9] {$\widetilde{F}_0 \cup (G-\{g\})$};
\draw (-0.86,-0.5) node {$\bullet$};
\draw (-0.86,-0.5) node[left, scale=0.9] {$Z_1$};
\draw (0.86,-0.5) -- (0.86,-1.5-0.5);
\draw (0.86,-0.75-0.5) node[right, scale= 0.9] {$F_0-\{f\}$};
\draw (0.86,-1.5-0.5) node {$\bullet$};
\draw (0.86,-1.5-0.5) node[above left, scale=0.9] {$Y$};
\draw (0.86,-1.5-0.5) -- (1.71,-2-0.5);
\draw (1.71,-2-0.5) node {$\bullet$};
\draw (1.71,-2-0.5) node[right, scale=0.9] {$X_2$};
\draw (1.2,-1.7-0.5) node[below] {$f$};
\draw (0.86,-1.5-0.5) -- (0,-2-0.5);
\draw (0,-2-0.5) node {$\bullet$};
\draw (0.45,-1.7-0.5) node[below] {$e$};
\draw (0,-2-0.5) node[left, scale=0.9] {$Z_2$};
\end{tikzpicture}
\caption{The adjacency of the homothety classes constructed in the second case of Lemma~\ref{decomposition join}.}\label{decomposition join figure case 2}
\end{figure}

Now suppose that $x_1,x_2 \in \lk_-(x)$. We use the same strategy as when $x_1,x_2 \in \lk_+(x)$ (see Figure~\ref{decomposition join figure case 2}).

Since $x_1 \in \lk_-(x)$, up to changing the representative $X$ of $\mathcal{X}$, there exists a forest $G$ in $\overline{X}$ such that $\overline{X}_1$ is obtained from $\overline{X}$ by collapsing $G$. Let $g \in EG$. Let $\widetilde{F}_0$ be the canonical lift of $F_0$ in $\overline{X}$. Remark that $\widetilde{F}_0 \cap G$ does not contain any edge. Let $\mathcal{Z}_1$ be the homothety class of the marked graph of groups $Z_1$ obtained from $X$ by collapsing  $\widetilde{F}_0 \cup (G-\{g\})$ and let $z_1$ be the corresponding vertex. Let $p_{x_1} \colon \overline{X} \to \overline{X}_1$ and $p_{z_1} \colon \overline{X} \to \overline{Z}_1$  be the natural projections. We claim that $x_1$ and $z_1$ are not adjacent. Indeed, suppose that $x_1$ and $z_1$ are adjacent. As $x_1$ and $z_1$ are distinct, we see that $|E\overline{X}_1| \neq |E\overline{Z}_1|$. Therefore, as $|E(\widetilde{F}_0 \cup (G-\{g\}))| \geq |EG|$, we see that $|E\overline{X}_1| > |E\overline{Z}_1|$, and a representative of $\mathcal{Z}_1$ is obtained from a representative of $\mathcal{X}_1$ by collapsing a forest. Let $X_1'$ be a representative of $\mathcal{X}_1$ obtained from $Z_1$ by blowing-up a forest. As $p_{z_1}(g)$ is an edge in $\overline{Z}_1$, Remark~\ref{free factor decomposition and canonical lift}~$(1)$ implies that there exists a unique edge $\widetilde{g}$ in $\overline{X}_1$ such that $\widetilde{g}$ induces the same free factor decomposition as $p_{z_1}(g)$ and $g$. But since $p_{x_1}(g)$ is a point, Remark~\ref{free factor decomposition and canonical lift}~$(1)$ implies that there is no edge in $\overline{X}_1'$ which induces the same free factor decomposition as $g$. So $x_1$ and $z_1$ are not adjacent and $z_1 \in A$.

In order to construct $z_2$, let $f \in EF_0$. Let $Y$ be the marked graph of groups obtained from $X_1$ by collapsing $F_0-\{f\}$. Let $p_0 \colon \overline{X}_1 \to \overline{Y}$ be the natural projection. Let $a_1$ and $a_2$ be the endpoints of $p_0(f)$. Let $X_1/F_0$ be the marked graph of groups obtained from $X_1$ by collapsing $F_0$. Since the homothety class of  $X_1/F_0$ is an element of $K_n$ (namely it is $\mathcal{X}_2$), one of the endpoints of $p_0(f)$ has trivial associated group. Suppose without loss of generality that $a_1$ has trivial associated group. In particular, $\operatorname{deg}(a_1) \geq 3$. Let $a_3$ and $a_4$ be two distinct vertices adjacent to $a_1$ other than $a_2$ and let $e$ be the edge between $a_1$ and $a_3$. Finally let $\mathcal{Z}_2$ be the homothety class of the marked graph of groups $Z_2$ obtained from $Y$ by collapsing $\{e\}$. Let $z_2$ be the corresponding vertex in $K_n$. Then, since $|EX_2|=|EZ_2|$ and since $X_2$ and $Z_2$ are obtained from $X$ by collapsing two distinct forests, we see that $z_2$ and $x_2$ are not adjacent in $K_n$. So $z_2 \in B$.

Let us prove that $z_1$ and $z_2$ are not adjacent in $\lk(x)$. Suppose towards a contradiction that $z_1$ and $z_2$ are adjacent. As $G$ contains at least one edge, we have that $Z_1$ is obtained from $X$ by collapsing $|F_0|+|G|-1$ edges. Moreover, $Z_2$ is obtained from $X$ by collapsing $|F_0|+|G|$ edges. This implies that the number of edges of a representative of $z_1$ is greater than the number of edges of a representative of $z_2$. Thus, there exists a representative of $z_1$ that collapses onto a representative of $z_2$. Let $p_{z_2} \colon \overline{X} \to \overline{Z}_2$ be the natural projection. Let $\widetilde{f} \in \widetilde{F}_0$ be the canonical lift of $f$ in $\overline{X}$. Since  $p_{z_2}(\widetilde{f})$ is an edge in $\overline{Z}_2$, Remark~\ref{free factor decomposition and canonical lift}~$(1)$ implies that there exists an edge in $\overline{Z}_2$ which induces the same free free factor decomposition as $\widetilde{f}$. But, as $p_{z_1}(\widetilde{f})$ is a point in $\overline{Z}_1$, Remark~\ref{free factor decomposition and canonical lift}~$(1)$ shows that there is no edge in $\overline{Z}_2$ that induces the same free factor decomposition as $\widetilde{f}$. Thus, $z_1$ and $z_2$ are not adjacent. 

This contradicts the fact that $\lk(x)=A \ast B$ is a join decomposition. This concludes the proof of the proposition.
\hfill\qedsymbol

\medskip

\begin{coro}\label{0star F star preserved}
Let $n \geq 4$ and $f \in \Aut(K_n)$. Then $f$ preserves the set of $\{0\}$-stars and the set of $F$-stars.
\end{coro}

\dem Let $\rho$ be a $\{0\}$-star. Since $\lk_-(\rho)$ has no edge and is of cardinal equal to $n$, Proposition~\ref{decomposition join} tells us that either $\lk_+(f(\rho))$ has no edge and its cardinal is equal to $n$, or $\lk_-(f(\rho))$ has no edge and its cardinal is equal to $n$. Since $n \geq 4$, Lemma~\ref{quasi maximal vertices lk_+=3}~$(1)$ tells us that the first case is not possible. So $\lk_-(f(\rho))$ has no edge and its cardinal is equal to $n$. Then Lemma~\ref{quasi minimal number of vertices}~$(2)$ shows that $f(\rho)$ is a $\{0\}$-star.

\medskip

Let $\rho'$ be an $F$-star. Then there exists a $\{0\}$-star $\rho$ such that $\rho' \in \lk_-(\rho)$. Therefore, $f(\rho') \in \lk_-(f(\rho))$. As $f(\rho)$ is a $\{0\}$-star and since the negative link of a $\{0\}$-star is composed of $F$-stars, we see that $f(\rho')$ is an $F$-star.
\hfill\qedsymbol

\bigskip

Thus, there exists a homomorphism $\Aut(K_n) \to \Aut(L_n)$ defined by restriction. We now prove that this homomorphism is in fact injective.

\begin{lem}\label{quasi minimaux et minimaux fixes}
Let $n \geq 4$. Let $f \in \Aut(K_n)$ be such that $f|_{O_n}=\mathrm{id}_{O_n}$ and $f|_{F_n}=\mathrm{id}_{F_n}$. Let $y=\mathcal{Y} \in VK_n$ be such that $\lk_-(v)$ is trivial. Then $f(y)=y$.
\end{lem}

\dem In order to prove Lemma~\ref{quasi minimaux et minimaux fixes}, we prove the following claim.

\medskip

\noindent{\bf Claim. } Let $0 \leq k \leq n-3$. Let $\mathcal{X}$ and $\mathcal{Y}$ be vertices of $K_n$. Let $X$ and $Y$ be representatives of $\mathcal{X}$ and $\mathcal{Y}$. We write $\overline{X}$ and $\overline{Y}$ for their underlying graphs. Suppose that $\mathcal{X}$ has a nontrivial negative link with no edge and that $\mathcal{Y}$ has a trivial negative link. 
If $\overline{X}$ has $k$ vertices with nontrivial associated group that are not leaves, and if $\overline{Y}$ has $k+1$ vertices with nontrivial associated group that are not leaves, then $f(\mathcal{X})=\mathcal{X}$ and $f(\mathcal{Y})=\mathcal{Y}$.

\medskip

Lemma~\ref{quasi minimaux et minimaux fixes} then follows from the claim because for every vertex $y \in VK_n$ with trivial negative link, there exists $k \in \{0,\ldots,n-3\}$ such that $y$ has a representative $Y$ whose underlying graph has exactly $k+1$ vertices with nontrivial associated group that are not leaves. 

We prove the claim by induction on $k$. When $k=0$, $\overline{X}$ has $n$ leaves, so by Lemma~\ref{quasi minimal number of vertices}~$(2)$, we have that $|\lk_-(v)|=n$. Thus, by Lemma~\ref{quasi minimal number of vertices}~$(2)$, we see that $\mathcal{X}$ is a $\{0\}$-star. Moreover, $\overline{Y}$ has $n-1$ leaves and $n$ vertices, so $\mathcal{Y}$ is an $F$-star. Thus, when $k=0$, the claim is a restatement of the fact that $f$ fixes the $\{0\}$-stars and the $F$-stars. 

Now suppose that the claim is true for some $0 \leq k \leq n-4$. Let $\mathcal{X}$ and $\mathcal{Y}$ be such that $\mathcal{X}$ has a nontrivial negative link with no edge and that $\mathcal{Y}$ has a trivial negative link. Let $X$ and $Y$ be representatives of $\mathcal{X}$ and $\mathcal{Y}$, and let $\overline{X}$ and $\overline{Y}$ their underlying graphs. Suppose that $\overline{X}$ has $k+1$ vertices with nontrivial associated group that are not leaves, and that $\overline{Y}$ has $k+2$ vertices with nontrivial associated group that are not leaves.

\bigskip

We start by showing that $f(\mathcal{X})=\mathcal{X}$. First, by Proposition~\ref{decomposition join}, the homothety class $f(\mathcal{X})$ has either a nontrivial negative link with no edge or a nontrivial positive link with no edge. 

\medskip

\noindent{\bf Claim. } The homothety class $\mathcal{X}$ cannot be sent by $f$ to a homothety class $z=\mathcal{Z}$ such that $\lk_+(z)$ has no edge. 

\medskip

\dem Suppose towards a contradiction that it is the case. By Lemma~\ref{quasi minimal number of vertices}, $|\lk_-(\mathcal{X})| \geq 3$, while by Lemma~\ref{quasi maximal vertices lk_+=3}~$(1)$, $|\lk_+(z)| \leq 3$. Thus, $|\lk_-(\mathcal{X})|=|\lk_+(z)|=3$. But then, Lemma~\ref{quasi maximal vertices lk_+=3}~$(1)$ implies that the underlying graph of any representative of $\mathcal{Z}$ has $n$ leaves. However, such a vertex $z$ is adjacent to $n$ $F$-stars whereas $\mathcal{X}$ is adjacent to at most one $F$-star. Indeed if $k+1=1$, the homothety class $\mathcal{X}$ is adjacent to exactly one $F$-star obtained from $X$ by collapsing the unique edge between the vertex with trivial associated group (the uniqueness of this vertex follows from Lemma~\ref{quasi minimal number of vertices}~$(1)$) and the non-leaf vertex with nontrivial associated group. If $k+1 \geq 2$, then $\mathcal{X}$ is not adjacent to an $F$-star because $\overline{X}$ has at least two vertices with nontrivial associated group that are not leaves, whereas any $F$-star has exactly one such vertex. As the set of $F$-stars is fixed by $f$, we get a contradiction.
\hfill\qedsymbol

\medskip

So $f(\mathcal{X})$ has a nontrivial negative link with no edge. Let $v$ be the unique vertex of $\overline{X}$ with trivial associated group given by Lemma~\ref{quasi minimal number of vertices}~$(1)$.

\medskip

\noindent{\bf Claim. } The underlying graph of any representative of $f(\mathcal{X})$ has exactly $n-k-1$ leaves. 

\medskip

\dem By the induction hypothesis, the automorphism $f$ fixes all vertices of $K_n$ whose negative link is nontrivial and has no edges and such that the underlying graph of any representative has at least $n-k$ leaves. Thus, the underlying graph of any representative of $f(\mathcal{X})$ has at most $n-k-1$ leaves. 

Now, suppose that $\mathcal{Z}$ is the homothety class of a marked graph of groups $Z$ whose underlying graph has at most $n-k-2$ leaves and such that $\lk_-(\mathcal{Z})$ is nontrivial and has no edge. Then $\lk_-(\mathcal{Z})$ does not contain any homothety class of marked graphs of groups whose underlying graph has $n-k-1$ leaves. But $\lk_-(\mathcal{X})$ contains one such homothety class, namely the homothety class of a marked graph of groups obtained from $X$ by collapsing an edge between $v$ and a vertex that is not a leaf. As $f$ fixes all vertices of $K_n$ with trivial negative link and such that the underlying graph of any representative has at least $n-k-1$ leaves and as $f(\lk_-(\mathcal{X}))=\lk_-(f(\mathcal{X}))$, it follows that $f(\mathcal{X}) \neq \mathcal{Z}$. Thus, the underlying graph of any representative of $f(\mathcal{X})$ has at least $n-k-1$ leaves. Therefore the underlying graph of any representative of $f(\mathcal{X})$ has exactly $n-k-1$ leaves.
\hfill\qedsymbol

\medskip

To prove that, in fact, $f(\mathcal{X})=\mathcal{X}$, we distinguish between two cases, according to the vertices adjacent to $v$. Note that, as $\overline{X}$ is connected, the vertex $v$ is adjacent to at least one vertex that is not a leaf.

\medskip

\noindent{\bf Case 1. } Suppose that $v$ is adjacent to at least two vertices $w_1$ and $w_2$ that are not leaves. 

\medskip

For $i \in \{1,2\}$, let $e_i$ be the edge between $v$ and $w_i$, and let $\mathcal{Y}_i$ be the homothety class of the marked graph of groups $Y_i$ obtained from $X$ by collapsing $e_i$. Then $\mathcal{Y}_1$ and $\mathcal{Y}_2$ are homothety classes of marked graphs of groups with trivial negative link and such that the underlying graphs of $Y_1$ and $Y_2$ have $k+1$ vertices with nontrivial associated group that are not leaves. By induction hypothesis, $f(\mathcal{Y}_1)=\mathcal{Y}_1$ and $f(\mathcal{Y}_2)=\mathcal{Y}_2$. Let $p_1 \colon X \to Y_1$ and $p_2 \colon X \to Y_2$ be the natural projections. In Case $1$, the fact that $f(\mathcal{X})=\mathcal{X}$ is a consequence of the following claim.

\medskip

\noindent{\bf Claim. } The homothety class $\mathcal{X}$ is the only vertex in $\lk(\mathcal{Y}_1) \cap \lk(\mathcal{Y}_2)$ whose negative link is nontrivial and has no edge.

\medskip

\dem Let $\mathcal{Z} \in \lk(\mathcal{Y}_1) \cap \lk(\mathcal{Y}_2)$ be such that $\lk_-(\mathcal{Z})$ is nontrivial and has no edge. Assume towards a contradiction that $\mathcal{Z} \neq \mathcal{X}$. As $\mathcal{Y}_1$ has trivial negative link, for all $\mathcal{Z}' \in VK_n$ such that $\mathcal{Z}' \in \lk(\mathcal{Y}_1)$, we have in fact $\mathcal{Z}' \in \lk_+(\mathcal{Y}_1)$. Thus, there exists a representative $Z$ of $\mathcal{Z}$ such that $Z$ is obtained from $Y_1$ by blowing-up a forest $F_0$. Let $\overline{Z}$ be the underlying graph of $Z$, and let $p_1^{Z} \colon \overline{Z} \to \overline{Y}_1$ be the natural projection.

We claim that there exists a unique edge in $F_0$. Indeed, otherwise there would exist two vertices in $\overline{Z}$ with trivial associated groups. As $\lk_-(\mathcal{Z})$ has no edge, this would contradict Lemma~\ref{quasi minimal number of vertices}~$(1)$. Thus, there exists a unique edge $f \in EF_0$.

Since $\mathcal{Z} \in \lk(\mathcal{Y}_2)$ and since $\lk_-(\mathcal{Z})$ is nontrivial and has no edge, Lemma~\ref{quasi minimal number of vertices}~$(1)$ implies that there exists an edge $g$ such that the homothety class of the marked graph of groups $Z/\{g\}$ obtained from $Z$ by collapsing $g$ is $\mathcal{Y}_2$. Let $p_2^{Z} \colon \overline{Z} \to \overline{Z/\{g\}}$ be the natural projection. By Remark~\ref{free factor decomposition and canonical lift}~$(1)$, there exists a unique edge $h \in E\overline{Z}$ such that $p_2^Z(h)$ induces the same free factor decomposition of $W_n$ as $p_2(e_1)$. But since $Z$ is a blow-up of $Y_1$ by an edge, and since $Y_1$ is obtained from $X$ by collapsing $e_1$, Lemma~\ref{partition link_-(x)}~$(2)$ implies that $p_1^Z(h)$ is reduced to a point. Therefore $f=h$ and $\overline{Z}$ is obtained from $Y_1$ by blowing-up the edge $e_1$. It follows that the graph $\overline{Z}$ is isomorphic to the graph $\overline{X}$. Thus, we can suppose that $\overline{X}=\overline{Z}$. We can also suppose, by Lemma~\ref{partition link_-(x)}~$(2)$, that $g=e_2$. As $v$ has trivial associated group, $\operatorname{deg}(v) \geq 3$. If $\mathcal{X} \neq \mathcal{Z}$, since both $X$ and $Z$ are obtained from $Y_1$ by blowing-up the edge $e_1$, there exist an integer $\ell \in \{0,1\}$ and a vertex $w_3 \in V\overline{X}$ distinct from $w_1$ and $w_2$ and adjacent to $v$ such that:

\begin{enumerate}
\item For $i \in \{1,2,3\}$, the preimage by the marking of $X$ of the generator of the group associated with $w_i$ is $x_i$~; 
\item The preimage by the marking of $Z$ of the generator of the group associated with $w_2$ is $x_1^{\ell}x_2x_1^{\ell}$ and the preimage by the marking of $Z$ of the generator of the group associated with $w_3$ is $x_1^{\ell +1}x_3x_1^{\ell +1}$.
\end{enumerate}

As $p_2(w_2)$ and $p_2(w_3)$ are in the same connected component of $p_2(\overline{X})-\{p_2(w_1)\}$, it follows that $p(e_1)$ and $p_2(e_1)$ induces distinct free factor decompositions of $W_n$. This contradicts the fact that $Z/\{e_1\} \in \mathcal{Y}_2$ by Remark~\ref{free factor decomposition and canonical lift}~$(1)$. The claim follows.
\hfill\qedsymbol

\medskip

\noindent{\bf Case 2. } Suppose that $v$ is adjacent to only one vertex $w$ that is not a leaf. 

\medskip

Let $e$ be the edge between $v$ and $w$ and let $\mathcal{Y}'$ be the homothety class of the marked graph of groups $Y'$ obtained from $X$ by collapsing $e$. Let $\overline{Y}'$ be the underlying graph of $Y'$. Let $p_{X} \colon \overline{X} \to \overline{Y}'$ be the natural projection. Then, as $\lk_-(\mathcal{Y}')$ is trivial and as $\mathcal{Y}'$ has a representative whose underlying graph has $n-k-1$ leaves, by induction hypothesis, we see that $f(\mathcal{Y}')=\mathcal{Y}'$. So $f(\mathcal{X}) \in \lk(\mathcal{Y}')$. Thus a representative $Z$ of $f(\mathcal{X})$ is obtained from $Y'$ by blowing-up a forest $F_0 \subseteq E\overline{Z}$. As $\lk_-(f(\mathcal{X}))$ has no edge, the forest $F_0$ contains a unique edge $e'$. Let $\overline{Z}$ be the underlying graph of $Z$, and $p_Z \colon \overline{Z} \to \overline{Y}'$ be the canonical projection. 

Suppose towards a contradiction that $f(\mathcal{X}) \neq \mathcal{X}$. By the claim above Case~$1$, the underlying graph of any representative of $f(\mathcal{X})$ has exactly $n-k-1$ leaves. Therefore none of the two endpoints of $e'$ is a leaf. Thus, as one of the endpoints of $e'$ has trivial associated group, there exists a vertex $a \in  V\overline{Y}'$ such that $\operatorname{deg}(a) \geq 3$ and such that $e'$ collapses onto $a$. As $\mathcal{Y}'$ has trivial negative link, we see that the group associated with $a$ is nontrivial. Let $y_i$ be the preimage by the marking of $Y'$ of the generator of the group associated with $a$. Let $\widetilde{a}$ be the lift of $a$ in $\overline{Z}$ such that $\widetilde{a}$  has nontrivial associated group. Then $y_i$ is the preimage by the marking of $Z$ of the generator of the group associated with $\widetilde{a}$. Let $\widetilde{b}$ be the endpoint of $e'$ distinct from $\widetilde{a}$ (see Figure~\ref{construction Z1'}). As $Z/\{e'\} \in \mathcal{Y}'$, the vertex $\widetilde{b}$ has trivial associated group and $\operatorname{deg}(\widetilde{b}) \geq 3$. Moreover, by the previous case, the vertex $\widetilde{b}$ cannot be adjacent to two vertices that are not leaves.

\medskip

Suppose first that $Z$ is not a blow-up of $Y'$ at $p_X(e)$. This implies that $p_X^{-1}(a)$ is a vertex.

\medskip

As $\widetilde{b}$ is not adjacent to two vertices that are not leaves, there exist two distinct leaves $w_1$ and $w_2$ of $\overline{Y}'$ adjacent to $a$ such that $w_1$ and $w_2$ have lifts $\widetilde{w}_1$ and $\widetilde{w}_2$ in $\overline{Z}$ with nontrivial associated group that are adjacent to $\widetilde{b}$. Let $y_j$ and $y_k$ be the preimages by the marking of $Y'$ of the groups associated with $w_1$ and $w_2$. Then there exist $\alpha_j,\alpha_k \in \{0,1\}$ such that  $y_i^{\alpha_j}y_jy_i^{\alpha_j}$ and $y_i^{\alpha_k}y_ky_i^{\alpha_k}$ are the preimages by the marking of $Z$ of the groups associated with $\widetilde{w}_1$ and $\widetilde{w}_2$.

Let $\mathcal{Z}_1'$ be the homothety class of the marked graph of groups $Z_1'$ defined as follows (see Figure~\ref{construction Z1'}):

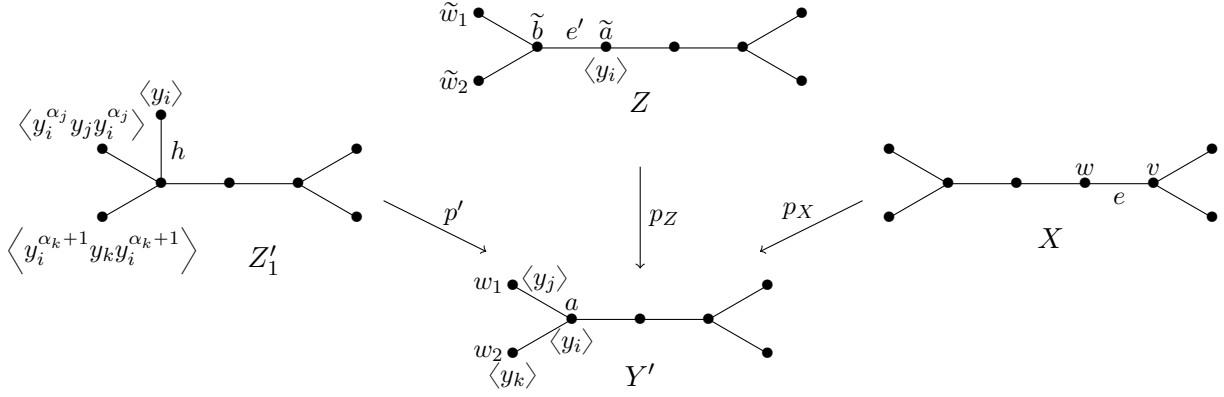
\begin{figure}
\hspace{-1cm}
\begin{tikzpicture}[scale=0.90]
\draw (0.5-7,0) -- (1.5-7,0);
\draw (0.5-7,0) node {$\bullet$};
\draw (1.5-7,0) node {$\bullet$};
\draw (1.5-7,0) -- (2.5-7,0);
\draw (2.5-7,0) node {$\bullet$};
\draw (0.5-7,0) -- (-0.36-7,0.5);
\draw (-0.36-7,0.5) node {$\bullet$};
\draw (0.5-7,0) -- (-0.36-7,-0.5);
\draw (-0.36-7,-0.5) node {$\bullet$};
\draw (0.5-7,0) -- (0.5-7,1);
\draw (0.5-7,1) node {$\bullet$};
\draw (2.5-7,0) -- (3.36-7,0.5);
\draw (3.36-7,0.5) node {$\bullet$};
\draw (2.5-7,0) -- (3.36-7,-0.5);
\draw (3.36-7,-0.5) node {$\bullet$};
\draw (0.5-7,1) node[above,scale= 0.90] {$\left\langle y_i \right\rangle$};
\draw (0.5-7,0.5) node[right,scale= 0.90] {$h$};
\draw (-0.36-7,-0.5) node[below,scale= 0.90] {$\left\langle y_i^{\alpha_{k}+1}y_ky_i^{\alpha_k+1} \right\rangle$};
\draw (-0.66-7,0.5) node[above,scale= 0.90] {$\left\langle y_i^{\alpha_j}y_jy_i^{\alpha_j} \right\rangle$};
\draw (2-7,-0.75) node[below, scale=1] {$Z_1'$};

\draw (4.5-7,-0.75) node[above right, scale=0.9] {$p'$};
\draw [->] (3.75-7,-0.25) -- (5.25-7,-1);

\draw (-1,2) -- (0,2);
\draw (0,2) node {$\bullet$};
\draw (-1,2) node {$\bullet$};
\draw (-0,2) -- (1,2);
\draw (1,2) node {$\bullet$};
\draw (1,2) -- (2,2);
\draw (2,2) node {$\bullet$};
\draw (-1,2) -- (-1.86,2.5);
\draw (-1.86,2.5) node {$\bullet$};
\draw (-1,2) -- (-1.86,1.5);
\draw (-1.86,1.5) node {$\bullet$};
\draw (2,2) -- (2.86,2.5);
\draw (2.86,2.5) node {$\bullet$};
\draw (2,2) -- (2.86,1.5);
\draw (2.86,1.5) node {$\bullet$};
\draw (-1,2) node[above,scale= 0.90] {$\widetilde{b}$};
\draw (0,2) node[above,scale= 0.90] {$\widetilde{a}$};
\draw (-0.45,2) node[above,scale= 0.90] {$e'$};
\draw (-0,2) node[below,scale= 0.90] {$\left\langle y_i \right\rangle$};
\draw (-1.86,2.5) node[left,scale= 0.90] {$\widetilde{w}_1$};
\draw (-1.86,1.5) node[left,scale= 0.90] {$\widetilde{w}_2$};
\draw (0.5,1.5) node[below, scale=1] {$Z$};

\draw [->] (0.5, 0.25) -- (0.5,-1.25);
\draw (0.5,-0.5) node[right, scale=0.9] {$p_{Z}$};

\draw (-0.5,-2) -- (0.5,-2);
\draw (-0.5,-2) node {$\bullet$};
\draw (0.5,-2) node {$\bullet$};
\draw (0.5,-2) -- (1.5,-2);
\draw (1.5,-2) node {$\bullet$};
\draw (-0.5,-2) -- (-1.36,-2.5);
\draw (-1.36,-2.5) node {$\bullet$};
\draw (-0.5,-2) -- (-1.36,-1.5);
\draw (-1.36,-1.5) node {$\bullet$};
\draw (1.5,-2) -- (2.36,-1.5);
\draw (2.36,-1.5) node {$\bullet$};
\draw (1.5,-2) -- (2.36,-2.5);
\draw (2.36,-2.5) node {$\bullet$};
\draw (-0.5,-2) node[above,scale= 0.90] {$a$};
\draw (-0.5,-2) node[below,scale= 0.90] {$\left\langle y_i \right\rangle$};
\draw (-1.36,-1.5) node[left,scale= 0.90] {$w_1$};
\draw (-1.36,-1.4) node[right,scale= 0.90] {$\left\langle y_j \right\rangle$};
\draw (-1.36,-2.5) node[left,scale= 0.90] {$w_2$};
\draw (-1.36,-2.5) node[below,scale= 0.90] {$\left\langle y_k \right\rangle$};
\draw (0.5,-2.5) node[below, scale=1] {$Y'$};

\draw [<-] (2.25,-1) -- (3.75,-0.25);
\draw (3.2,-0.65) node[above left, scale=0.9] {$p_{X}$};

\draw (-1.5+6.5,0) -- (-0.5+6.5,0);
\draw (-0.5+6.5,0) node {$\bullet$};
\draw (-1.5+6.5,0) node {$\bullet$};
\draw (-0.5+6.5,0) -- (0.5+6.5,0);
\draw (0.5+6.5,0) node {$\bullet$};
\draw (0.5+6.5,0) -- (1.5+6.5,0);
\draw (1.5+6.5,0) node {$\bullet$};
\draw (-1.5+6.5,0) -- (-2.36+6.5,0.5);
\draw (-2.36+6.5,0.5) node {$\bullet$};
\draw (-1.5+6.5,0) -- (-2.36+6.5,-0.5);
\draw (-2.36+6.5,-0.5) node {$\bullet$};
\draw (1.5+6.5,0) -- (2.36+6.5,0.5);
\draw (2.36+6.5,0.5) node {$\bullet$};
\draw (1.5+6.5,0) -- (2.36+6.5,-0.5);
\draw (2.36+6.5,-0.5) node {$\bullet$};
\draw (1.5+6.5,0) node[above,scale= 0.90] {$v$};
\draw (1+6.5,0) node[below,scale= 0.90] {$e$};
\draw (0.5+6.5,0) node[above,scale= 0.90] {$w$};
\draw (0+6.5,-0.5) node[below, scale=1] {$X$};

\end{tikzpicture}
\caption{The construction of $Z_1'$ in Lemma~\ref{quasi minimaux et minimaux fixes} when $p_X(e) \neq a$.}\label{construction Z1'}
\end{figure} 

$\bullet$ The underlying graph of $Z_1'$ is obtained from $\overline{Y}'$ by pulling-up an edge $h$ at $a$ so that one of the two endpoints of $h$ is a leaf. Let $p' \colon \overline{Z}_1' \to \overline{Y}'$ be the projection. Let $x$ be a vertex of the underlying graph of $Z_1'$. Remark that, as $w_2$ is a leaf, $p'^{-1}(w_2)$ is a leaf.

$\bullet$ If $x$ is distinct from $p'^{-1}(w_2)$ and is such that $p'(x) \neq a$, then the group associated with $x$ in $Z_1'$ is the same one as the group associated with $p'(x)$. 

$\bullet$ If $p'(x)=a$ and if $x$ is a leaf, then the group associated with $x$ is the same one as the group associated with $a$. 

$\bullet$ If $p'(x)=a$ and if $x$ is not a leaf, then $x$ has trivial associated group. 

$\bullet$ Finally, the preimage by the marking of $Z_1'$ of the generator of the group associated with $p'^{-1}(w_2)$ is $y_i^{\alpha_k +1}y_ky_i^{\alpha_k +1}$.

\noindent By the induction hypothesis, $\mathcal{Z}_1'$ is fixed by $f$. What is more, $d_{\lk(\mathcal{Y}')}(\mathcal{X},\mathcal{Z}_1')=2$. Indeed, a common refinement of $X$ and $Z_1'$ is obtained from $X$ by pulling-up the edge $h$ at $p_X^{-1}(a)$ (this is possible since $p_X(e) \neq a$).

\medskip

\noindent{\bf Claim. } In $\lk(\mathcal{Y}')$, we have $d_{\lk(\mathcal{Y}')}(\mathcal{Z},\mathcal{Z}_1')>2$.

\medskip

\dem Since both $\mathcal{Z}$ and $\mathcal{Z}_1'$ have nontrivial negative link with no edge, we see by Lemma~\ref{quasi minimal number of vertices}~$(1)$ that $|V\overline{Z}|=|V\overline{Z}_1'|$. As both $\overline{Z}$ and $\overline{Z}_1'$ are trees, we have $|E\overline{Z}|=|E\overline{Z}_1'|$. Thus, as $\mathcal{Z} \neq \mathcal{Z}_1'$, we have \mbox{$d_{\lk(\mathcal{Y}')}(\mathcal{Z},\mathcal{Z}_1')>1$}. As $\mathcal{Y}'$ has trivial negative link, the only way $d_{\lk(\mathcal{Y}')}(\mathcal{Z},\mathcal{Z}_1')=2$ is that $Z$ and $Z_1'$ have a common refinement. Let $z$ be the leaf of $\overline{Z}_1'$ such that $p'(z)=a$. Then $p'^{-1}(w_1)$ and $p'^{-1}(w_2)$ are in the same connected component of $\overline{Z}_1'-\{z\}$. Let $Z_1^{(2)}$ be a refinement of $Z_1$, and let $y_1^{(2)},\ldots,y_n^{(2)}$ be the preimages by the marking of $Z_1^{(2)}$ of the generators of the nontrivial vertex groups of $Z_1^{(2)}$. Suppose that, for all $m \in \{1,\ldots,n\}$, there exists $\alpha_m \in \{0,1\}$ such that $y_m^{(2)}=y_i^{\alpha_m}y_my_i^{\alpha_m}$ and that there exist $m_0$ and $m_1$ such that $\alpha_{m_0}=0$ and $\alpha_{m_1}=1$. Since the preimage by the marking of the generator of the group associated with $z$ is $y_i$, we see that $Z_1^{(2)}$ is obtained from $Z_1'$ by blowing-up a forest and applying a twist at an edge whose terminal point is $z$. As a consequence, since $p'^{-1}(w_1)$ and $p'^{-1}(w_2)$ are in the same connected component of $\overline{Z}_1'-\{z\}$, there does not exist a refinement of $Z_1'$ such that the preimages by the marking of the generator of the group associated with lifts of $p'^{-1}(w_1)$ and $p'^{-1}(w_2)$ are respectively $y_i^{\alpha_j}y_jy_i^{\alpha_j}$ and $y_i^{\alpha_k}y_ky_i^{\alpha_k}$. Thus, $Z$ and $Z_1'$ do not have any common refinement and $d_{\lk(\mathcal{Y}')}(\mathcal{Z},\mathcal{Z}_1')>2$.
\hfill\qedsymbol

\medskip

Since $d_{\lk(f(\mathcal{Y}'))}(f(\mathcal{X}),f(\mathcal{Z}_1'))=d_{\lk(\mathcal{Y}')}(\mathcal{X},\mathcal{Z}_1')=2$, the last claim implies that $f(\mathcal{X})=\mathcal{X}$ when $p_X(e) \neq a$.

\bigskip

Suppose now that $p_X(e)=a$. Then, as $|\lk_-(\mathcal{X})|=|\lk_-(f(\mathcal{X}))|$, and as $\overline{X}$ and $\overline{Z}$ both have a unique vertex with trivial associated group by Lemma~\ref{quasi minimal number of vertices}~$(1)$ (namely $v$ and $\widetilde{b}$), we see that $\operatorname{deg}(v)=\operatorname{deg}(\widetilde{b})=m$. Moreover, both $v$ and $\widetilde{b}$ are adjacent to a unique vertex that is not a leaf (namely $w$ and $\widetilde{a}$). Thus, both $v$ and $\widetilde{b}$ are adjacent to exactly $m-1$ leaves. Note that, as $v$ and $\widetilde{b}$ have trivial associated group, $m-1 \geq 2$. Let $v_1,\ldots,v_{m-1}$ be the leaves of $\overline{X}$ adjacent to $v$, and let $\widetilde{b}_1,\ldots,\widetilde{b}_{m-1}$ be the leaves of $\overline{Z}$ adjacent to $\widetilde{b}$. For $j \in \{1,\ldots,m-1\}$, let $y_j^X$ be the preimage by the marking of the generator of the group associated with $v_j$ and let $y_j^Z$ be the preimage by the marking of the generator of the group associated with $\widetilde{b}_j$. As we suppose that $\mathcal{X} \neq \mathcal{Z}$, up to reordering and composing by an inner automorphism, one of the following holds.

\medskip

\noindent{$(i)$ } There exist $j,k \in \{1,\ldots,m-1\}$ distinct such that $y_j^X=y_iy_{j}^Zy_i$ and $y_{k}^X=y_{k}^Z$ (see Figure~\ref{construction Z2'}).

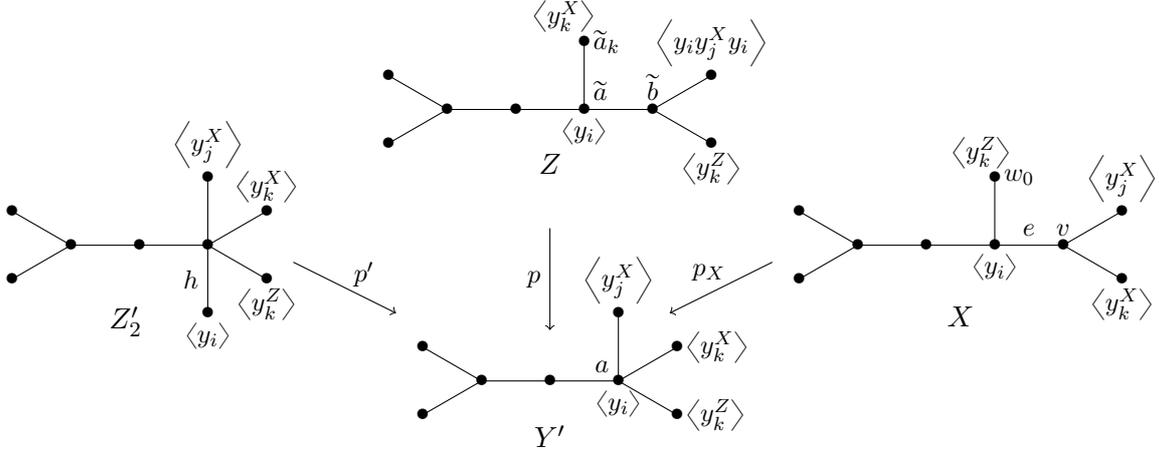
\begin{figure}
\hspace*{-0.5cm}
\begin{tikzpicture}[scale=0.90]
\draw (0.5-7,0) -- (1.5-7,0);
\draw (0.5-7,0) node {$\bullet$};
\draw (1.5-7,0) node {$\bullet$};
\draw (1.5-7,0) -- (2.5-7,0);
\draw (2.5-7,0) node {$\bullet$};
\draw (0.5-7,0) -- (-0.36-7,0.5);
\draw (-0.36-7,0.5) node {$\bullet$};
\draw (0.5-7,0) -- (-0.36-7,-0.5);
\draw (-0.36-7,-0.5) node {$\bullet$};
\draw (2.5-7,0) -- (2.5-7,1);
\draw (2.5-7,1) node {$\bullet$};
\draw (2.5-7,0) -- (2.5-7,-1);
\draw (2.5-7,-1) node {$\bullet$};
\draw (2.5-7,0) -- (3.36-7,0.5);
\draw (3.36-7,0.5) node {$\bullet$};
\draw (2.5-7,0) -- (3.36-7,-0.5);
\draw (3.36-7,-0.5) node {$\bullet$};
\draw (2.5-7,-1) node[below,scale= 0.90] {$\left\langle y_i \right\rangle$};
\draw (2.5-7,-0.5) node[left,scale= 0.90] {$h$};
\draw (2.5-7,1) node[above,scale= 0.90] {$\left\langle y_j^X \right\rangle$};
\draw (3.36-7,0.5) node[above,scale= 0.90] {$\left\langle y_{k}^X \right\rangle$};
\draw (3.36-7,-0.5) node[below,scale= 0.90] {$\left\langle y_{k}^Z \right\rangle$};
\draw (1.3-7,-0.75) node[below, scale=1] {$Z_2'$};

\draw [->] (3.75-7,-0.25) -- (5.25-7,-1);
\draw (4.5-7,-0.75) node[above right, scale=0.9] {$p'$};

\draw (-1,2) -- (0,2);
\draw (0,2) node {$\bullet$};
\draw (-1,2) node {$\bullet$};
\draw (-0,2) -- (1,2);
\draw (1,2) node {$\bullet$};
\draw (1,2) -- (2,2);
\draw (2,2) node {$\bullet$};
\draw (-1,2) -- (-1.86,2.5);
\draw (-1.86,2.5) node {$\bullet$};
\draw (-1,2) -- (-1.86,1.5);
\draw (-1.86,1.5) node {$\bullet$};
\draw (2,2) -- (2.86,2.5);
\draw (2.86,2.5) node {$\bullet$};
\draw (2,2) -- (2.86,1.5);
\draw (2.86,1.5) node {$\bullet$};
\draw (1,2) -- (1,3);
\draw (1,3) node {$\bullet$};
\draw (2,2) node[above,scale= 0.90] {$\widetilde{b}$};
\draw (1,2) node[below,scale= 0.90] {$\left\langle y_i \right\rangle$};
\draw (1,2) node[above right,scale= 0.90] {$\widetilde{a}$};
\draw (0.7,3) node[above,scale= 0.90] {$\left\langle y_k^X \right\rangle$};
\draw (2.86,1.5) node[below,scale= 0.90] {$\left\langle y_k^Z \right\rangle$};
\draw (1,3) node[right,scale= 0.90] {$\widetilde{a}_k$};
\draw (2.86,2.5) node[above,scale= 0.90] {$\left\langle y_iy_{j}^Xy_i \right\rangle$};
\draw (0.5,1.5) node[below, scale=1] {$Z$};

\draw [->] (0.5, 0.25) -- (0.5,-1.25);
\draw (0.5,-0.5) node[left, scale=0.9] {$p$};

\draw (-0.5,-2) -- (0.5,-2);
\draw (-0.5,-2) node {$\bullet$};
\draw (0.5,-2) node {$\bullet$};
\draw (0.5,-2) -- (1.5,-2);
\draw (1.5,-2) node {$\bullet$};
\draw (-0.5,-2) -- (-1.36,-2.5);
\draw (-1.36,-2.5) node {$\bullet$};
\draw (-0.5,-2) -- (-1.36,-1.5);
\draw (-1.36,-1.5) node {$\bullet$};
\draw (1.5,-2) -- (2.36,-1.5);
\draw (2.36,-1.5) node {$\bullet$};
\draw (1.5,-2) -- (2.36,-2.5);
\draw (2.36,-2.5) node {$\bullet$};
\draw (1.5,-2) -- (1.5,-1);
\draw (1.5,-1) node {$\bullet$};
\draw (1.5,-2) node[below,scale= 0.90] {$\left\langle y_i \right\rangle$};
\draw (1.5,-2) node[above left,scale= 0.90] {$a$};
\draw (1.5,-1) node[above,scale= 0.90] {$\left\langle y_j^X \right\rangle$};
\draw (2.36,-1.5) node[right,scale= 0.90] {$\left\langle y_{k}^X \right\rangle$};
\draw (2.36,-2.5) node[right,scale= 0.90] {$\left\langle y_{k}^Z \right\rangle$};
\draw (0.5,-2.5) node[below, scale=1] {$Y'$};

\draw [<-] (2.25,-1) -- (3.75,-0.25);
\draw (3.2,-0.65) node[above left, scale=0.9] {$p_{X}$};

\draw (-1.5+6.5,0) -- (-0.5+6.5,0);
\draw (-0.5+6.5,0) node {$\bullet$};
\draw (-1.5+6.5,0) node {$\bullet$};
\draw (-0.5+6.5,0) -- (0.5+6.5,0);
\draw (0.5+6.5,0) node {$\bullet$};
\draw (0.5+6.5,0) -- (1.5+6.5,0);
\draw (1.5+6.5,0) node {$\bullet$};
\draw (-1.5+6.5,0) -- (-2.36+6.5,0.5);
\draw (-2.36+6.5,0.5) node {$\bullet$};
\draw (-1.5+6.5,0) -- (-2.36+6.5,-0.5);
\draw (-2.36+6.5,-0.5) node {$\bullet$};
\draw (1.5+6.5,0) -- (2.36+6.5,0.5);
\draw (2.36+6.5,0.5) node {$\bullet$};
\draw (1.5+6.5,0) -- (2.36+6.5,-0.5);
\draw (2.36+6.5,-0.5) node {$\bullet$};
\draw (1.5+6.5,0) node[above,scale= 0.90] {$v$};
\draw (1+6.5,0) node[above,scale= 0.90] {$e$};
\draw (0.5+6.5,0) node[below,scale= 0.90] {$\left\langle y_i \right\rangle$};
\draw (0.5+6.5,1) node[right,scale= 0.90] {$w_0$};
\draw (0.3+6.5,1) node[above,scale= 0.90] {$\left\langle y_k^Z \right\rangle$};
\draw (2.36+6.5,-0.5) node[below,scale= 0.90] {$\left\langle y_{k}^X \right\rangle$};
\draw (0.5+6.5,0) -- (0.5+6.5,1);
\draw (0.5+6.5,1) node {$\bullet$};
\draw (2.36+6.5,0.5) node[above,scale= 0.90] {$\left\langle y_j^X \right\rangle$};
\draw (0+6.5,-0.75) node[below, scale=1] {$X$};

\end{tikzpicture}
\caption{The construction of $Z_2'$ in Lemma~\ref{quasi minimaux et minimaux fixes} in Case~$(ii)$.}\label{construction Z2' case with a lift}
\end{figure}

\noindent{$(ii)$ } There exist $j,k \in \{1,\ldots,m-1\}$ distinct and a leaf $\widetilde{a}_{k}$ adjacent to $\widetilde{a}$ such that $y_{j}^X=y_iy_{j}^Zy_i$ and such that the preimage by the marking of the generator of the group associated with $\widetilde{a}_k$ is $y_k^X$. Moreover, there exists a leaf $w_0$ in $\overline{X}$ adjacent to $w$ such that the preimage by the marking of the generator of the group associated with $w_0$ is $y_k^Z$ (see Figure~\ref{construction Z2' case with a lift}).

\noindent{$(iii)$ } There exist $j,\ell \in \{1,\ldots,m-1\}$ distinct and a leaf $\widetilde{a}_j$ adjacent to $\widetilde{a}$ such that $y_{\ell}^X=y_{\ell}^Z$ and such that the preimage by the marking of the generator of the group associated with $\widetilde{a}_j$ is $y_j^X$ (see Figure~\ref{construction Z3'}).

\noindent{$(iv)$ } For all $j \in \{1,\ldots,m-1\}$, there exists a leaf $\widetilde{a}_j$ adjacent to $\widetilde{a}$ such that the preimage by the marking of the generator of the group associated with $\widetilde{a}_j$ is $y_j^X$.

\medskip

\noindent We then distinguish two cases, according to whether $y_{j}^X=y_iy_{j}^Zy_i$ or not.

\bigskip

Suppose first that $y_j^Z=y_iy_j^Xy_i$ (Cases~$(i)$ and $(ii)$). Let $\mathcal{Z}_2'$ be the homothety class of the marked graph of groups $Z_2'$ defined as follows (see Figures~\ref{construction Z2'} and \ref{construction Z2' case with a lift}):

\begin{figure}
\hspace{-0.5cm}
\begin{tikzpicture}[scale=0.90]
\draw (0.5-7,0) -- (1.5-7,0);
\draw (0.5-7,0) node {$\bullet$};
\draw (1.5-7,0) node {$\bullet$};
\draw (1.5-7,0) -- (2.5-7,0);
\draw (2.5-7,0) node {$\bullet$};
\draw (0.5-7,0) -- (-0.36-7,0.5);
\draw (-0.36-7,0.5) node {$\bullet$};
\draw (0.5-7,0) -- (-0.36-7,-0.5);
\draw (-0.36-7,-0.5) node {$\bullet$};
\draw (2.5-7,0) -- (2.5-7,1);
\draw (2.5-7,1) node {$\bullet$};
\draw (2.5-7,0) -- (3.36-7,0.5);
\draw (3.36-7,0.5) node {$\bullet$};
\draw (2.5-7,0) -- (3.36-7,-0.5);
\draw (3.36-7,-0.5) node {$\bullet$};
\draw (2.5-7,1) node[above,scale= 0.90] {$\left\langle y_i \right\rangle$};
\draw (3.36-7,0.5) node[above,scale= 0.90] {$\left\langle y_j^X \right\rangle$};
\draw (3.36-7,-0.5) node[below,scale= 0.90] {$\left\langle y_k^X \right\rangle$};
\draw (1.5-7,-0.5) node[below, scale=1] {$Z_2'$};

\draw [->](3.75-7,-0.25) -- (5.25-7,-1);
\draw (4.5-7,-0.75) node[above right, scale=0.9] {$p'$};

\draw (-1,2) -- (0,2);
\draw (0,2) node {$\bullet$};
\draw (-1,2) node {$\bullet$};
\draw (-0,2) -- (1,2);
\draw (1,2) node {$\bullet$};
\draw (1,2) -- (2,2);
\draw (2,2) node {$\bullet$};
\draw (-1,2) -- (-1.86,2.5);
\draw (-1.86,2.5) node {$\bullet$};
\draw (-1,2) -- (-1.86,1.5);
\draw (-1.86,1.5) node {$\bullet$};
\draw (2,2) -- (2.86,2.5);
\draw (2.86,2.5) node {$\bullet$};
\draw (2,2) -- (2.86,1.5);
\draw (2.86,1.5) node {$\bullet$};
\draw (2,2) node[above,scale= 0.90] {$\widetilde{b}$};
\draw (1,2) node[above,scale= 0.90] {$\widetilde{a}$};
\draw (1.5,2) node[above,scale= 0.90] {$e'$};
\draw (1,2) node[below,scale= 0.90] {$\left\langle y_i \right\rangle$};
\draw (2.86,2.5) node[above,scale= 0.90] {$\left\langle y_iy_j^Xy_i \right\rangle$};
\draw (2.86,2.5) node[right,scale= 0.90] {$\widetilde{b}_j$};
\draw (2.86,1.5) node[below,scale= 0.90] {$\left\langle y_k^X \right\rangle$};
\draw (2.86,1.7) node[right,scale= 0.90] {$\widetilde{b}_k$};
\draw (0.5,1.5) node[below, scale=1] {$Z$};

\draw [->] (0.5, 0.25) -- (0.5,-1.25);
\draw (0.5,-0.5) node[right, scale=0.9] {$p$};

\draw (-0.5,-2) -- (0.5,-2);
\draw (-0.5,-2) node {$\bullet$};
\draw (0.5,-2) node {$\bullet$};
\draw (0.5,-2) -- (1.5,-2);
\draw (1.5,-2) node {$\bullet$};
\draw (-0.5,-2) -- (-1.36,-2.5);
\draw (-1.36,-2.5) node {$\bullet$};
\draw (-0.5,-2) -- (-1.36,-1.5);
\draw (-1.36,-1.5) node {$\bullet$};
\draw (1.5,-2) -- (2.36,-1.5);
\draw (2.36,-1.5) node {$\bullet$};
\draw (1.5,-2) -- (2.36,-2.5);
\draw (2.36,-2.5) node {$\bullet$};
\draw (1.5,-2) node[above,scale= 0.90] {$a$};
\draw (1.5,-2) node[below,scale= 0.90] {$\left\langle y_i \right\rangle$};

\draw (0.5,-2.5) node[below, scale=1] {$Y'$};

\draw [<-] (2.25,-1) -- (3.75,-0.25);
\draw (3.2,-0.65) node[above left, scale=0.9] {$p_{X}$};

\draw (-1.5+6.5,0) -- (-0.5+6.5,0);
\draw (-0.5+6.5,0) node {$\bullet$};
\draw (-1.5+6.5,0) node {$\bullet$};
\draw (-0.5+6.5,0) -- (0.5+6.5,0);
\draw (0.5+6.5,0) node {$\bullet$};
\draw (0.5+6.5,0) -- (1.5+6.5,0);
\draw (1.5+6.5,0) node {$\bullet$};
\draw (-1.5+6.5,0) -- (-2.36+6.5,0.5);
\draw (-2.36+6.5,0.5) node {$\bullet$};
\draw (-1.5+6.5,0) -- (-2.36+6.5,-0.5);
\draw (-2.36+6.5,-0.5) node {$\bullet$};
\draw (1.5+6.5,0) -- (2.36+6.5,0.5);
\draw (2.36+6.5,0.5) node {$\bullet$};
\draw (1.5+6.5,0) -- (2.36+6.5,-0.5);
\draw (2.36+6.5,-0.5) node {$\bullet$};
\draw (1.5+6.5,0) node[above,scale= 0.90] {$v$};
\draw (1+6.5,0) node[below,scale= 0.90] {$e$};
\draw (0.5+6.5,0) node[above,scale= 0.90] {$w$};
\draw (2.25+6.5,0.5) node[above,scale= 0.90] {$\left\langle y_j^X \right\rangle$};
\draw (2.46+6.5,0.5) node[right,scale= 0.90] {$v_j$};
\draw (2.25+6.5,-0.5) node[below,scale= 0.90] {$\left\langle y_k^X \right\rangle$};
\draw (2.46+6.5,-0.5) node[right,scale= 0.90] {$v_k$};
\draw (0+6.5,-0.5) node[below, scale=1] {$X$};

\end{tikzpicture}
\caption{The construction of $Z_2'$ in Lemma~\ref{quasi minimaux et minimaux fixes} in Case~$(i)$.}\label{construction Z2'}
\end{figure}
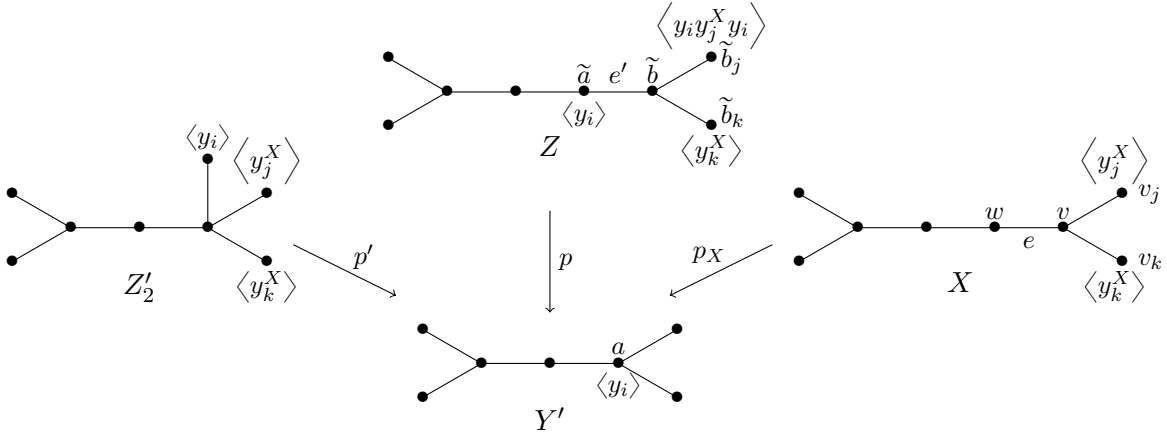 

$\bullet$ The underlying graph $\overline{Z}_2'$ of $Z_2'$ is obtained from $\overline{Y}'$ by blowing-up an edge $h$ at $a$ so that one of the two endpoints of $h$ is a leaf. Let $p' \colon \overline{Z}_2' \to \overline{Y}'$ be the projection. Let $x$ be a vertex of the underlying graph of $Z_2'$. Let $y$ be the preimage by the marking of the generator of the group associated with $p'(x)$. 

$\bullet$ If $p'(x) \neq a$ or if $p'(x)=a$ and $x$ is a leaf, then the preimage by the marking of $Z_2'$ of the generator of the group associated with $x$ is $y$.

$\bullet$ If $p'(x)=a$ and $x$ is not a leaf, then $x$ has trivial associated group.

\noindent By the induction hypothesis, as $\overline{Z}_2'$ has $k$ vertices with nontrivial associated vertex groups that are not leaves, the homothety class $\mathcal{Z}_2'$ is fixed by $f$. What is more, $d_{\lk(\mathcal{Y}')}(\mathcal{X},\mathcal{Z}_2')=2$. Indeed, a common refinement of $Z$ and $Z_2'$ is obtained from $X$ by blowing-up the edge $h$ at the vertex $w$ of $p_X^{-1}(a)$ with nontrivial associated group.

\medskip

\noindent{\bf Claim. } In $\lk(\mathcal{Y}')$, we have $d_{\lk(\mathcal{Y}')}(\mathcal{Z},\mathcal{Z}_2')>2$.

\medskip

\dem Since both $\mathcal{Z}$ and $\mathcal{Z}_2'$ have nontrivial negative link with no edge, we see by Lemma~\ref{quasi minimal number of vertices}~$(1)$ that $|V\overline{Z}|=|V\overline{Z}_2'|$. As both $\overline{Z}$ and $\overline{Z}_2'$ are trees, we have $|E\overline{Z}|=|E\overline{Z}_2'|$. Thus, as $\mathcal{Z} \neq \mathcal{Z}_2'$, we have \mbox{$d_{\lk(\mathcal{Y}')}(\mathcal{Z},\mathcal{Z}_2')>1$}. As $\mathcal{Y}'$ has trivial negative link, the only way $d_{\lk(\mathcal{Y}')}(\mathcal{Z},\mathcal{Z}_2')=2$ is that $Z$ and $Z_2'$ have a common refinement. Let $z$ be the leaf of $\overline{Z}_2'$ such that $p'(z)=a$. Then the preimage by the marking of the group associated with $z$ is $\left\langle y_i \right\rangle$.

Let $Z_2^{(2)}$ be a refinement of $Z_2'$, let $\overline{Z}_2^{(2)}$ be its underlying graph and let $y_1^{(2)},\ldots,y_n^{(2)}$ be the preimages by the marking of $Z_2^{(2)}$ of the generators of the nontrivial vertex groups of $Z_2^{(2)}$. Since both $Z$ and $Z_2'$ are obtained from $Y'$ by blowing-up an edge at $a$ while applying a twist around an edge adjacent to $a$, a potential common refinement of $Z$ and $Z_2'$ is obtained from $Y'$ by blowing-up a forest while applying a twist around an edge adjacent to $a$. Thus, we may assume that, for all $m \in \{1,\ldots,n\}$, there exists $\alpha_m \in \{0,1\}$ such that $y_m^{(2)}=y_i^{\alpha_m}y_my_i^{\alpha_m}$.

Suppose first that $y_k^X=y_k^Z$ (Case~$(i)$). Let $\widetilde{v}_j$ and $\widetilde{v}_k$ be the lifts in $\overline{Z}_2^{(2)}$ of respectively $p'^{-1}(p_X(v_j))$ and $p'^{-1}(p_X(v_k))$ with nontrivial associated group. Since $p'^{-1}(p_X(v_j))$ and $p'^{-1}(p_X(v_{k}))$ are contained in the same connected component of $\overline{Z}_2'-\{z\}$, there exists $\alpha \in \{0,1\}$ such that the preimages by the marking of $Z_2^{(2)}$ of the generators of the groups associated with $\widetilde{v}_j$ and $\widetilde{v}_k$ are respectively $y_i^{\alpha}y_j^Xy_i^{\alpha}$ and $y_i^{\alpha}y_k^Xy_i^{\alpha}$. As a consequence, there does not exist a refinement of $Z_2'$ such that the preimages by the marking of the generators of the groups associated with the lifts of $p'^{-1}(p_X(v_j))$ and $p'^{-1}(p_X(v_k))$ with nontrivial associated groups are respectively $y_iy_j^Xy_i$ and $y_k^X$. Thus, $Z$ and $Z_2'$ do not have any common refinement and $d_{\lk(\mathcal{Y}')}(\mathcal{Z},\mathcal{Z}_2')>2$.

Suppose now that there exists a leaf $w_0$ in $\overline{X}$ adjacent to $w$ such that the preimage by the marking of the generator of the group associated with $w_0$ is $y_k^Z$ (Case~$(ii)$). Let $\widetilde{v}_j$, $\widetilde{v}_k$ and $\widetilde{w}_0$ be the lifts in $\overline{Z}_2^{(2)}$ of respectively $p'^{-1}(p_X(v_j))$, $p'^{-1}(p_X(v_k))$ and $p'^{-1}(p_X(w_0))$ with nontrivial associated group.  Since $p'^{-1}(p_X(v_j))$, $p'^{-1}(p_X(v_{k}))$ and $p'^{-1}(p_X(w_{0}))$ are contained in the same connected component of $\overline{Z}_2'-\{z\}$, there exists $\alpha \in \{0,1\}$ such that the preimages by the marking of $Z_2^{(2)}$ of the generators of the groups associated with $\widetilde{v}_j$ and $\widetilde{v}_k$ are respectively $y_i^{\alpha}y_j^Xy_i^{\alpha}$, $y_i^{\alpha}y_k^Xy_i^{\alpha}$ and $y_i^{\alpha}y_k^Zy_i^{\alpha}$. As a consequence, there does not exist a refinement of $Z_2'$ such that the preimages by the marking of the generators of the groups associated with the lifts of $p'^{-1}(p_X(v_j))$, $p'^{-1}(p_X(v_k))$ and $p'^{-1}(p_X(w_0))$ with nontrivial associated groups are respectively $y_iy_j^Xy_i$, $y_k^X$ and $y_k^Z$. Thus, $Z$ and $Z_2'$ do not have any common refinement and $d_{\lk(\mathcal{Y}')}(\mathcal{Z},\mathcal{Z}_2')>2$.
\hfill\qedsymbol

\bigskip

Suppose now that $y_j^Z \neq y_iy_j^Xy_i$ (Cases~$(iii)$ and $(iv)$). Then there exists a leaf $\widetilde{a}_j$ adjacent to $\widetilde{a}$ such that the preimage by the marking of the generator of the group associated with $\widetilde{a}_j$ is $y_j^X$. Let $a_j=p(\widetilde{a}_j)$. Moreover, as $\mathcal{Z} \neq \mathcal{X}$, either there exists $\ell \in \{1,\ldots,m-1\}$ such that either $y_{\ell}^Z=y_iy_{\ell}^Xy_i$ (Case~$(ii)$) or $y_{\ell}^Z=y_{\ell}^X$ (Case~$(iii))$ or there exist $\ell \in \{1,\ldots,m-1\}-\{j\}$ and a leaf $\widetilde{a}_{\ell}$ of $\overline{Z}$ adjacent to $\widetilde{a}$ such that the preimage by the marking of the generator of the group associated with $\widetilde{a}_{\ell}$ is $y_{\ell}^X$ (Case~$(iv)$). By the claim above (see Case~$(ii)$), we can suppose that $y_{\ell}^Z \neq y_iy_{\ell}^Xy_i$.

Let $\mathcal{Z}_3'$ be the homothety class of the marked graph of groups $Z_3'$ defined as follows (see Figure~\ref{construction Z3'} with $y_{\ell}^X=y_{\ell}^Z$):

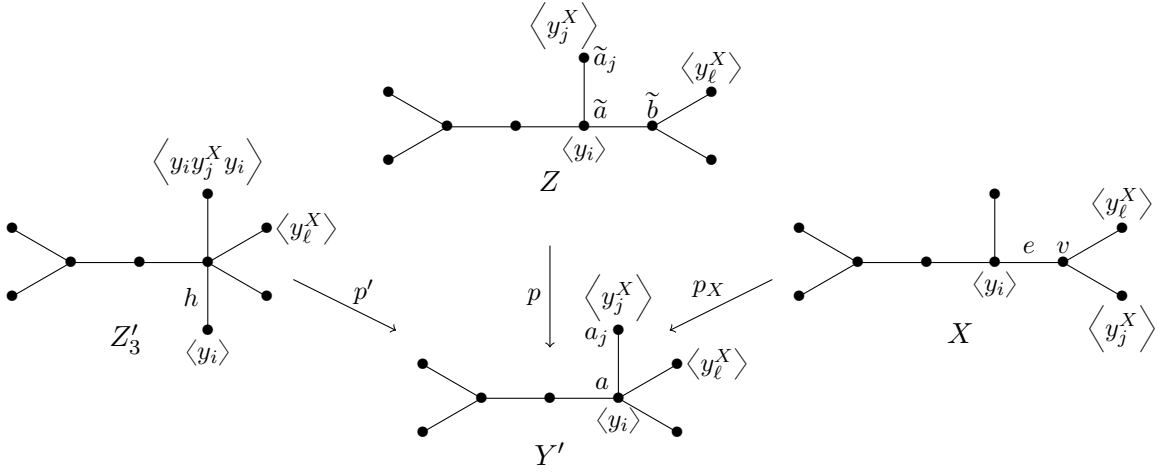
\begin{figure}
\hspace*{-0.5cm}
\begin{tikzpicture}[scale=0.90]
\draw (0.5-7,0) -- (1.5-7,0);
\draw (0.5-7,0) node {$\bullet$};
\draw (1.5-7,0) node {$\bullet$};
\draw (1.5-7,0) -- (2.5-7,0);
\draw (2.5-7,0) node {$\bullet$};
\draw (0.5-7,0) -- (-0.36-7,0.5);
\draw (-0.36-7,0.5) node {$\bullet$};
\draw (0.5-7,0) -- (-0.36-7,-0.5);
\draw (-0.36-7,-0.5) node {$\bullet$};
\draw (2.5-7,0) -- (2.5-7,1);
\draw (2.5-7,1) node {$\bullet$};
\draw (2.5-7,0) -- (2.5-7,-1);
\draw (2.5-7,-1) node {$\bullet$};
\draw (2.5-7,0) -- (3.36-7,0.5);
\draw (3.36-7,0.5) node {$\bullet$};
\draw (2.5-7,0) -- (3.36-7,-0.5);
\draw (3.36-7,-0.5) node {$\bullet$};
\draw (2.5-7,-1) node[below,scale= 0.90] {$\left\langle y_i \right\rangle$};
\draw (2.5-7,-0.5) node[left,scale= 0.90] {$h$};
\draw (2.5-7,1) node[above,scale= 0.90] {$\left\langle y_iy_j^Xy_i \right\rangle$};
\draw (3.36-7,0.5) node[right,scale= 0.90] {$\left\langle y_{\ell}^X \right\rangle$};
\draw (1.3-7,-0.75) node[below, scale=1] {$Z_3'$};

\draw [->] (3.75-7,-0.25) -- (5.25-7,-1);
\draw (4.5-7,-0.75) node[above right, scale=0.9] {$p'$};

\draw (-1,2) -- (0,2);
\draw (0,2) node {$\bullet$};
\draw (-1,2) node {$\bullet$};
\draw (-0,2) -- (1,2);
\draw (1,2) node {$\bullet$};
\draw (1,2) -- (2,2);
\draw (2,2) node {$\bullet$};
\draw (-1,2) -- (-1.86,2.5);
\draw (-1.86,2.5) node {$\bullet$};
\draw (-1,2) -- (-1.86,1.5);
\draw (-1.86,1.5) node {$\bullet$};
\draw (2,2) -- (2.86,2.5);
\draw (2.86,2.5) node {$\bullet$};
\draw (2,2) -- (2.86,1.5);
\draw (2.86,1.5) node {$\bullet$};
\draw (1,2) -- (1,3);
\draw (1,3) node {$\bullet$};
\draw (2,2) node[above,scale= 0.90] {$\widetilde{b}$};
\draw (1,2) node[below,scale= 0.90] {$\left\langle y_i \right\rangle$};
\draw (1,2) node[above right,scale= 0.90] {$\widetilde{a}$};
\draw (0.7,3) node[above,scale= 0.90] {$\left\langle y_j^X \right\rangle$};
\draw (1,3) node[right,scale= 0.90] {$\widetilde{a}_j$};
\draw (2.86,2.5) node[above,scale= 0.90] {$\left\langle y_{\ell}^X \right\rangle$};
\draw (0.5,1.5) node[below, scale=1] {$Z$};

\draw [->] (0.5, 0.25) -- (0.5,-1.25);
\draw (0.5,-0.5) node[left, scale=0.9] {$p$};

\draw (-0.5,-2) -- (0.5,-2);
\draw (-0.5,-2) node {$\bullet$};
\draw (0.5,-2) node {$\bullet$};
\draw (0.5,-2) -- (1.5,-2);
\draw (1.5,-2) node {$\bullet$};
\draw (-0.5,-2) -- (-1.36,-2.5);
\draw (-1.36,-2.5) node {$\bullet$};
\draw (-0.5,-2) -- (-1.36,-1.5);
\draw (-1.36,-1.5) node {$\bullet$};
\draw (1.5,-2) -- (2.36,-1.5);
\draw (2.36,-1.5) node {$\bullet$};
\draw (1.5,-2) -- (2.36,-2.5);
\draw (2.36,-2.5) node {$\bullet$};
\draw (1.5,-2) -- (1.5,-1);
\draw (1.5,-1) node {$\bullet$};
\draw (1.5,-2) node[below,scale= 0.90] {$\left\langle y_i \right\rangle$};
\draw (1.5,-2) node[above left,scale= 0.90] {$a$};
\draw (1.5,-1) node[above,scale= 0.90] {$\left\langle y_j^X \right\rangle$};
\draw (1.5,-1.1) node[left,scale= 0.90] {$a_j$};
\draw (2.36,-1.5) node[right,scale= 0.90] {$\left\langle y_{\ell}^X \right\rangle$};
\draw (0.5,-2.5) node[below, scale=1] {$Y'$};

\draw [<-] (2.25,-1) -- (3.75,-0.25);
\draw (3.2,-0.65) node[above left, scale=0.9] {$p_{X}$};

\draw (-1.5+6.5,0) -- (-0.5+6.5,0);
\draw (-0.5+6.5,0) node {$\bullet$};
\draw (-1.5+6.5,0) node {$\bullet$};
\draw (-0.5+6.5,0) -- (0.5+6.5,0);
\draw (0.5+6.5,0) node {$\bullet$};
\draw (0.5+6.5,0) -- (1.5+6.5,0);
\draw (1.5+6.5,0) node {$\bullet$};
\draw (-1.5+6.5,0) -- (-2.36+6.5,0.5);
\draw (-2.36+6.5,0.5) node {$\bullet$};
\draw (-1.5+6.5,0) -- (-2.36+6.5,-0.5);
\draw (-2.36+6.5,-0.5) node {$\bullet$};
\draw (1.5+6.5,0) -- (2.36+6.5,0.5);
\draw (2.36+6.5,0.5) node {$\bullet$};
\draw (1.5+6.5,0) -- (2.36+6.5,-0.5);
\draw (2.36+6.5,-0.5) node {$\bullet$};
\draw (1.5+6.5,0) node[above,scale= 0.90] {$v$};
\draw (1+6.5,0) node[above,scale= 0.90] {$e$};
\draw (0.5+6.5,0) node[below,scale= 0.90] {$\left\langle y_i \right\rangle$};
\draw (2.36+6.5,-0.5) node[below,scale= 0.90] {$\left\langle y_{j}^X \right\rangle$};
\draw (0.5+6.5,0) -- (0.5+6.5,1);
\draw (0.5+6.5,1) node {$\bullet$};
\draw (2.36+6.5,0.5) node[above,scale= 0.90] {$\left\langle y_{\ell}^X \right\rangle$};
\draw (0+6.5,-0.75) node[below, scale=1] {$X$};

\end{tikzpicture}
\caption{The construction of $Z_3'$ in Lemma~\ref{quasi minimaux et minimaux fixes} in the case $y_{\ell}^X=y_{\ell}^Z$.}\label{construction Z3'}
\end{figure} 

$\bullet$ The underlying graph $\overline{Z}_3'$ of $Z_3'$ is obtained from $\overline{Y}'$ by blowing-up an edge $h$ at $a$ so that one of the two endpoints of $h$ is a leaf. Let $p' \colon \overline{Z}_3' \to \overline{Y}'$ be the projection. Let $x$ be a vertex of the underlying graph of $Z_3'$. Let $y$ be the preimage by the marking of the generator of the group associated with $p'(x)$. 

$\bullet$ If $p'(x) \neq a,a_j$ or if $p'(x)=a$ and $x$ is a leaf, then the preimage by the marking of $Z_3'$ of the generator of the group associated with $x$ is $y$.

$\bullet$ If $p'(x)=a$ and $x$ is not a leaf, then $x$ has trivial associated group.

$\bullet$ If $p'(x)=a_j$, then $x$ is a leaf. Moreover, the preimage by the marking of $Z_3'$ of the generator of the group associated with $x$ is $y_iy_j^Xy_i$.

\medskip

\noindent Since $Z_3'$ has one less vertex with nontrivial associated group which is not a leaf, by the induction hypothesis, $\mathcal{Z}_3'$ is fixed by $f$. What is more, $d_{\lk(\mathcal{Y}')}(\mathcal{Z},\mathcal{Z}_3')=2$. Indeed, a common refinement of $Z$ and $Z_3'$ is obtained from $Z$ by blowing-up the edge $h$ at the vertex of $p^{-1}(a)$ with nontrivial associated group.

\medskip

\noindent{\bf Claim. } In $\lk(\mathcal{Y}')$, we have $d_{\lk(\mathcal{Y}')}(\mathcal{X},\mathcal{Z}_3')>2$.

\medskip

\dem The proof is identical for Cases~$(iii)$ and $(iv)$. Since both $\mathcal{X}$ and $\mathcal{Z}_3'$ have nontrivial negative link with no edge we see by Lemma~\ref{quasi minimal number of vertices}~$(1)$ that $|V\overline{X}|=|V\overline{Z}_3'|$. As both $\overline{X}$ and $\overline{Z}_3'$ are trees, we have $|E\overline{X}|=|E\overline{Z}_3'|$. Thus, as $\mathcal{X} \neq \mathcal{Z}_3'$, we have \mbox{$d_{\lk(\mathcal{Y}')}(\mathcal{X},\mathcal{Z}_3')>1$}. As $\mathcal{Y}'$ has trivial negative link, the only way that $d_{\lk(\mathcal{Y}')}(\mathcal{X},\mathcal{Z}_3')=2$ is that $X$ and $Z_3'$ have a common refinement. Let $z$ be the leaf of $\overline{Z}_3'$ such that $p'(z)=a$. Then the preimage by the marking of the group associated with $z$ is $\left\langle y_i \right\rangle$.

Let $Z_3^{(2)}$ be a refinement of $Z_3'$, let $\overline{Z}_3^{(2)}$ be its underlying graph and let $y_1^{(2)},\ldots,y_n^{(2)}$ be the preimages by the marking of $Z_3^{(2)}$ of the generators of the nontrivial vertex groups of $Z_3^{(2)}$. Since both $Z$ and $Z_3'$ are obtained from $Y'$ by blowing-up an edge at $a$ while applying a twist around an edge adjacent to $a$, a potential common refinement of $Z$ and $Z_3'$ is obtained from $Y'$ by blowing-up a forest while applying a twist around an edge adjacent to $a$. Thus, we may assume that, for all $m \in \{1,\ldots,n\}$, there exists $\alpha_m \in \{0,1\}$ such that $y_m^{(2)}=y_i^{\alpha_m}y_my_i^{\alpha_m}$.
Let $\widetilde{v}_j$ and $\widetilde{v}_{\ell}$ be the lifts in $\overline{Z}_3^{(2)}$ of respectively $p'^{-1}(p_X(v_j))$ and $p'^{-1}(p_X(v_{\ell}))$ with nontrivial associated group. Since $p'^{-1}(p_X(v_j))$ and $p'^{-1}(p_X(v_{\ell}))$ are contained in the same connected component of $\overline{Z}_3'-\{z\}$, there exists $\alpha \in \{0,1\}$ such that the preimages by the marking of $Z_3^{(2)}$ of the generators of the groups associated with $\widetilde{v}_j$ and $\widetilde{v}_{\ell}$ are respectively $y_i^{\alpha+1}y_j^Xy_i^{\alpha+1}$ and $y_i^{\alpha}y_{\ell}^Xy_i^{\alpha}$. As a consequence, there does not exist a refinement $Z_0$ of $Z_3'$ such that the preimages by the marking of the generators of the groups associated with the lifts of $p'^{-1}(p_X(v_j))$ and $p'^{-1}(p_X(v_{\ell}))$ with nontrivial associated groups are respectively $y_j^X$ and $y_{\ell}^X$. Thus, $X$ and $Z_3'$ do not have any common refinement and $d_{\lk(\mathcal{Y}')}(\mathcal{X},\mathcal{Z}_3')>2$.
\hfill\qedsymbol

\medskip

Since $d_{\lk(f(\mathcal{Y}'))}(f(\mathcal{X}),f(\mathcal{Z}_3'))=d_{\lk(\mathcal{Y}')}(\mathcal{X},\mathcal{Z}_3')=2$, the two claims imply that $f(\mathcal{X})=\mathcal{X}$.

\bigskip

We now prove that $f(\mathcal{Y})=\mathcal{Y}$. Let $v_1$ be a vertex of $\overline{Y}$ that is adjacent to at least one leaf. Let $v_2,\ldots,v_{\ell}$ be the leaves adjacent to $v_1$ and, for $i \in \{1,\ldots,\ell\}$, let $y_i$ be the preimage by the marking of the generator of the group associated with $v_i$. Let $\mathcal{Y}'$ be the equivalence class of the marked graph of groups $Y'$ defined as follows (see Figure~\ref{construction Y'}):

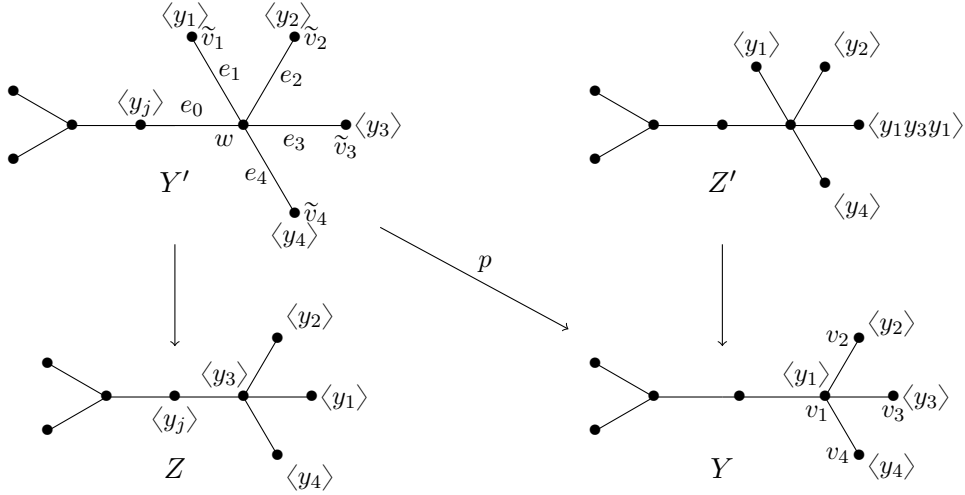
\begin{figure}
\centering
\captionsetup{justification=centering}
\begin{tikzpicture}[scale=0.90]

\draw (-1-4,2) node {$\bullet$};
\draw (-1-4,2) -- (0.5-4,2);
\draw (0-4,2) node {$\bullet$};
\draw (0-4,2) -- (1.5-4,2);
\draw (1.5-4,2) node {$\bullet$};
\draw (-1-4,2) -- (-1.36-4.5,2.5);
\draw (-1.36-4.5,2.5) node {$\bullet$};
\draw (-1-4,2) -- (-1.36-4.5,1.5);
\draw (-1.36-4.5,1.5) node {$\bullet$};
\draw (1.5-4,2) -- (-1.75,3.29);
\draw (-1.75,3.29) node {$\bullet$};
\draw (1.5-4,2) -- (3-4,2);
\draw (3-4,2) node {$\bullet$};
\draw (1.5-4,2) -- (-1.75,0.71);
\draw (-1.75,0.71) node {$\bullet$};
\draw (1.5-4,2) -- (-3.25,3.29);
\draw (-3.25,3.29) node {$\bullet$};
\draw (3-4,2) node[right,scale= 0.90] {$\left\langle y_3 \right\rangle$};
\draw (3-4,2) node[below,scale= 0.90] {$\widetilde{v}_3$};
\draw (3-4.75,2) node[below,scale= 0.90] {$e_3$};
\draw (-1.75,0.71) node[below,scale= 0.90] {$\left\langle y_4 \right\rangle$};
\draw (-1.75,0.71) node[right,scale= 0.90] {$\widetilde{v}_4$};
\draw (-2,1.25) node[left,scale= 0.90] {$e_4$};
\draw (-3.35,3.29) node[above,scale= 0.90] {$\left\langle y_1 \right\rangle$};
\draw (-3.25,3.29) node[right,scale= 0.90] {$\widetilde{v}_1$};
\draw (-3,2.79) node[right,scale= 0.90] {$e_1$};
\draw (-1.85,3.29) node[above,scale= 0.90] {$\left\langle y_2 \right\rangle$};
\draw (-1.75,3.29) node[right,scale= 0.90] {$\widetilde{v}_2$};
\draw (-2.1,2.69) node[right,scale= 0.90] {$e_2$};
\draw (1.5-4,2) node[below left,scale= 0.90] {$w$};
\draw (1.5-4.75,2) node[above,scale= 0.90] {$e_0$};
\draw (0-4,2) node[above,scale= 0.90] {$\left\langle y_j \right\rangle$};
\draw (0.5-4,1.5) node[below, scale=1] {$Y'$};

\draw (-0.5+4,2) node {$\bullet$};
\draw (-0.5+4,2) -- (0.5+4,2);
\draw (0.5+4,2) node {$\bullet$};
\draw (0.5+4,2) -- (1.5+4,2);
\draw (1.5+4,2) node {$\bullet$};
\draw (-0.5+4,2) -- (-1.36+4,2.5);
\draw (-1.36+4,2.5) node {$\bullet$};
\draw (-0.5+4,2) -- (-1.36+4,1.5);
\draw (-1.36+4,1.5) node {$\bullet$};
\draw (1.5+4,2) -- (2+4,2.86);
\draw (2+4,2.86) node {$\bullet$};
\draw (1.5+4,2) -- (2.5+4,2);
\draw (2.5+4,2) node {$\bullet$};
\draw (1.5+4,2) -- (2+4,1.14);
\draw (2+4,1.14) node {$\bullet$};
\draw (1.5+4,2) -- (1+4,2.86);
\draw (1+4,2.86) node {$\bullet$};
\draw (2.5+4,2) node[right,scale= 0.90] {$\left\langle y_1y_3y_1 \right\rangle$};
\draw (2+4,1.14) node[below right,scale= 0.90] {$\left\langle y_4 \right\rangle$};
\draw (1+4,2.86) node[above,scale= 0.90] {$\left\langle y_1 \right\rangle$};
\draw (2+4,2.86) node[above right,scale= 0.90] {$\left\langle y_2 \right\rangle$};
\draw (0.5+4,1.5) node[below, scale=1] {$Z'$};

\draw [->] (0.5+4, 0.25) -- (0.5+4,-1.25);
\draw [->] (0.5-4, 0.25) -- (0.5-4,-1.25);
\draw [->] (0.5-1, 0.5) -- (0.5+1.75,-1);
\draw (0.8,-0.3) node[above right,scale=0.9] {$p$};

\draw (-0.5+4,-2) -- (0.5+4,-2);
\draw (-0.5+4,-2) node {$\bullet$};
\draw (0.75+4,-2) node {$\bullet$};
\draw (0.5+4,-2) -- (1.5+4.5,-2);
\draw (1.5+4.5,-2) node {$\bullet$};
\draw (-0.5+4,-2) -- (-1.36+4,-2.5);
\draw (-1.36+4,-2.5) node {$\bullet$};
\draw (-0.5+4,-2) -- (-1.36+4,-1.5);
\draw (-1.36+4,-1.5) node {$\bullet$};
\draw (1.5+4.5,-2) -- (2.5+4.5,-2);
\draw (2.5+4.5,-2) node {$\bullet$};
\draw (1.5+4.5,-2) -- (1.5+0.5+4.5,-2-0.86);
\draw (1.5+0.5+4.5,-2-0.86) node {$\bullet$};
\draw (1.5+4.5,-2) -- (1.5+0.5+4.5,-2+0.86);
\draw (1.5+0.5+4.5,-2+0.86) node {$\bullet$};
\draw (1.5+4.7,-2) node[above left, scale= 0.90] {$\left\langle y_1 \right\rangle$};
\draw (1.5+4.7,-2) node[below left, scale= 0.90] {$v_1$};
\draw (2+4.5,-2+0.76) node[above right,scale= 0.90] {$\left\langle y_2 \right\rangle$};
\draw (2+4.5,-2+0.86) node[left,scale= 0.90] {$v_2$};
\draw (1.5+0.5+4.5,-2-0.7) node[below right, scale= 0.90] {$\left\langle y_4 \right\rangle$};
\draw (1.5+0.5+4.5,-2-0.86) node[left, scale= 0.90] {$v_4$};
\draw (2.5+4.5,-2) node[right,scale= 0.90] {$\left\langle y_3 \right\rangle$};
\draw (2.5+4.5,-2) node[below,scale= 0.90] {$v_3$};
\draw (0.5+4,-2.75) node[below, scale=1] {$Y$};

\draw (-0.5-4,-2) -- (0.5-4,-2);
\draw (-0.5-4,-2) node {$\bullet$};
\draw (0.5-4,-2) node {$\bullet$};
\draw (0.5-4,-2) -- (1.5-4,-2);
\draw (1.5-4,-2) node {$\bullet$};
\draw (-0.5-4,-2) -- (-1.36-4,-2.5);
\draw (-1.36-4,-2.5) node {$\bullet$};
\draw (-0.5-4,-2) -- (-1.36-4,-1.5);
\draw (-1.36-4,-1.5) node {$\bullet$};
\draw (1.5-4,-2) -- (2.5-4,-2);
\draw (2.5-4,-2) node {$\bullet$};
\draw (1.5-4,-2) -- (1.5+0.5-4,-2-0.86);
\draw (1.5+0.5-4,-2-0.86) node {$\bullet$};
\draw (1.5-4,-2) -- (1.5+0.5-4,-2+0.86);
\draw (1.5+0.5-4,-2+0.86) node {$\bullet$};
\draw (1.5-3.8,-2) node[above left, scale= 0.90] {$\left\langle y_3 \right\rangle$};
\draw (2-4,-2+0.86) node[above right,scale= 0.90] {$\left\langle y_2 \right\rangle$};
\draw (1.5+0.5-4,-2-0.86) node[below right, scale= 0.90] {$\left\langle y_4 \right\rangle$};
\draw (2.5-4,-2) node[right,scale= 0.90] {$\left\langle y_1 \right\rangle$};
\draw (0.5-4,-2) node[below,scale= 0.90] {$\left\langle y_j \right\rangle$};
\draw (0.5-4,-2.75) node[below, scale=1] {$Z$};

\end{tikzpicture}
\caption{The constructions of $Y'$ and $Z'$ in Lemma~\ref{quasi minimaux et minimaux fixes}.}\label{construction Y'}
\end{figure}

$\bullet$ The underlying graph $\overline{Y}'$ of $Y'$ is obtained from $Y$ by blowing-up an edge $e_{1}$ at $v_1$ such that one of the endpoint of $e_{1}$ is a leaf. Let $p \colon Y' \to Y$ be the natural projection. Let $x$ be a vertex of $\overline{Y}'$, and let $y$ be the preimage by the marking of the generator of the group associated with $p(x)$. 

$\bullet$ If $p(x) \neq v_1$ or if $p(x)=v_1$ and $x$ is a leaf, then the preimage by the marking of $Y'$ of the generator of the group associated with $x$ is $y$.

$\bullet$ If $p(x)=v_1$ and $x$ is not a leaf, then $v_1$ has trivial associated group.

\noindent By the previous step, as $Y'$ has $k$ vertices with nontrivial associated groups that are not leaves, and as $\lk_-(\mathcal{Y}')$ is nontrivial and has no edge, we see that $f(\mathcal{Y}')=\mathcal{Y}'$. By the second claim in the proof of Lemma~\ref{quasi minimaux et minimaux fixes}, the negative link of $f(\mathcal{Y}')$ is nontrivial and has no edge. Therefore, $f$ preserves the negative link of $\mathcal{Y}'$ and $f(\mathcal{Y}) \in \lk_-(\mathcal{Y}')$. Let $w$ be the endpoint of $e_{1}$ with trivial associated group. For $i \in \{1,\ldots,\ell\}$, let $\widetilde{v}_i$ be the leaf of $\overline{Y}'$ which lifts $v_i$. Let $e_2,\ldots,e_{\ell}$ be the edges of $\overline{Y}'$ such that for all $i \in \{2,\ldots,\ell\}$, the endpoints of $e_i$ are $\widetilde{v}_i$ and $w$. 

We claim that there exists a unique $i \in \{1,\ldots,\ell\}$ such that a representative of $f(\mathcal{Y})$ is obtained from $Y'$ by contracting $e_i$. Indeed, as $\lk_-(\mathcal{Y}')$ is nontrivial and has no edge, Lemma~\ref{quasi minimal number of vertices}~$(1)$ implies that $\overline{Y}'$ has exactly one vertex with trivial associated group, namely $w$. Therefore, a representative of $f(\mathcal{Y})$ is obtained from $Y'$ by contracting a unique edge adjacent to $w$. 

Suppose towards a contradiction that there exists an edge $e_0$ in $\overline{Y}'$ between $w$ and a vertex $w'$ with nontrivial associated group that is not a leaf and such that a representative of $f(\mathcal{Y})$ is obtained from $Y'$ by collapsing $e_0$. Let $y_j$ be the preimage by the marking of $Y'$ of the generator of the group associated with $w'$. Let $\mathcal{Z}_0$ be the homothety class of the marked graph of groups obtained from $Y'$ by contracting $e_0$. By induction hypothesis, $f(\mathcal{Z}_0)=\mathcal{Z}_0$. Thus $f(\mathcal{Y}) \neq \mathcal{Z}_0$. 

Thus, there exists a unique $i \in \{1,\ldots,\ell\}$ such that a representative of $f(\mathcal{Y})$ is obtained from $Y'$ by contracting $e_i$. We claim that $i=1$. Indeed, for $i \neq 1$, let $\mathcal{Z}$ be the equivalence class of the marked graph of groups $Z$ obtained from $Y'$ by collapsing $e_i$. Let $\mathcal{Z}'$ be the equivalence class of the marked graph of groups $Z'$ whose underlying graph is $\overline{Y}'$ and such that the preimage by the markings of the generators of the groups associated with $\widetilde{v}_1,\ldots,\widetilde{v}_{\ell}$ are $y_1,\ldots,y_{i-1},y_1y_iy_1,y_{i+1},\ldots,y_{\ell}$. Then $\mathcal{Z}' \in \lk_+(\mathcal{Y})$ because $Z'$ is obtained from $Y$ first by precomposing the marking of $Y$ by the automorphism which sends $y_i$ to $y_1y_iy_1$ and fixes all the other $y_i$ and then blowing-up an edge at $v_1$ such that one of the endpoints of this edge is a leaf and then. However, $\mathcal{Z}' \notin \lk(\mathcal{Z})$ because the vertex of $\overline{Z}$ whose preimage by the marking of the associated group is $\left\langle y_i \right\rangle$ is a leaf. Therefore there is no refinement of $Z$ such that there exist two vertex groups of the refinement such that the preimage by the marking of the generators of the vertex groups are respectively $y_1y_iy_1$ and $y_k$. As $f(\mathcal{Z}')=\mathcal{Z}'$ by the previous step, we have $f(\mathcal{Y}) \in \lk(\mathcal{Z}')$. Therefore, $f(\mathcal{Y}) \neq \mathcal{Z}$ and $f(\mathcal{Y})=\mathcal{Y}$.
\hfill\qedsymbol

\bigskip

We can now show the injectivity of the homomorphism $\Aut(K_n) \to \Aut(L_n)$.

\begin{prop}\label{injectivity Kn to Ln}
Let $n \geq 4$. Let $f \in \Aut(K_n)$ such that $f|_{O_n}=\mathrm{id}_{O_n}$ and $f|_{F_n}=\mathrm{id}_{F_n}$. Then $f=\mathrm{id}_{K_n}$.
\end{prop}

\dem Let $k \in \NN$ and let $\mathcal{X} \in VK_n$ be such that the underlying graph $\overline{X}$ of a representative $X$ of $\mathcal{X}$ has exactly $k$ vertices with trivial associated group. We prove by induction on $k$ that $f(\mathcal{X})=\mathcal{X}$. If $k=0$, then $\mathcal{X}$ has trivial negative link. Thus, by Lemma~\ref{quasi minimaux et minimaux fixes}, we have $f(\mathcal{X})=\mathcal{X}$. 

Suppose now that $k \geq 1$. Then, as any representative of an element of $\lk_-(\mathcal{X})$ is obtained from $X$ by collapsing at least an edge, by the induction hypothesis, we have $f|_{\lk_-(\mathcal{X})}=\mathrm{id}|_{\lk_-(\mathcal{X})}$ and $\lk_-(f(\mathcal{X}))=\lk_-(\mathcal{X})$. 

Suppose towards a contradiction that $f(\mathcal{X}) \neq \mathcal{X}$. Let $Y$ be a representative of $f(\mathcal{X})$, and let $\overline{Y}$ be the underlying graph of $Y$. By the induction hypothesis, $\overline{Y}$ has at least $k$ vertices with trivial associated group. Since $\mathcal{X} \neq f(\mathcal{X})$, there exists an edge $e \in E\overline{Y}$ such that the free factor decomposition of $W_n$ induced by $e$ is distinct from the free factor decomposition induced by any edge of $\overline{X}$. Let $\mathcal{Z} \in \lk_-(\mathcal{X})$. Let $Z$ be a representative of $\mathcal{Z}$ obtained from $X$ by collapsing a forest, and let $\overline{Z}$ be the underlying graph of $Z$. By Remark~\ref{free factor decomposition and canonical lift}~$(1)$, for any edge $f \in E\overline{Z}$, there exists an edge $\widetilde{f} \in E\overline{X}$ such that the free factor decomposition induced by $f$ is the same one as the free factor decomposition induced by $\widetilde{f}$. Thus, there does not exist any edge of $\overline{Z}$ which induces, up to global conjugation, the same free factor decomposition of $W_n$ as $e$. But, for any edge $f$ of $\overline{Y}$, there exists $\mathcal{Z}' \in \lk_-(\mathcal{Y})$, a representative $Z'$ of $\mathcal{Z'}$ with underlying graph $\overline{Z}'$ and an edge $g$ of $\overline{Z}'$ such that the free factor decomposition induced by $e$ is the same as the one induced by $g$ ($Z$ is obtained from $Y$ by contracting an edge distinct from $f$). This contradicts the fact that $\lk_-(f(\mathcal{X}))=\lk_-(\mathcal{X})$. Thus $f(\mathcal{X})=\mathcal{X}$ and $f=\mathrm{id}_{K_n}$.
\hfill\qedsymbol

\bigskip

\noindent{\bf Proof of Theorem~\ref{Rigidity Kn}. } Let $n \geq 4$. The injectivity is immediate since the homomorphism $\Out(W_n) \to \Aut(L_n)$ is injective by Theorem~\ref{rigidity Ln} and since $L_n$ is a subgraph of $K_n$. We now prove surjectivity. Let $f \in \Aut(K_n)$. By Proposition~\ref{0star F star preserved}, the automorphism $f$ induces an automorphism $\widetilde{f} \in \Aut(L_n)$. By Theorem~\ref{rigidity Ln}, $\widetilde{f}$ is induced by an element $\gamma \in \Out(W_n)$. Since the homomorphism $\Aut(K_n) \to \Aut(L_n)$ is injective by Proposition~\ref{injectivity Kn to Ln}, $f$ is induced by $\gamma$. This concludes the proof.
\hfill\qedsymbol

\section{Rigidity of the simplicial completion of $K_n$}\label{section FS} 

Let $n \geq 4$. A \emph{splitting} of $W_n$ is a minimal, simplicial $W_n$-action on a simplicial tree $S$ and such that:

\medskip

\noindent{$(1)$ } The finite graph $W_n \backslash S$ is not empty and not reduced to a point.

\medskip

\noindent{$(2)$ } Vertices of $S$ with trivial stabilizer have degree at least $3$.

\medskip

Here \emph{minimal} means that $W_n$ does not preserve any proper subtree of $S$. A splitting $S$ of $W_n$ is \emph{free} if all edge stabilizers are trivial. A splitting $S'$ is a \emph{blow-up}, or equivalently a \emph{refinement}, of a splitting $S$ if $S$ is obtained from $S'$ by collapsing some edge orbits in $S'$. Two splittings are \emph{compatible} if they have a common refinement. If $k \geq 1$ is an integer, a free splitting $S$ is a \emph{$k$-edge free splitting} if $W_n \backslash S$ has exactly $k$ edges. An \emph{$F$-one-edge free splitting} is a one-edge free splitting $S$ such that one of the vertex groups of $W_n \backslash S$ is isomorphic to $W_{n-1}$ while the other vertex group is isomorphic to $F$.

\bigskip

The simplicial completion of $K_n$, denoted by $\overline{K}_n$, is the flag complex such that:

$\bullet$~ The vertices of $\overline{K}_n$ are the equivalence classes of free splittings of $W_n$, where two free splittings $S$ and $S'$ are equivalent if there exists a $W_n$-equivariant homeomorphism between them.

$\bullet$~ Two equivalence classes of free splittings $\mathcal{S}$ and $\mathcal{S}'$ are adjacent in $\overline{K}_n$ if there exist $S \in \mathcal{S}$ and $S' \in \mathcal{S}'$ such that $S$ refines $S'$ or conversely. 

\medskip

In the literature, this complex is also called the \emph{free splitting complex}. The free splitting complex appears as well in the study of the outer automorphism group of a free group of finite rank and more generally in the study of the outer automorphism group of a free product of groups (see \cite{AramayonaSouto2011,HandelMosher14,HandelMosher13}). In particular, Handel-Mosher (\cite{HandelMosher14,HandelMosher13}) in the case of $\Out(F_N)$ and Handel-Mosher and Horbez (\cite{HandelMosher14,Horbez18}) in the case of the outer automorphism group of a free product of groups proved that this complex is Gromov hyperbolic. 

\medskip

We have a canonical injective homomorphism $K_n \hookrightarrow \overline{K}_n$ defined as follows. Let \mbox{$\mathcal{X} \in VK_n$} be the equivalence class of a marked graph of groups, and let $X$ be a representative of $\mathcal{X}$. Let $S$ be a Bass-Serre tree corresponding to $X$, and let $\mathcal{S}$ be the equivalence class of $S$. Then the map 
$$\begin{array}{ccccl}
\Phi & \colon & K_n & \to & \overline{K}_n \\
{} & {} & \mathcal{X} & \mapsto & \mathcal{S}
\end{array}
$$ 
is a well-defined injective homomorphism. From now on, we identify $K_n$ with its image in $\overline{K}_n$.

The group $\Aut(W_n)$ acts on $\overline{K}_n$ by precomposition of the action. For any \mbox{$\alpha \in \Inn(W_n)$} and for any $S \in \overline{K}_n$, we have $\alpha(S)=S$. Therefore the action of $\Aut(W_n)$ induces an action of $\Out(W_n)$. 

\bigskip

In this section, we prove Theorem~\ref{Rigidity FS}. In order to do so, we first show that any automorphism of $\overline{K}_n$ preserves $K_n$. Thus, we have a restriction homomorphism $\Aut(\overline{K}_n) \to \Aut(K_n)$ which, as we will see, turns out to be injective. Theorem~\ref{Rigidity FS} then follows from Theorem~\ref{Rigidity Kn}.

We first characterise the vertices of $K_n$ in $\overline{K}_n$.

\begin{prop}\label{infinite adjacency vertices FS}
Let $n \geq 4$. Let $\mathcal{S} \in V\overline{K}_n$. If $\mathcal{S} \in VK_n$, then $\mathcal{S}$ has finite valence in $\overline{K}_n$. 
If $\mathcal{S} \in V\overline{K}_n-VK_n$, then $\mathcal{S}$ has infinite valence in $\overline{K}_n$.
\end{prop}

\dem Suppose that $\mathcal{S} \in VK_n$, and let $\mathcal{S}' \in \lk(\mathcal{S})$. Let $S$ and $S'$ be representatives of $\mathcal{S}$ and $\mathcal{S}'$. If $S$ refines $S'$, then $W_n\backslash S'$ is obtained from $W_n \backslash S$ by collapsing a forest. Since $W_n \backslash S$ is a finite tree, there are only finitely many possibilities for $W_n \backslash S'$, hence finitely many possibilities for $S'$. If $S'$ refines $S$, then, since $\mathcal{S} \in K_n$, the equivalence class $\mathcal{S}'$ also belongs to $K_n$. Thus, we have $\mathcal{S}' \in \lk_+^{K_n}(\mathcal{S})$ where $\lk_+^{K_n}(\mathcal{S})$ is the positive link of $\mathcal{S}$ in $K_n$. Since $\lk_+^{K_n}(\mathcal{S})$ is finite, there are only finitely many possibilities for $\mathcal{S}'$. Hence $\lk(\mathcal{S})$ is finite.

Now suppose that $\mathcal{S} \in V\overline{K}_n-VK_n$. Let $S$ be a representative of $\mathcal{S}$. Since we have $\mathcal{S} \in V\overline{K}_n-VK_n$, there exists a vertex of $S$ whose stabilizer contains a subgroup $G$ of $W_n$ isomorphic to $W_2$. Since $\Aut(W_2)$ is isomorphic to $W_2$ (see e.g. \cite[Lemma 1.4.2]{thomasautotower}), we see that $\Stab(\mathcal{S})$ is infinite by Proposition~\ref{Levitt stab}. Moreover, we claim that there exists $\mathcal{S'} \in \lk(\mathcal{S}) \cap VK_n$. Indeed, let $v$ be a vertex of $W_n \backslash S$ whose associated group is isomorphic to $W_i$ with $i \geq 2$. then one can construct an element $T \in K_i$ and then blow-up $T$ at $v$. The equivalence class of the result is an element in $\lk(\mathcal{S})$. Applying the process to every vertex of $W_n \backslash S$ with infinite associated vertex group gives an element $\mathcal{S'} \in \lk(\mathcal{S}) \cap VK_n$. As $\Stab(\mathcal{S})$ acts on $\lk(\mathcal{S})$, the orbit of $\mathcal{S}'$ under the action of $\Stab(\mathcal{S})$ is infinite (recall that $\Stab(\mathcal{S}')$ is finite by Proposition~\ref{Levitt stab}). Thus $\lk(\mathcal{S})$ is infinite.
\hfill\qedsymbol

\bigskip

Thus, Proposition~\ref{infinite adjacency vertices FS} tells us that any automorphism of $\overline{K}_n$ preserves $K_n$. This gives a restriction homomorphism $$\Aut(\overline{K}_n) \to \Aut(K_n).$$ In the rest of the section, we prove that this homomorphism is injective. In order to show this, we first prove that any automorphism of $\overline{K}_n$ which fixes $K_n$ pointwise also fixes the set of one-edge free splittings pointwise. We will then conclude by the following proposition, due to Scott and Swarup.

\begin{theo}\cite[Theorem 2.5]{ScottSwarup}\label{theorem scott swarup}
Let $n \geq 4$. Any set $\{S_1,\ldots,S_k\}$ of pairwise distinct, pairwise compatible, one-edge free splittings of $W_n$ has a unique refinement $S$ such that $W_n \backslash S$ has exactly $k$ edges. If $S$ is a free splitting such that $W_n \backslash S$ has exactly $k$ edges, then $S$ refines exactly $k$ distinct one-edge free splittings.
\end{theo}

The next lemma is inspired by \cite[Lemma 2.3]{HorbezWade15} due to Horbez and Wade.

\begin{lem}\label{F one edge free splitting preserved}
Let $n \geq 4$. For all $\mathcal{S} \in V\overline{K}_n$, the following assertions are equivalent.
\begin{enumerate}
\item There exists $S \in \mathcal{S}$ such that $S$ is an $F$-one-edge free splitting.
\item The equivalence class $\mathcal{S}$ satisfies the following properties.
\begin{enumerate}
\item The link of $\mathcal{S}$ is infinite.
\item There exists a $\{0\}$-star $\mathcal{X}$ such that $\mathcal{X} \in \lk(\mathcal{S})$.
\item There exist $\mathcal{S}_1,\mathcal{S}_2 \in \lk(\mathcal{S})$ such that $d_{\overline{K}_n}(\mathcal{S}_1,\mathcal{S}_2)=2$ and such that $\mathcal{S}_1-\mathcal{S}-\mathcal{S}_2$ is the unique path of length $2$ joining $\mathcal{S}_1$ and $\mathcal{S}_2$.
\end{enumerate} 
\end{enumerate}
\end{lem}

\dem We first prove that $(1)$ implies $(2)$. Let $S \in \mathcal{S}$ be an $F$-one-edge free splitting. Then $\mathcal{S} \notin VK_n$ and Proposition~\ref{infinite adjacency vertices FS} implies that $\lk(\mathcal{S})$ is infinite, which proves Property~$(a)$.

In order to prove Property~$(b)$, let 
$$W_n=\left\langle x_1,\ldots, x_{n-1} \right\rangle \ast \left\langle x_n \right\rangle$$ 
be the free factor decomposition of $W_n$ induced by $S$. Let $X$ be the $\{0\}$-star such that, if $w_1,\ldots,w_n$ are the leaves of $\overline{X}$, and if $i \in \{1,\ldots,n\}$, then the stabilizer of $w_i$ is $\left\langle x_i \right\rangle$. Let $\mathcal{X}$ be the equivalence class of $X$. Then $\mathcal{X} \in \lk(\mathcal{S})$. 

In order to prove Property~$(c)$, let $S_1$ be the $2$-edge free splitting induced by the decomposition 
$$W_n=\left\langle x_1,x_2 \right\rangle \ast \left\langle x_3,\ldots, x_{n-1} \right\rangle \ast \left\langle x_n \right\rangle,$$ 
where the preimage by the marking of the group associated with the central vertex of $W_n \backslash S_1$ is $\left\langle x_3\ldots, x_{n-1} \right\rangle$. Let $S_2$ be the $2$-edge free splitting induced by the decomposition 
$$W_n=\left\langle x_1,x_3 \right\rangle \ast \left\langle x_2,\widehat{x_3},\ldots, x_{n-1} \right\rangle \ast \left\langle x_n \right\rangle,$$ where the preimage by the marking of the group associated with the central vertex of $W_n \backslash S_2$ is $\left\langle x_2,\widehat{x_3},\ldots, x_{n-1} \right\rangle$.
For $i \in \{1,2\}$, let $\mathcal{S}_i$ be the equivalence class of $S_i$. Then $\mathcal{S}_1,\mathcal{S}_2 \in \lk(\mathcal{S})$. Moreover, the equivalence classes $\mathcal{S}_1$ and $\mathcal{S}_2$ are not adjacent in $\overline{K}_n$ since both $S_1$ and $S_2$ are $2$-edge free splittings, thus, there does not exist $i \in \{1,2\}$, $j \in \{1,2\}-\{i\}$ such that $S_i$ collapses onto $S_{j}$. So $d_{\overline{K}_n}(\mathcal{S}_1,\mathcal{S}_2)=2$.

\medskip

\noindent{\bf Claim. } Let $\mathcal{T} \in V\overline{K}_n$ be such that $\mathcal{S}_1-\mathcal{T}-\mathcal{S}_2$ is a path of length $2$ joining $\mathcal{S}_1$ and $\mathcal{S}_2$. Then $\mathcal{S}=\mathcal{T}$.
\medskip

\dem Suppose towards a contradiction that there exists a representative $T$ of $\mathcal{T}$ such that $T$ is a common refinement of $S_1$ and $S_2$. For $i \in \{1,\ldots,n\}$, let $v_i$ be the only vertex of $T$ fixed by $x_i$. Note that, for $i \neq j$, the vertices $v_i$ and $v_j$ may not be distinct. Since $T$ refines $S_1$, for every edge $e \in ET$, one of the following holds:

\noindent{$\bullet$ } the vertices $v_1$, $v_2$ and $v_3$ belong to the same connected component of $T-\{\mathring{e}\}$,

\noindent{$\bullet$ } the vertices $v_1$ and $v_2$ belong to a connected component of $T-\{\mathring{e}\}$ distinct from the one that contains $v_3$,

\noindent{$\bullet$ } there exist $i \in \{1,2\}$ and $j \in \{1,2\}-\{i\}$ such that $v_i$ is in a connected component  of $T-\{\mathring{e}\}$ distinct from the one containing $v_j$, $v_3$ and $v_n$. 

But if $T$ refines $S_2$, there exists $e \in ET$ such that $v_1$ and $v_3$ belong to a connected component of $T-\{\mathring{e}\}$ distinct from the one that contains $v_2$ and $v_n$. This leads to a contradiction.

Thus, there exists a representative $T$ of $\mathcal{T}$ such that both $S_1$ and $S_2$ collapse to $T$. As $S$ is the only such one-edge free splitting, the claim follows.
\hfill\qedsymbol

\bigskip

We now prove that $(2)$ implies $(1)$. Suppose that $\mathcal{S}$ satisfies the properties of Assertion~$(2)$ of the lemma.

\medskip

\noindent{\bf Claim. } Property~$(c)$ implies that $\mathcal{S}$ has a representative $S$ that is either a one-edge free splitting or is such that there is no free splitting of $W_n$ that properly refines $S$.

\medskip

\dem Let $\mathcal{S}_1$ and $\mathcal{S}_2$ be as in Property~$(c)$, and, for $i \in \{1,2\}$, let $S_i$ be a representative of $\mathcal{S}_i$.  Let $S$ be a representative of $\mathcal{S}$. There are three cases to distinguish. 

\begin{itemize}
\item If $S$ refines $S_1$ and if $S_2$ refines $S$, then $S_2$ refines $S_1$, so that $d_{\overline{K}_n}(\mathcal{S}_1,\mathcal{S}_2) \leq 1$. This leads to a contradiction.
\item If $S$ refines both $S_1$ and $S_2$, then there does not exist any proper refinement of $S$ as this would contradict the uniqueness of the path of length $2$ between $\mathcal{S}_1$ and $\mathcal{S}_2$.
\item If $S$ is refined by both $S_1$ and $S_2$, then $S$ is a one-edge free splitting as otherwise there would exist a splitting $S'$ that is properly refined by $S$. This would contradict the uniqueness of the path.
\end{itemize}
The claim follows.
\hfill\qedsymbol

\medskip

Since a free splitting which has no proper refinement is in $K_n$, the above claim, Property~$(a)$ and Proposition~\ref{infinite adjacency vertices FS} imply that $S$ is a one-edge free splitting. Property~$(b)$ implies in fact that $S$ is an $F$-one-edge free splitting as the $F$-one-edge-free splittings are the only one-edge free splittings that are adjacent to a $\{0\}$-star. The lemma follows.
\hfill\qedsymbol

\begin{lem}\label{F one edge free splitting fixed}
Let $n \geq 4$. Let $f \in \Aut(\overline{K}_n)$ be such that $f|_{L_n}=\mathrm{id}_{L_n}$. Let $\mathcal{S}$ be the equivalence class of an $F$-one-edge free splitting $S$. Then $f(\mathcal{S})=\mathcal{S}$.
\end{lem}

\begin{figure}[ht]
\centering
\begin{tikzpicture}[scale=2]
\draw (0:0) -- (150:1);
\draw (0:0) -- (210:1);
\draw (0:0) node {$\bullet$};
\draw (210:1) node {$\bullet$};
\draw (150:1) node {$\bullet$};
\draw[dotted] (160:0.93) -- (200:0.93);
\draw (0:0) node[right, scale=0.9] {$\;\;\left\langle x_1 \right\rangle$};

\draw (150:1) node[above left, scale=0.9] {$\left\langle x_2 \right\rangle$};
\draw (210:1) node[below left, scale=0.9] {$\left\langle x_n \right\rangle$};
\end{tikzpicture}
\begin{tikzpicture}[scale=2]
\draw (0:0) -- (150:1);
\draw (0:0) -- (210:1);
\draw (0:0) node {$\bullet$};
\draw (210:1) node {$\bullet$};
\draw (150:1) node {$\bullet$};
\draw[dotted] (160:0.93) -- (200:0.93);
\draw (0:0) node[right, scale=0.9] {$\;\;\left\langle x_n \right\rangle$};

\draw (150:1) node[above left, scale=0.9] {$\left\langle x_1 \right\rangle$};
\draw (210:1) node[below left, scale=0.9] {$\left\langle x_{n-1} \right\rangle$};
\end{tikzpicture}
\caption{The $F$-stars $X$ and $X'$ of the proof of Lemma~\ref{F one edge free splitting fixed}.}\label{proof X}
\end{figure}
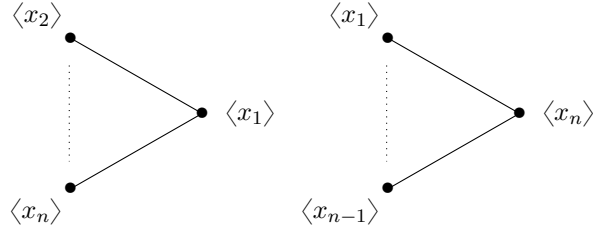

\dem As $f \in \Aut(\overline{K}_n)$, Corollary~\ref{0star F star preserved} and Lemmas~\ref{infinite adjacency vertices FS} and \ref{F one edge free splitting preserved} imply that $f(\mathcal{S})$ is the equivalence class of an $F$-one-edge free splitting $S'$. Let 
$$W_n=\left\langle x_1,\ldots, x_{n-1} \right\rangle \ast \left\langle x_n \right\rangle$$ 
be the free factor decomposition of $W_n$ induced by $S$. Let $\mathcal{X}$ be the equivalence class of the $F$-star $X$ represented in Figure~\ref{proof X} on the left. 

Since $f(\mathcal{X})=\mathcal{X}$, the free splitting $S'$ is an $F$-one-edge free splitting obtained from $X$ by collapsing $n-1$ edges. But if $T$ is an $F$-one-edge free splitting obtained from $X$ by collapsing $n-1$ edges, then there exists $i \in \{1,\ldots,n\}$ such that the free factor decomposition of $W_n$ induced by $T$ is 
$$W_n=\left\langle x_1,\ldots, \widehat{x}_i,\ldots, x_{n} \right\rangle \ast \left\langle x_i \right\rangle.
$$
For $i \in \{1,\ldots,n\}$, we will denote by $T_i$ the $F$-one-edge free splitting with associated free factor decomposition $\left\langle x_1,\ldots, \widehat{x}_i,\ldots, x_{n} \right\rangle \ast \left\langle x_i \right\rangle$, and by $\mathcal{T}_i$ its equivalence class. For $i \neq n$, the free splitting $T_i$ is a collapse of the $F$-star $X'$ depicted in Figure~\ref{proof X} on the right, whereas $S$ is not a collapse of $X'$.

Let $\mathcal{X}'$ be the equivalence class of $X'$. Since $f(\mathcal{X}')=\mathcal{X}'$, we have that $f(\mathcal{S})$ is not adjacent to $\mathcal{X}'$. But, for all $i \neq n$, the equivalence class $\mathcal{T}_i$ is adjacent to $\mathcal{X}'$. Thus, for all $i \neq n$, we have $f(\mathcal{S}) \neq \mathcal{T}_i$. Therefore, as $\mathcal{S}=\mathcal{T}_n$, we conclude that $f(\mathcal{S})=\mathcal{S}$.
\hfill\qedsymbol

\bigskip

\noindent{\bf Proof of Theorem~\ref{Rigidity FS}. } By Proposition~\ref{infinite adjacency vertices FS}, there exists a homomorphism $$\Aut(\overline{K}_n) \to \Aut(K_n)$$ induced by the restriction to $K_n$. In order to prove Theorem~\ref{Rigidity FS}, it suffices to prove that this homomorphism is injective. Let $f \in \Aut(\overline{K}_n)$ be such that $f|_{K_n}=\mathrm{id}_{K_n}$. Let us prove that $f=\mathrm{id}$. By Theorem~\ref{theorem scott swarup}, it suffices to prove that, for any equivalence class $\mathcal{S}$ of a one-edge free splitting $S$, we have $f(\mathcal{S})=\mathcal{S}$. Indeed, let $\mathcal{S}$ be the equivalence class of a free splitting. Then, by Theorem~\ref{theorem scott swarup}, there exist $k$ one-edge free splittings $\mathcal{S}_1,\ldots,\mathcal{S}_k$ such that $\mathcal{S}$ is the unique vertex of $\overline{K}_n$ such that, for all $i \in \{1,\ldots,k\}$, $\mathcal{S}$ is adjacent to $\mathcal{S}_i$. Thus, if, for any equivalence class $\mathcal{S}$ of a one-edge free splitting $S$, we have $f(\mathcal{S})=\mathcal{S}$, then $f=\mathrm{id}$. 

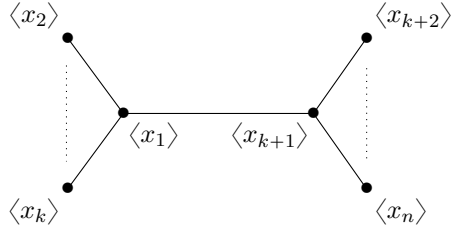
\begin{figure}[ht]
\centering
\captionsetup{justification=centering}
\begin{tikzpicture}[scale=2]
\draw (180:0.5) -- (0:0.75);
\draw (0:0.75) -- (1.1,0.5);
\draw (0:0.75) -- (1.1,-0.5);
\draw (180:0.5) -- (210:1);
\draw (180:0.5) -- (150:1);
\draw[dotted] (160:0.93) -- (200:0.93);
\draw[dotted] (1.1,0.3) -- (1.1,-0.3);
\draw (180:0.5) node {$\bullet$};
\draw (0:0.75) node {$\bullet$};
\draw (210:1) node {$\bullet$};
\draw (150:1) node {$\bullet$};
\draw (1.1,0.5) node {$\bullet$};
\draw (1.1,-0.5) node {$\bullet$};
\draw (0:0.77) node[below left, scale=0.9] {$\left\langle x_{k+1} \right\rangle$};
\draw (180:0.52) node[below right, scale=0.9] {$\left\langle x_{1} \right\rangle$};
\draw (150:1) node[above left, scale=0.9] {$\left\langle x_{2} \right\rangle$};
\draw (210:1) node[below left, scale=0.9] {$\left\langle x_{k} \right\rangle$};
\draw (1.1,0.5) node[above right, scale=0.9] {$\left\langle x_{k+2} \right\rangle$};
\draw (1.1,-0.5) node[below right, scale=0.9] {$\left\langle x_{n} \right\rangle$};
\end{tikzpicture}
\caption{The free splitting $X$ of the proof of Theorem~\ref{Rigidity FS}.}\label{free splitting proof FS}
\end{figure}

Suppose that $\mathcal{S}$ is the equivalence class of a one-edge free splitting $S$. The case where $S$ is an $F$-one-edge free splitting was proved in Lemma~\ref{F one edge free splitting fixed}. If $S$ is not an $F$-one-edge free splitting, let $W_n=\left\langle x_1,\ldots,x_k \right\rangle \ast \left\langle x_{k+1},\ldots,x_n\right\rangle$ be the free factor decomposition of $W_n$ induced by $S$, with $2 \leq k \leq n-2$. Let $X$ be the free splitting of $W_n$ depicted in Figure~\ref{free splitting proof FS}, and let $\mathcal{X}$ be its equivalence class. 

Then $\mathcal{X} \in VK_n$, so $f(\mathcal{X})=\mathcal{X}$. As $\mathcal{S} \in \lk(\mathcal{X})$, we also have that $f(\mathcal{S}) \in \lk(\mathcal{X})$. Moreover, $f(\mathcal{S}) \notin VK_n$ by Proposition~\ref{infinite adjacency vertices FS}. Thus, a representative of $f(\mathcal{S})$ is obtained from $X$ by collapsing a forest $F$.

\medskip

\noindent{\bf Claim. } Any splitting $S'$ distinct from $S$ and obtained from $X$ by collapsing a forest is either an $F$-one-edge free splitting or is adjacent to an $F$-one-edge free splitting. 

\medskip

\dem If $S' \neq S$ is obtained from $X$ by collapsing a forest, and if $S'$ is not an $F$-one-edge free splitting, there exists an edge $e \in V(W_n \backslash S')$ such that $e$ is adjacent to a leaf. This edge determines an $F$-one-edge free splitting adjacent to $S'$.
\hfill\qedsymbol

\medskip

Thus, by Lemma~\ref{F one edge free splitting fixed}, any equivalence class $\mathcal{S}' \in \lk_-(\mathcal{X})$ is determined by the equivalence classes of $F$-one edge free splittings that are adjacent to $\mathcal{S}'$. Therefore we have $f(\mathcal{S})=\mathcal{S}$ and the equivalence class of any one-edge free splitting is fixed by $f$. Theorem~\ref{theorem scott swarup} then implies that $f=\mathrm{id}$. This concludes the proof of Theorem~\ref{Rigidity FS}.
\hfill\qedsymbol

\bibliographystyle{alphanum}
\bibliography{bibliographie}

\noindent \begin{tabular}{l}
Laboratoire de mathématique d'Orsay\\
UMR 8628 CNRS \\
Université Paris-Saclay\\
91405 ORSAY Cedex, FRANCE\\
{\it e-mail: yassine.guerch@universite-paris-saclay.fr}
\end{tabular}
\end{document}

%% file: figurecheminLn.tex
\centering
\captionsetup{justification=centering}
\hspace{-0.8cm}
\begin{tikzpicture}[scale=1.4]
\draw (0:0) -- (90:1);
\draw (0:0) -- (150:1);
\draw (0:0) -- (210:1);
\draw (0:0) -- (270:1);
\draw (0:0) -- (330:1);
\draw (0:0) -- (30:1);
\draw (0:0) node {$\bullet$};
\draw (90:1) node {$\bullet$};
\draw (210:1) node {$\bullet$};
\draw (150:1) node {$\bullet$};
\draw (270:1) node {$\bullet$};
\draw (330:1) node {$\bullet$};
\draw (30:1) node {$\bullet$};
\draw (0:0) node[right, scale=0.75] {$\;\;\{0\}$};
\draw (90:1) node[above, scale=0.75] {$\left\langle x_2x_1x_2 \right\rangle$};
\draw (150:1) node[above left, scale=0.75] {$\left\langle x_2 \right\rangle$};
\draw (210:1) node[below left, scale=0.75] {$\left\langle x_3 \right\rangle$};
\draw (270:1) node[below, scale=0.75] {$\left\langle x_4 \right\rangle$};
\draw (330:1) node[below right, scale=0.75] {$\left\langle x_5 \right\rangle$};
\draw (30:1) node[above right, scale=0.75] {$\left\langle x_6 \right\rangle$};
\draw (270:1.5) node[scale=1] {$\mathcal{X}'$};

\draw (0:1.25) -- (0:2.25);

\draw (0,-3.5) -- (1.71/2,-3);
\draw (0,-3.5) -- (1.71/2,-4);
\draw (0,-3.5) -- (0,-4.5);
\draw (0,-3.5) -- (-1.71/2,-3);
\draw (0,-3.5) -- (-1.71/2,-4);

\draw (0,-3.5) node {$\bullet$};
\draw (1.71/2,-3) node {$\bullet$};
\draw (1.71/2,-4) node {$\bullet$};
\draw (0,-4.5) node {$\bullet$};
\draw (-1.71/2,-3) node {$\bullet$};
\draw (-1.71/2,-4) node {$\bullet$};

\draw (0,-3.5) node[right, scale=0.75] {$\;\;\left\langle x_6 \right\rangle$};
\draw (1.71/2,-3) node[above right, scale=0.75] {$\left\langle x_1 \right\rangle$};
\draw (1.71/2,-4) node[below right, scale=0.75] {$\left\langle x_5 \right\rangle$};
\draw (0,-4.5) node[below, scale=0.75] {$\left\langle x_4 \right\rangle$};

\draw (-1.71/2,-3) node[above left, scale=0.75] {$\left\langle x_2 \right\rangle$};
\draw (-1.71/2,-4) node[below left, scale=0.75] {$\left\langle x_3 \right\rangle$};
\draw (1.25, -3.5) -- (2.25,-3.5);
\draw (0,-5) node[scale=1] {$\mathcal{Y}_6$};
\end{tikzpicture}
\hspace{-0.5cm}
\begin{tikzpicture}[scale=1.4]
\draw (0,-3.5) -- (0,-2.5);
\draw (0,-3.5) -- (1.71/2,-3);
\draw (0,-3.5) -- (1.71/2,-4);
\draw (0,-3.5) -- (0,-4.5);
\draw (0,-3.5) -- (-1.71/2,-3);
\draw (0,-3.5) -- (-1.71/2,-4);

\draw (0,-2.5) node {$\bullet$};
\draw (0,-3.5) node {$\bullet$};
\draw (1.71/2,-3) node {$\bullet$};
\draw (1.71/2,-4) node {$\bullet$};
\draw (0,-4.5) node {$\bullet$};
\draw (-1.71/2,-3) node {$\bullet$};
\draw (-1.71/2,-4) node {$\bullet$};

\draw (0,-3.5) node[right, scale=0.75] {$\;\;\{0\}$};
\draw (0,-2.5) node[above, scale=0.75] {$\;\;\left\langle x_1 \right\rangle$};
\draw (1.71/2,-3) node[above right, scale=0.75] {$\left\langle x_6 \right\rangle$};
\draw (1.71/2,-4) node[below right, scale=0.75] {$\left\langle x_5 \right\rangle$};
\draw (0,-4.5) node[below, scale=0.75] {$\left\langle x_6x_4x_6 \right\rangle$};

\draw (-1.71/2,-3) node[above left, scale=0.75] {$\left\langle x_2 \right\rangle$};
\draw (-1.71/2,-4) node[below left, scale=0.75] {$\left\langle x_6x_3x_6 \right\rangle$};
\draw (0,-5) node[scale=1] {$\mathcal{X}_6^{(3)}$};

\draw (1.25, -3.5) -- (2.25,-3.5);

\draw (0:0) -- (270:1);
\draw (0:0) -- (150:1);
\draw (0:0) -- (210:1);
\draw (0:0) -- (330:1);
\draw (0:0) -- (30:1);
\draw (0:0) node {$\bullet$};
\draw (270:1) node {$\bullet$};
\draw (210:1) node {$\bullet$};
\draw (150:1) node {$\bullet$};
\draw (330:1) node {$\bullet$};
\draw (30:1) node {$\bullet$};
\draw (0:0) node[right, scale=0.75] {$\;\;\left\langle x_6 \right\rangle$};
\draw (270:1) node[below, scale=0.75] {$\left\langle x_4 \right\rangle$};
\draw (150:1) node[above left, scale=0.75] {$\left\langle x_2 \right\rangle$};
\draw (210:1) node[below left, scale=0.75] {$\left\langle x_3 \right\rangle$};

\draw (330:1) node[below right, scale=0.75] {$\left\langle x_5 \right\rangle$};
\draw (30:1) node[above right, scale=0.75] {$\left\langle x_2x_1x_2 \right\rangle$};
\draw (270:1.5) node[scale=1] {$\mathcal{Y}_6'$};
\draw (0:1.25) -- (0:2.25);
\end{tikzpicture}
\begin{tikzpicture}[scale=1.4]
\draw (0:0) -- (90:1);
\draw (0:0) -- (150:1);
\draw (0:0) -- (210:1);
\draw (0:0) -- (270:1);
\draw (0:0) -- (330:1);
\draw (0:0) -- (30:1);
\draw (0:0) node {$\bullet$};
\draw (90:1) node {$\bullet$};
\draw (210:1) node {$\bullet$};
\draw (150:1) node {$\bullet$};
\draw (270:1) node {$\bullet$};
\draw (330:1) node {$\bullet$};
\draw (30:1) node {$\bullet$};
\draw (0:0) node[right, scale=0.75] {$\;\;\{0\}$};
\draw (90:1) node[above, scale=0.75] {$\left\langle x_2x_1x_2 \right\rangle$};
\draw (150:1) node[above left, scale=0.75] {$\left\langle x_2 \right\rangle$};
\draw (210:1) node[below left, scale=0.75] {$\left\langle x_6x_3x_6 \right\rangle$};
\draw (270:1) node[below, scale=0.75] {$\left\langle x_6x_4x_6 \right\rangle$};
\draw (330:1) node[below right, scale=0.75] {$\left\langle x_5 \right\rangle$};
\draw (30:1) node[above right, scale=0.75] {$\left\langle x_6 \right\rangle$};
\draw (270:1.5) node[scale=1] {$\mathcal{X}_6^{(2)}$};

\draw (0,-1.75-0.25) -- (0,-2.75-0.25);

\draw (0,-3.5) -- (1.71/2,-3);
\draw (0,-3.5) -- (1.71/2,-4);
\draw (0,-3.5) -- (0,-4.5);
\draw (0,-3.5) -- (-1.71/2,-3);
\draw (0,-3.5) -- (-1.71/2,-4);

\draw (0,-3.5) node {$\bullet$};
\draw (1.71/2,-3) node {$\bullet$};
\draw (1.71/2,-4) node {$\bullet$};
\draw (0,-4.5) node {$\bullet$};
\draw (-1.71/2,-3) node {$\bullet$};
\draw (-1.71/2,-4) node {$\bullet$};

\draw (0,-3.5) node[right, scale=0.75] {$\;\;\left\langle x_2 \right\rangle$};
\draw (1.71/2,-3) node[above right, scale=0.75] {$\left\langle x_1 \right\rangle$};
\draw (1.71/2,-4) node[below right, scale=0.75] {$\left\langle x_6 \right\rangle$};
\draw (0,-4.5) node[below, scale=0.75] {$\left\langle x_5 \right\rangle$};

\draw (-1.71/2,-3) node[above left, scale=0.75] {$\left\langle x_6x_3x_6 \right\rangle$};
\draw (-1.71/2,-4) node[below left, scale=0.75] {$\left\langle x_6x_4x_6 \right\rangle$};
\draw (0,-5) node[scale=1] {$\mathcal{Y}_6^{(2)}$};
\end{tikzpicture}